\documentclass[review,hidelinks]{elsarticle}
\pdfoutput=1 
\usepackage[displaymath,mathlines,pagewise]{lineno}
\usepackage{amsmath}
\usepackage{amssymb} 
\usepackage{booktabs}
\usepackage{multirow}

\newcommand{\pd}[2]{\frac{\partial #1}{\partial #2}}

\newcommand{\pdl}[2]{\partial #1 / \partial #2}
\newcommand{\pdls}[3]{\partial^2 #1 / \partial #2 \partial #3}

\newcommand{\Sec}[1]{Sec.\ \ref{#1}}

\newcommand{\Eq}[1]{Eq.\ (\ref{#1})}
\newcommand{\Eqs}[2]{Eqs.\ (\ref{#1}) and (\ref{#2})}
\newcommand{\Fig}[1]{Fig.~\ref{#1}}
\newcommand{\Figs}[2]{Figs.\ \ref{#1} -- \ref{#2}}

\newcommand{\Tab}[1]{Tab.\ \ref{#1}}

\newcommand{\lr}{\left(}
\newcommand{\rr}{\right)}

\newcommand{\ls}{\left[}
\newcommand{\rs}{\right]}
\newcommand{\lc}{\left\{}
\newcommand{\rc}{\right\}}

\newcommand{\be}{\begin{equation}}
\newcommand{\blnm}{\begin{linenomath*}}
\newcommand{\ee}{\end{equation}}
\newcommand{\elnm}{\end{linenomath*}}

\newcommand{\bspl}{\begin{split}}
\newcommand{\espl}{\end{split}}
\newcommand{\bea}{\begin{eqnarray}}
\newcommand{\eea}{\end{eqnarray}}
\newcommand{\bagn}{\begin{align}}
\newcommand{\eagn}{\end{align}}
\newcommand{\bs}{\begin{split}}
\newcommand{\es}{\end{split}}
\newcommand{\bc}{\begin{center}}
\newcommand{\ec}{\end{center}}
\newcommand{\bp}{\begin{picture}(0,0)}
\newcommand{\ep}{\end{picture}}
\newcommand{\bfl}{\begin{flushleft}}
\newcommand{\efl}{\end{flushleft}}

\newcommand{\bx}{\mathbf{x}}

\newcommand{\bu}{\mathbf{u}}

\newcommand{\bn}{\mathbf{n}}
\newcommand{\br}{\mathbf{r}}
\newcommand{\bng}{\mathbf{n}_\Gamma}

\newcommand{\eph}{\epsilon_h}
\newcommand{\alpx}{\alpha \lr  \bx,t  \rr}
\newcommand{\alpp}{\alpha \lr  \psi_0 \lr \bx,t \rr \rr}
\newcommand{\alpps}{\alpha \lr  \psi_0 \rr}
\newcommand{\deltaa}{\delta \lr  \alpha \rr}

\newcommand{\psio}{\psi_0  \left( \alpha \right)}
\newcommand{\psix}{\psi_0  \left( \bx,t \right)}
\newcommand{\psizw}{\psi_0'}
\newcommand{\psip}{\psi_0' \lr \psi_1 \rr}
\newcommand{\psipf}{\psi_0' \lr \psigs \rr}

\newcommand{\psig}{\psi_1 \left( \alpha,\gamma \right)}
\newcommand{\psigf}{\psi_1 \left( \alpha,\gamma_1 \right)}
\newcommand{\psigs}{\psi_1 \left( \alpha,\gamma_2 \right)}

\journal{Journal of Computational Physics}

\begin{document}

\begin{frontmatter}

\title{A consistent solution of the re-initialization equation\\ in the conservative level-set method}

\author{Tomasz Wac{\l}awczyk\fnref{}\corref{mycorrespondingauthor}}
\address{Institute of Numerical Methods in Mechanical Engineering,\\ Technische Universit{\"a}t Darmstadt,\\ Dolivostr. 15, 64293 Darmstadt, Germany}


\cortext[mycorrespondingauthor]{Tomasz Wac{\l}awczyk}
\ead{twacl@fnb.tu-darmstadt.de}


\begin{abstract}
In this paper, a new 
re-initialization method 
for the  conservative level-set function
is put forward.
First, it has been shown 
that the re-initialization and advection equations 
of the conservative level-set function 
are  mathematically equivalent to 
the re-initialization and advection equations 
of the localized signed distance function.
Next, a new discretization 
for the spatial derivatives 
of the conservative level-set function 
has been proposed.
This new discretization  
is consistent with 
the re-initialization procedure 
and it guarantees 
a second-order convergence rate of 
the interface curvature 
on gradually refined grids.
The new re-initialization method
does not introduce artificial deformations 
to stationary and non-stationary interfaces,
even when the number of re-initialization 
steps is large.
\end{abstract}


\end{frontmatter}

\section{Introduction}
\label{sec1}
%
In the conservative level-set method
introduced by Olsson and Kreiss  \cite{olsson05},  
a solution of the transport equation
of the characteristic function $\alpx$ 
is divided into two steps: 
advection and re-initialization.
These two steps are
carried out consecutively 
in only one single time step $\Delta t$.
The solution of 
the $\alpx$ advection equation
is typically performed using 
an explicit or implicit time discretization
with a TVD MUSCL flux limiter,
in order to keep
$0 \le \alpx \le 1$ 
a bounded function
\cite{olsson05, shukla10, so11, tiwari13, balcazar2014}.
The present work 
provides improvements 
in the second step of 
the above interface capturing procedure. 
Let us consider 
the re-initialization equation
 \blnm
 \be
   \pd{\alpha}{\tau} = \nabla \cdot \ls D |\nabla\alpha| \bng - C \alpha \lr 1 - \alpha \rr \bng \rs,
   \label{eq1}
 \ee
 \elnm
where $\alpx$ denotes 
a regularized Heaviside function, 
$\tau$ is the artificial time, 
$\bng=\nabla \alpha/|\nabla \alpha|$
is a vector normal to 
the iso-lines 
(iso-surfaces) of $\alpx$. 
We assume that $C=1\,m/s$, $D=\Delta x/2 C=\eph C$ 
where $\Delta x=\Delta y=\Delta z=1/N_c$ 
is the size a control volume and $N_c$ 
is the total number of control volumes 
in the  $x,\,y$ or $z$ direction.

In \cite{olsson05} 
it has been shown   
the solution of \Eq{eq1} 
after advection of $\alpx$,
reduces artificial 
deformations of the interface 
induced by numerical errors 
during the advection step.
Therein it was also
noticed that the analytical solution 
to the stationary equation (\ref{eq1}) 
is given by 
 \blnm
 \be
  \alpp = 1 - \frac{1}{1+\exp{\lr \psix/\eph \rr}}=
  \frac{1}{2} \ls 1+\tanh{\lr \frac{\psix}{2\eph} \rr} \rs,
  \label{eq2}
 \ee
 \elnm
where $\psix$ 
is the signed distance function 
\cite{balcazar2014,glasner2001,sun2007}.
When $\eph \to 0$, 
then $\alpps$ in equation (\ref{eq2}) 
tends to the phase indicator function 
which is given by 
the exact Heaviside function 
$H \lr \psi_0 \rr $ as is
demonstrated in \ref{appA}. 
The phase indicator function $H \lr \psi_0 \rr $
is typically discretized in 
the volume of fluid (VOF) 
family of methods  
satisfying the 
law of conservation  
of mass 
\cite{trygg11}.
In real simulations $\eph \ne 0$, 
and hence $\alpps$ in equation (\ref{eq2})
is not represented
by a sharp jump  
localized at 
the interface $\Gamma$.
The level-set function $\alpps$ 
is a Lipschitz continuous function
and, therefore, resembles 
the signed distance function
in the standard level-set (LS) method 
introduced by Osher and Sethian 
\cite{osher1988} and extended by
Sussman et al. \cite{sussman94,sussman98}.
For these reasons, Olsson and Kreiss  \cite{olsson05}
called their interface capturing technique 
the conservative level-set (CLS) method. 
The signed distance function
derived from \Eq{eq2} is given 
by the equation
 \blnm
 \be
  \psi_0 \lr \alpha \rr = \eph \ln{\ls \frac{\alpha\lr \psi_0 \rr}
                                  {1 - \alpha \lr  \psi_0 \rr} \rs}.
 \label{eq3}
 \ee
 \elnm
We note that in \Eq{eq2} 
the interface $\Gamma$ 
is localized at $\alpha \lr \bx_\Gamma,t \rr=1/2$
whereas in \Eq{eq3} the set of points
where the signed distance function 
$\psi_0 \lr \alpha \lr \bx_\Gamma,t \rr = 1/2\rr=0$ 
represents a position 
of the interface $\Gamma$.

A mapping 
between $\alpps$ and $\psio$ in 
equations (\ref{eq2}) and (\ref{eq3})
suggests 
a closer relation 
between the CLS method and 
the standard LS method 
exists.
In the present paper,
we show how this relation 
can be established. 
Moreover, 
we put forward
a new method 
for computation 
of higher-order spatial 
derivatives of $\alpps$,
which is  consistent with
the new re-initialization procedure.
The spatial derivatives of $\alpps$ 
obtained with our new method  
are later used to  approximate  
the interface curvature $\kappa$
with second-order accuracy.

A relation between 
the regularized Heaviside function 
and the signed distance function 
was first observed by Glasner \cite{glasner2001}, 
and was  used 
for a non-linear preconditioning
of the phase-field equation.
The non-linear preconditioning 
was later exploited by Sun and Beckermann \cite{sun2007}
in order to solve the phase-field equation 
in a context of the interface capturing.
Therein, it was mentioned 
that the stationary solutions 
to the phase-field equation 
and to equation (\ref{eq1}) 
are different.
Later in this paper,
the key differences between
the present results and
the results reported 
in \cite{sun2007} 
are addressed.

The main difficulty 
in using $\alpps$ and $\psio$
interchangeably is 
the lack of a correct numerical solution 
to the re-initialization equation (\ref{eq1}).
Although vast literature 
concerning the numerical solution 
of the re-initialization equation 
for the signed distance function exists, 
see \cite{osher03,hartmann2010,rocca14} 
to name only a few recent works,
the solution of equation (\ref{eq1}) 
has drawn less attention. 
In \cite{olsson05} 
the re-initialization equation (\ref{eq1})
is solved directly, 
in \cite{olsson07} 
to reduce artificial interface deformations 
due to discretization errors;
the diffusive term in \Eq{eq1} was projected on $\bng$,
however as it is shown by Shukla et al. \cite{shukla10}, 
this reformulation leads to numerical instabilities.
Recently,  
McCaslin and Desjardins \cite{mccaslin2014}
proposed to multiply 
diffusive and compressive fluxes
on the right-hand side (RHS) of \Eq{eq1} 
by a new function $\beta \lr \bx,t \rr$.
This allows them to vary
the number of steps
in their re-initialization
procedure depending on 
the local flow conditions, 
thus reducing artificial
deformations of the interface.

Shukla et al. in \cite{shukla10} 
assume that 
\Eq{eq1} 
has no physical meaning,
and thus 
it can be solved in 
the non-conservative form
without the term containing
the interface curvature 
$\kappa=-\nabla \cdot \bng$.
With such an assumption, 
the counteracting diffusive 
and compressive 
fluxes in \Eq{eq1} 
are projected 
only on the direction normal 
to the interface $\bng$
 \blnm
 \be
   \pd{\alpha}{\tau} = \bng \cdot \nabla  \ls \eph |\nabla \alpha | -  \alpha \lr 1 - \alpha \rr  \rs.
   \label{eq4}
 \ee
 \elnm
Moreover,
in \cite{shukla10} 
it has been shown
the key element which 
guarantees the successful 
numerical solution of \Eq{eq4} is 
the discretization of $|\nabla \alpha|=|\nabla\psi| F\lr \alpha,\gamma \rr$,
where $\psi$ is a mapping function 
which smooths $\alpx$ and 
allows to compute $|\nabla \alpha|$ with 
a smaller error.

Since in 
\cite{olsson05,shukla10,tiwari13,
            balcazar2014,olsson07,mccaslin2014}
the discretization and solution of \Eq{eq1} 
are only briefly addressed,
in this paper we mainly focus 
on the discretization and 
solution of \Eq{eq1} in 
the framework 
of the second-order
accurate finite 
volume method.
In particular, we are interested
in the case where the number of 
re-initialization steps in 
the numerical solution of \Eq{eq1} 
is   $N_\tau \gg 1$
and the interface $\Gamma$ 
is stationary.
As it is described by McCaslin and Desjardins \cite{mccaslin2014},
in such circumstances $\Gamma$ is especially 
prone to artificial deformations caused 
by errors in calculations
of $\alpps$ and $\bng$.

The outline of this paper
is as follows.
In Section \ref{sec2},
the influence 
of the mapping function 
$\psi \lr \alpha,\gamma \rr$ 
introduced in \cite{shukla10}
on a convergence rate
of numerical solutions 
to \Eq{eq1}
is analyzed.
In Section \ref{sec3}, we show
that under certain conditions
the mapping function $\psi \lr \alpha,\gamma \rr$  
approximates the signed distance function
$\psio$ up to higher-order terms.
For this reason, 
in the present work,
it is proposed
to use the signed distance 
function $\psio$
as the new mapping function
in the discretization
of $|\nabla \alpha|$
in \Eq{eq1}.
Consequently, 
Section \ref{sec4}
presents the selection
of the mapping function $\psi$
further leads to a new, 
mathematically consistent
method for computation of
spatial derivatives of $\alpps$,
and thus,
the interface curvature.
This allows us to reformulate
the re-initialization 
and advection equations
of the conservative level-set function $\alpps$ 
in Section \ref{sec5}.
Moreover, in Section \ref{sec5}, 
we show mathematical equivalence 
between the CLS method
where the interface $\Gamma$ 
is represented by $\alpps$ 
and the standard LS method
where the interface is 
represented by $\psio$,
localized at
the interface  
by Dirac's delta.
In Section \ref{sec6}, 
we investigate
properties of the newly formulated
re-initialization and advection equations,
in particular their convergence rates and errors
in approximation of spatial derivatives
of $\alpps$ used to compute the interface
curvature.
Finally, 
the new re-initialization
method is examined 
in several  
test cases with advection
in order to closely inspect 
its conservative properties.
%

\section{Selection of the mapping function}
\label{sec2}
%
%
To assure convergence 
of equation (\ref{eq4}) during integration
in time $\tau$, in \cite{shukla10}
the mapping function 
$\psi \lr \alpha,\gamma \rr$ 
that smooths $\alpx$ 
was introduced
for discretization of
$|\nabla \alpha|$. 
Therein it was also
noticed that $\psi \lr \alpha,\gamma \rr$
has to satisfy two conditions,
the first condition 
is given by the equality
\blnm
\be
 \frac{\nabla \alpha}{|\nabla \alpha|} =
 \frac{\nabla \psi}{|\nabla \psi|},
 \label{eq5}
\ee
\elnm
as the mapping cannot change 
directions of the vectors 
normal to $\alpx$ iso-surfaces.
The second condition 
demands that 
the linear relation
between $\nabla \psi$ 
and  $\nabla \alpha$ exists
\blnm
\be
 \nabla \alpha =  F\lr \alpha, \gamma \rr \nabla \psi,
 \label{eq6}
\ee
\elnm
where $F \lr \alpha,\gamma \rr$ 
is a known function
and $0<\gamma<1$ is a constant.
In this work, 
we show that 
the condition given
by \Eq{eq6} can be also 
used to compute
the second-order spatial derivatives
of the conservative level-set function 
$\alpha \lr \psi \lr \bx,t \rr \rr$
\blnm
 \be
  \alpha_{,ij} 
   =  F\lr \alpha, \gamma \rr 
      \ls \psi_{,ij}  
   +      \psi_{,i}\psi_{,j} \pd{F}{\alpha} \rs,
 \label{eq7}
 \ee
\elnm
when the mapping 
function $\psi$ has been chosen 
properly,  
$\alpha_{,ij}=\pdls{\alpha}{x_i}{x_j}$
and $i,j=1,2,3$. 
Equation (\ref{eq7}) is required 
for a consistent approximation 
of the interface 
curvature $\kappa$.

\begin{figure}[h!] \nonumber
 \includegraphics[angle=-90]{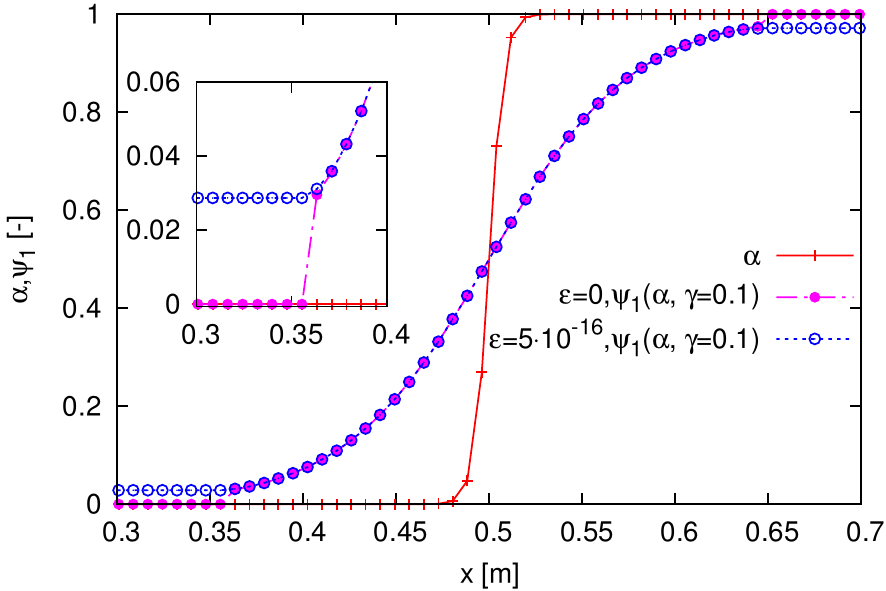}
 \includegraphics[angle=-90]{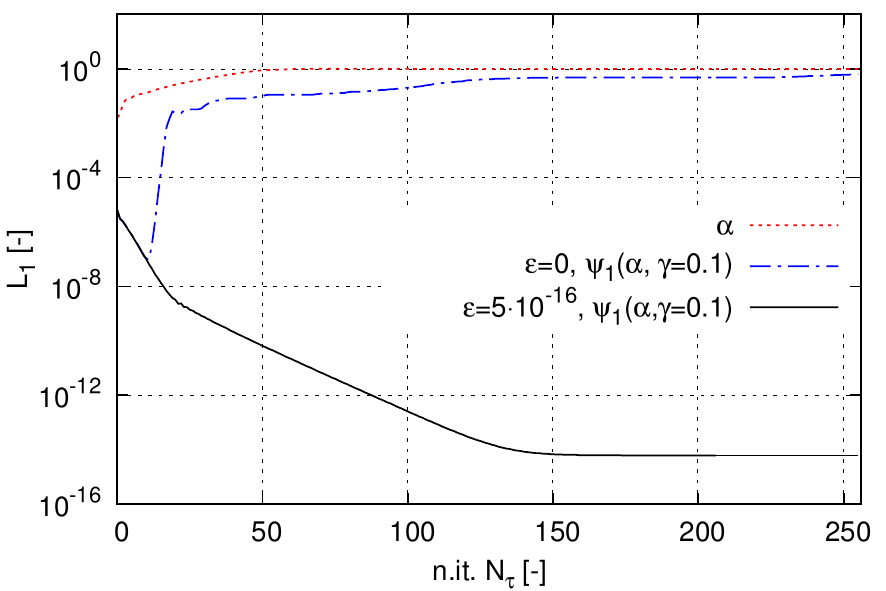}
 \caption{\small{
          Re-initialization 
          of the 1D regularized Heaviside function
          $H\lr x-0.5 \rr$:  
          (a) influence of $\epsilon \ne 0$ 
          in \Eq{eq8} on the presence of discontinuities  
          appearing due to arithmetical underflow, 
          (b) convergence of the solution to \Eq{eq1} 
          with different discretizations of $|\nabla\alpha|$:  
          $\alpha$ without smoothing, 
          $\psig$ where $\epsilon =0$ and
          $\psig$ where $\epsilon =5 \cdot 10^{-16}$,
          $\gamma=0.1$ in all cases.   
          $L_1$ norm is defined by \Eq{eq9}.}}
 \label{fig1}
 \end{figure}

Unlike in
works \cite{shukla10,tiwari13}, 
in this paper 
we use the mapping functions $\psi$
for discretization of $|\nabla \alpha|$ 
in equation (\ref{eq1}).
First, we note that the definition 
of the mapping function $\psi_1 \lr \alpha,\gamma \rr$ 
from \cite{shukla10}
\blnm
 \be
   \psig  = 
   \frac{\lr \alpha+\epsilon \rr^\gamma}
        {\lr \alpha+\epsilon \rr^\gamma+\lr 1-\alpha+\epsilon \rr^\gamma},
   \label{eq8}
 \ee
\elnm
where originally $\epsilon=0$, 
introduces discontinuities 
in the initial condition 
to \Eq{eq1} as it
is depicted in \Fig{fig1}(a).
Two discontinuities are 
caused by the arithmetic 
underflow when 
$\alpha \to 0$ and 
$1-\alpha \to 0$.
In order to avoid this, 
a straightforward modification 
of the mapping function from \cite{shukla10} 
is introduced. 
In figure \ref{fig1}(a)  
we show that setting
$\epsilon=5\cdot 10^{-16}$ 
allows avoiding jumps 
in $\psi_1$.
Since the arithmetic 
or floating point underflow 
is a purely numerical phenomenon,
we always set $\epsilon=0$ when
analytical operations using \Eq{eq8} 
are performed. 
 \blnm
 \elnm

This minor modification in \Eq{eq8})
has a great impact on convergence 
of the numerical solution 
to \Eq{eq1}.
In figure \ref{fig1}(b) convergence
of the solutions to
\Eq{eq1} 
in the case of 
re-initialization of 
the 1D regularized Heaviside 
function
is presented, in this figure 
the distance between solutions
on two different time levels
is measured by the first-order norm 
\blnm
\be
 L_1= \frac{1}{N_c} \sum_{l=1}^{N_{c}} |\alpha_l^{n+1}-\alpha_l^{n}|,
\label{eq9}
\ee
\elnm
where $N_c$ is the number of
control volumes and $n+1$ denotes
a new time level.
In this study, 
the initial condition  
to \Eq{eq1} is 
given by \Eq{eq2} 
and we use three different functions in 
the second-order central differencing
discretization of $|\nabla \alpha|$:
$\alpha$ alone without smoothing, 
the original mapping function from \cite{shukla10}
where $\epsilon=0$, and the modified mapping function 
given by \Eq{eq8} where $\epsilon \ne 0$.
As in \cite{shukla10,tiwari13} 
we use the value $\gamma=0.1$ in \Eq{eq8}.
Re-initialization of 
the 1D regularized Heaviside function 
is performed 
in the computational domain 
$\Omega = <0,1>\,m$ where
the interface $\Gamma$ 
is located at $x_\Gamma=0.5\,m$;
the mesh distribution 
is uniform $\Delta x=1/N_c$, 
$N_c=128$ 
is the number 
of control volumes.
At all boundaries 
of the computational 
domain $\Omega$,
the Neumann 
boundary condition
for $\alpps$ is used.
In our second-order accurate 
finite volume solver Fastest,
the third order TVD Runge-Kutta
method introduced in \cite{gottlieb98}
is used to integrate \Eq{eq1}
in the time $\tau$;
the time step size is set to
$\Delta \tau =D/C^2 = \eph$.
More details
concerning discretization
of \Eq{eq1} in the Fastest
flow solver can be found
in \ref{appB}. 

Since the initial condition
to \Eq{eq1} is given by \Eq{eq2}, 
we expect an immediate convergence
of its solution to numerical zero  
because equation (\ref{eq1}) 
is initialized with 
its own analytical solution.
However, in \Fig{fig1}(b) 
it is observed
that only the solution with 
$\epsilon \ne 0$ in \Eq{eq8}
allows convergence 
during all $N_{\tau}=256$
time steps. 
 \begin{figure}[h!] \nonumber
 \includegraphics[angle=-90]{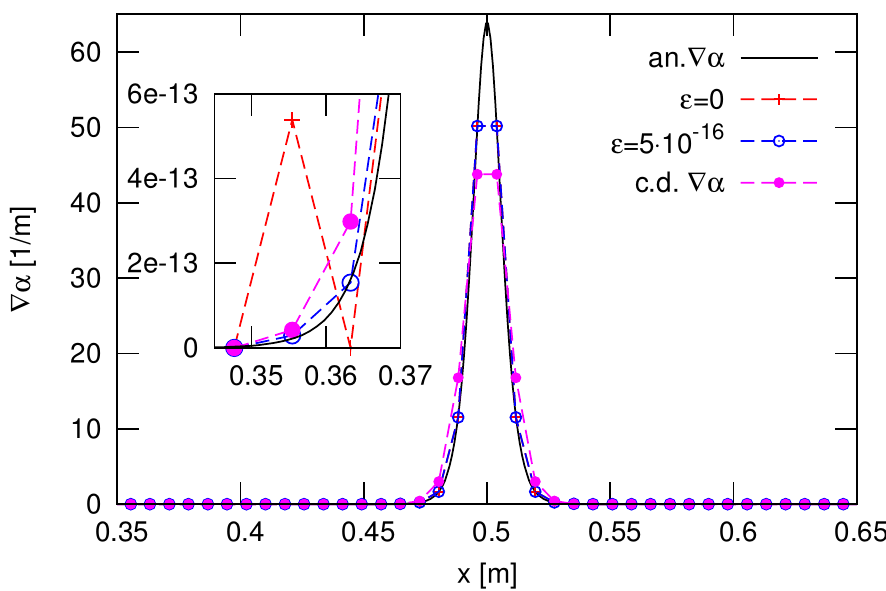}
 \caption{\small{The comparison of $|\nabla \alpha|$ after  $N_\tau=1$
                 re-initialization steps of the 1D regularized Heaviside function
                 with the central difference gradient approximation (c.d. $\nabla \alpha$) 
                 and analytical gradient (black solid line).
                 The $|\nabla \alpha|$ in equation (\ref{eq1}) is 
                 discretized using the mapping function $\psig$  with $\epsilon=0$ 
                 and  $\epsilon=5\cdot10^{-16}$, $\gamma=0.1$.}}
 \label{fig2}
 \end{figure}
 \begin{figure}[h!] \nonumber
 \includegraphics[angle=-90]{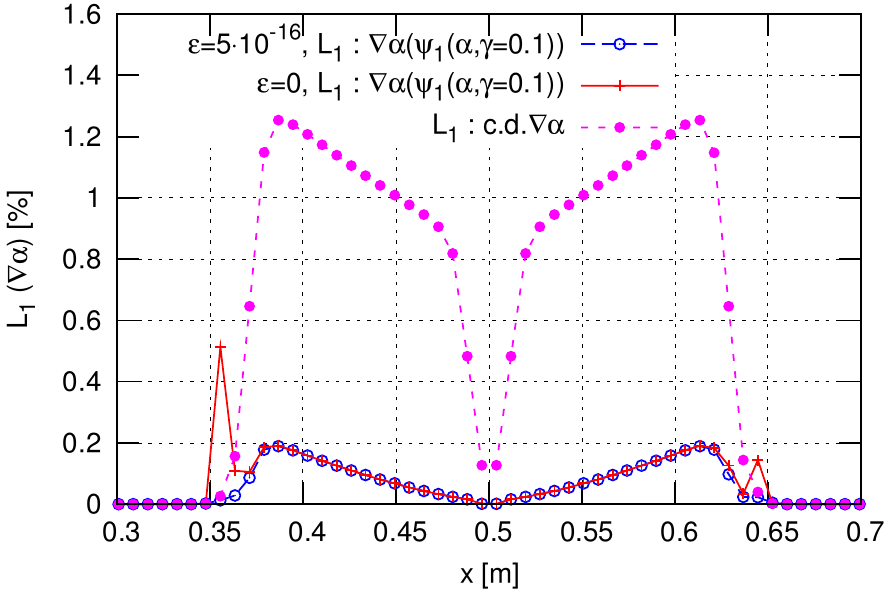}
 \caption{\small{Distributions of $L_1 \lr \nabla \alpha \rr$ norm 
                 defined by \Eq{eq10} after  $N_\tau=1$
                 re-initialization steps during 
                 re-initialization of the 1D regularized 
                 Heaviside function.  
                 The error $L_1$ is calculated 
                 for three $|\nabla \alpha|$
                 approximations depicted in \Fig{fig1}.}}
 \label{fig3}
 \end{figure}

To explain differences 
in the convergence rates 
which are observed in \Fig{fig1}(b), 
in figures \ref{fig2} and \ref{fig3} 
we compare the first-order
derivatives of $\alpps$ 
and $L_1 \lr \nabla \alpha \rr$ norms 
after $N_\tau=1$ re-initialization steps.
In figure \ref{fig3}, 
the $L_1 \lr \nabla \alpha \rr$ 
norm is defined by the equation
\blnm
\be
 L_1 \lr \phi \rr 
 = \frac{| \phi_{an} - \phi_{num}|}{|\phi_{an} |+\epsilon}, 
\label{eq10}
\ee
\elnm
where $\phi_{an}$, $\phi_{num}$ 
are functions calculated,
respectively,
analytically and numerically 
in each control volume,
$\epsilon=5\cdot10^{-16}$ and
$\phi= \nabla \alpha$.

In figures \ref{fig2} and \ref{fig3}, 
it is observed that 
both original $\epsilon=0$
and modified $\epsilon=5\cdot10^{-16}$ 
mapping functions provide very good 
approximations of $|\nabla \alpha|$ 
when compared with central difference
gradient approximation.
Differences between
these two gradient approximations
are visible only 
around $x=0.35\,m$ and $x=0.65\,m$
where the jumps in $\psi_1$
are present (compare 
results presented in figure 
\ref{fig1} and in
figures \ref{fig2}-\ref{fig3}).
If the mapping function 
is defined by \Eq{eq8}
with $\epsilon = 5 \cdot 10^{-16}$, 
the artificial oscillations 
are absent since $\psi_1$ 
is continuous everywhere.
The continuity of $\psi_1$  guarantees 
convergence of the solution to equation (\ref{eq1})
during all re-initialization steps
as shown in \Fig{fig1}(b).
From now on, 
when referring
to the mapping 
function $\psi_1$,
we reference its 
definition given by \Eq{eq8}
with $\epsilon=5\cdot 10^{-16}$.
%

\section{Relation between the mapping function and the signed distance function}
\label{sec3}
%
%
Figure \ref{fig1}(b)
shows the minor modification
of the mapping function $\psi_1$
guarantees convergence 
of the numerical solution to \Eq{eq1}
during all re-initialization steps.
Here, we also note that  
this discretization requires about
$N_\tau \approx 140$ time steps $\Delta \tau$ 
to achieve the stationary solution in spite
of the fact that \Eq{eq1} is initialized with
its  stationary
solution given by \Eq{eq2}. 
Therefore, the mapping function 
$\psi_1 \lr \alpha,\gamma=0.1 \rr$ where 
$\epsilon=5 \cdot 10^{-16}$ is not 
the best possible choice for 
discretization of $|\nabla \alpha|$  in \Eq{eq1}.
Hence, there is a need 
for an improved discretization 
of $|\nabla \alpha|$, 
guaranteeing immediate 
convergence  of the numerical 
solution to \Eq{eq1}  
towards a steady state.

In this section, we have shown 
when $0< N_c \gamma  \ll 1$
or equivalently $0< \gamma  \ll \Delta x$,
the mapping function $\psig$
given by \Eq{eq8} approximates 
the signed distance function $\psio$
up to higher-order terms.
Let us first note 
that the stationary solution (\ref{eq2}) 
can be rewritten as 
\blnm
\be
 \alpps =   \frac{1}{2} \ls 1+\tanh{\lr \frac{\psi_0}{2\eph} \rr} \rs
      =   \frac{\exp{\lr N_c \psi_0 \rr}}{\exp{\lr N_c \psi_0 \rr}+\exp{\lr-N_c \psi_0\rr}}.
\label{eq11}
\ee
\elnm
Next, we observe that
\blnm
\bea
\label{eq12}
 \alpha^\gamma = 
 \frac{\exp{\lr N_c \gamma \psi_0 \rr}}{\ls \exp{\lr N_c \psi_0 \rr}+\exp{\lr-N_c \psi_0\rr }\rs^\gamma},\\
 \lr 1-\alpha \rr^\gamma = 
 \frac{\exp{\lr -N_c \gamma \psi_0 \rr}}{\ls \exp{\lr N_c \psi_0 \rr}+\exp{\lr-N_c \psi_0\rr }\rs^\gamma},
\label{eq121}
\eea
\elnm
and we substitute 
\Eqs{eq12}{eq121}
into \Eq{eq8} 
to obtain
\blnm
\be
 \psi_1 =  \frac{\exp{\lr N_c \psi_0 \gamma \rr}}
                {\exp{\lr N_c \psi_0 \gamma \rr}+\exp{\lr-N_c \psi_0 \gamma \rr}}.
\label{eq13}
\ee
\elnm
Now, we expand 
exponents in \Eq{eq13}
in the Taylor series
\blnm
\bea
\label{eq14}
 \exp{\lr N_c \psi_0 \gamma \rr} = 1 + N_c \psi_0 \gamma 
                     + \frac{\lr N_c \psi_0 \gamma \rr^2}{2!}  
                     + \frac{\lr N_c \psi_0 \gamma \rr^3}{3!} + \dots \\
 \exp{\lr -N_c \psi_0 \gamma \rr} = 1 - N_c \psi_0 \gamma 
                  + \frac{\lr N_c \psi_0 \gamma \rr^2}{2!} 
                  - \frac{\lr N_c \psi_0 \gamma \rr^3}{3!} + \dots
\label{eq15}
\eea 
\elnm
If $0< N_c \gamma \ll 1$ 
is sufficiently small
then the higher-order 
terms in  
\Eqs{eq14}{eq15}
can be neglected
what gives
\blnm
\bea
\label{eq16}
 \exp{\lr  N_c \psi_0 \gamma \rr} = 1 + N_c \psi_0' \gamma,\\ 
 \exp{\lr -N_c \psi_0 \gamma \rr} = 1 - N_c \psi_0' \gamma,
\label{eq17}
\eea 
\elnm
where $\psi_0' \approx \psi_0$.
After substitution of \Eqs{eq16}{eq17} 
into \Eq{eq13} one obtains  
\blnm
\be
 \psi_1 = \frac{1+N_c \psi_0' \gamma}{2} 
      = \frac{1}{2} \lr 1 + \frac{\gamma}{2\eph} \psi_0' \rr.
\label{eq18}
\ee
\elnm
\Eq{eq18} is 
the exact relation between 
$\psi_1$ and $\psi_0'$
when $0 < N_c \gamma \ll 1$. 
Next, from \Eq{eq18}
we derive 
the relation 
between
the signed
distance function $\psi_0$
its approximation
$\psizw$
and the mapping 
function $\psig$
\blnm
 \be
 \psi_0 \approx  \psizw  = \frac{2\eph}{\gamma} \lr2 \psi_1 -1 \rr.
\label{eq19}
 \ee
\elnm
We note that 
the absolute 
value of the gradient 
of $\psi_0 \approx \psizw$
\blnm
 \be
 \Big|\pd{\psi_0}{x_i}\Big| \approx 
 \Big|\pd{\psizw}{x_i}\Big| = 
 \Big|\pd{\psizw}{\psi_1} \pd{\psi_1}{x_i}\Big| = 
 \Big|\frac{4\eph}{\gamma}\pd{\psi_1}{x_i}\Big| = 1
 \label{eq20}
 \ee
\elnm
since $|\nabla \psi_0|=1$ 
is the property  of 
the signed-distance function, see \cite{osher1988,osher03}. 
In 1D the second-order 
derivative of $\psi_0 \approx \psizw$
is equal to zero
\blnm
 \be
  \pd{}{x_1} \lr \pd{\psi_0 }{x_1} \rr \approx
  \pd{}{x_1} \lr \pd{\psizw }{x_1} \rr =
  \pd{}{x_1} \lr \frac{4\eph}{\gamma} \pd{\psi_1}{x_1} \rr = 0.
 \label{eq21}
 \ee
\elnm
The remaining question
to be considered
is how to select a value
of the constant $\gamma$
in the above equations;
the answer to 
this problem is given 
in Section \ref{sec5.1}.
Next, 
we show how to use 
the mapping functions
$\psi_0 \approx \psi_0'$
to compute  
spatial derivatives of
$\alpha \lr \psi_0 \rr$
and
$\alpha \lr \psi_0' \rr$, 
as this
leads to reformulation
of equation (\ref{eq1}).

\section{Computation of spatial derivatives with mapping functions}
\label{sec4}

In this section 
we derive 
formulas for 
the 
first and 
second-order
spatial derivatives 
of the conservative 
level-set function 
$\alpx$ ,
see \Eqs{eq6}{eq7}.
These new formulations  
exploit dependence of
the level-set function
$\alpha$ on the signed
distance function $\psio$
or its approximation $\psi_0' \lr \psi_1 \rr$,
see \Eq{eq3} or \Eq{eq19}, respectively.

Let us first calculate 
$\alpha_{,i}=\pdl{\alpha}{x_i}$,
$i=1,2,3$ using \Eq{eq8}, 
in this case the  
first-order 
spatial derivative 
is given by the equation 
\blnm
\be
\alpha_{,i} = \frac{\zeta^2 \delta^{1-\gamma}}{\gamma} \psi_{1,i} = F\lr \alpha,\gamma\rr  \psi_{1,i},
\label{eq22}
\ee
\elnm
where $\zeta \lr \alpha, \gamma \rr$ 
and $\delta \lr \alpha \rr$ 
are two auxiliary functions
\blnm
\bea
 && \zeta \lr \alpha,\gamma \rr = \alpha^\gamma 
                                + \lr 1-\alpha \rr^\gamma,
\label{eq23} \\
 && \delta \lr \alpha \rr =  \alpha \lr 1-\alpha \rr.
\label{eq24}
\eea
\elnm
The second-order spatial  
derivatives of $\alpha \lr \psi_0' \lr \psi_1 \rr \rr$
are obtained
directly from \Eq{eq22} 
and they read
\blnm
\begin{align}
  \bs
    \alpha_{,ij} = 
                 \frac{\delta^{1-\gamma} \zeta^2 }{\gamma} 
                & \lc \psi_{1,ij} 
               + \frac{\zeta}{\gamma} \psi_{1,i} \psi_{1,j} 
                 \right. \\
                & \left. 
                 \ls 2 \gamma \lr \lr 1-\alpha \rr^{1-\gamma} 
               - \alpha^{1-\gamma} \rr 
               + \lr 1-\gamma           \rr 
                 \lr 1-2\alpha          \rr 
                 \zeta \delta^{-\gamma} \rs \rc,
\label{eq25}
  \es
\end{align}
\elnm
where 
$\alpha_{,ij}=\pdls{\alpha}{x_i}{x_j}$
$i,j=1,2,3$.

Next,  
derivatives of $\alpp$  
in terms of the signed distance function 
$\psio$ are calculated, 
see \Eq{eq3}, which leads to
\blnm
\be
 \alpha_{,i} = \frac{ \delta \lr \alpha \rr  }{\eph} \psi_{0,i} 
\label{eq26}
\ee
\elnm
where $i=1,2,3$ and $\deltaa$
is defined by \Eq{eq24}. 
The second-order 
derivative  of $\alpp$
is calculated directly
from \Eq{eq26} as 
\blnm
\be
 \alpha_{,ij} = \frac{ \delta \lr \alpha \rr }{\eph} 
              \ls \psi_{0,ij} +  
              \frac{1}{\eph} \psi_{0,i} \psi_{0,j}
              \lr 1- 2 \alpha \rr \rs,
\label{eq27}
\ee
\elnm
where
$i,j=1,2,3$.

We observe that 
interesting similarities 
between \Eq{eq22} and \Eq{eq26} or
\Eq{eq25} and \Eq{eq27} do exist
when $0 < \gamma \ll \Delta x$.
For instance,
in the case of both
mapping functions,
the right-hand sides
of derived formulations
are multiplied by 
$\delta \lr \alpha \rr /\eph$,
and it is possible to show that
\blnm
 \be
 \frac{1}{\eph} \int_{-\infty}^{\infty} \delta \lr \alpha \rr d \psi_0 = 1.
 \label{ddelta}
 \ee
\elnm
Hence,
$\delta \lr \alpha \rr /\eph$
approximates 
Dirac's delta 
localized 
at the interface $\Gamma$
as it has 
a compact support 
(on the given grid 
 as shown in  
 figures \ref{fig2} and \ref{fig3}) 
and for $\eph \ne 0$
it is $C^\infty$.
Consequently, 
$\delta \lr \alpha \rr/ \eph$ 
restricts the support 
of $\alpha_{,i}$, $\alpha_{,ij}$
derivatives to the region localized
in the vicinity of the interface $\Gamma$.

Next, we note 
that for 
$0 < \gamma \ll \Delta x$
(see condition required 
to derive \Eq{eq19})
$\zeta \lr \alpha,\gamma \rr \approx 2$   
inside the support of $\deltaa/\eph$, see \Eq{eq23}.
With
this latter observation
and from a comparison between 
\Eq{eq22} and 
\Eq{eq26} the 
relation between $\gamma$ and $\eph$
is obtained 
\blnm
\be
 \frac{\gamma}{4} \approx \eph.
\label{eq28}
\ee
\elnm
\Eq{eq28} provides  
the condition of equality between 
the first and 
the second
order spatial derivatives
of $\alpps$ 
computed 
using the 
mapping function $\psig$
or $\psio$, 
see also \Eqs{eq20}{eq21}.

Finally,
we recognize the relation
between the present results and
the result from the theory of distributions,
see \cite{bronstein2012} page \textit{788}.
When $\eph \to 0$ then $\alpps \to H \lr \psi_0 \rr $
and $\deltaa/\eph \to \delta \lr \psi_0 \rr $,
in such case equation (\ref{eq22}) (for $0 < \gamma \ll \Delta x$) and 
equation (\ref{eq26}) read
\blnm
\be
  \nabla H \lr \psi_0 \rr = \mathcal{\delta} \lr \psi_0 \rr \nabla \psi_0,
\label{eq28b}
\ee
\elnm
where $H \lr \psi_0 \rr$ is the exact Heaviside function
and $\delta \lr \psi_0 \rr$ is the exact Dirac delta function,
both functions are localized at the interface 
$\Gamma \lr \psi_0 \approx \psi_0' =0 \rr$.
In particular, in the direction
normal to the interface $\bng$,
this definition holds only when 
$|\nabla \psi_0| \approx |\nabla \psi_0'| = 1$.
For this reason, $\psi_0 \approx \psi_0'$ 
has to be the signed distance function 
to satisfy the relation between 
$\nabla H \lr \psi_0 \rr$
and $\delta \lr \psi_0 \rr$
in equation (\ref{eq28b}).
%

\section{Reformulation of the re-initialization equation}
\label{sec5}

After substitution 
of  \Eq{eq26} into \Eq{eq1}
the new form of the re-initialization
equation is obtained
\blnm
\be
 \pd{\alpha}{\tau} = 
  \nabla \cdot \ls \delta \lr \alpha \rr 
  \lr |\nabla \psi_0| - 1 \rr \bng \rs,
\label{eq29}
\ee
\elnm
where $\delta \lr \alpha \rr$ 
is defined by \Eq{eq24} and 
$\bng = \nabla \psi_0 /|\nabla \psi_0|$,
see \Eq{eq5}.
We note that $\pdl{\alpha}{\tau} \equiv 0$ 
when $|\nabla \psi_0|=1$ or when 
$\delta \lr H \rr = 0$.
In the non-conservative form,
equation (\ref{eq29}) reads 
\blnm
 \begin{align}
 \bs
 \pd{\alpha}{\tau} =
  \bng \cdot \nabla \delta \lr \alpha \rr \lr |\nabla \psi_0 | - 1 \rr  \\
 +\bng \cdot \nabla \lr |\nabla \psi_0 | - 1 \rr \delta \lr \alpha \rr  \\
  - \kappa \lr |\nabla \psi_0 | - 1 \rr   \delta \lr \alpha \rr.
  \label{eq30}
 \es
 \end{align}
\elnm
Since  
$sign \ls \bng \cdot \nabla \delta \lr \alpha \rr\ \rs = - sign \ls \psi_0 \rs$,
the first term on the RHS of equation (\ref{eq30}) resembles 
the term in the re-initialization
equation of the signed distance
function introduced by Sussman et al. 
\cite{sussman94,sussman98}.
The second RHS term 
in equation (\ref{eq30})
contains information about 
the spatial distribution 
of a difference between 
$|\nabla \psio|$ 
and the solution 
to the eikonal equation 
$|\nabla \psi_0|=1$. 
The third RHS term in equation (\ref{eq30})
contains the interface
curvature $\kappa$ and expresses
its influence on the pair of re-initialized 
functions  $\alpps$ and $\psio$. 
Re-initialization
of $\alpps$
and $\psio$
in equations (\ref{eq29})-(\ref{eq30})
is  restricted to 
the region of support
of $\deltaa/\eph$, and thus
it is localized 
in the vicinity
of the interface 
$\Gamma$.
We emphasize 
that the solution 
to equation (\ref{eq1})
allowing re-initialization
of the level-set function $\alpps$
is mathematically equivalent
to the solution 
to equations (\ref{eq29})
or (\ref{eq30})
allowing re-initialization of 
$\psio$.

\subsection{Determination of the $\gamma$ value}
\label{sec5.1}

Shukla et al. in \cite{shukla10}
and Tiwari et al. in \cite{tiwari13} 
set the value of the constant $\gamma$
in equation (\ref{eq8})
to $\gamma=0.1$.
They argue
that this value is justified
because when $\gamma \to 0$
the mapping function given
by equation (\ref{eq8}) 
tends to $\psi_1 \to 1/2$
as shown in \ref{appA}.

In this section,
we investigate 
numerically
how to
select  
the value 
of the constant 
$0 < \gamma \ll \Delta x$.
We want to find 
$\gamma$ such that
$\psig$ in \Eq{eq19}
accurately approximates 
the signed distance function  $\psio$ 
and thus assures 
the solution to \Eq{eq1} 
with minimal error.
During tests in this section,
the number of grid nodes $N_c=128$ and
the support width 
$\eph = \Delta x / 2$ 
are both kept constant. 
 \begin{figure}[h!] \nonumber
 \includegraphics[angle=-90]{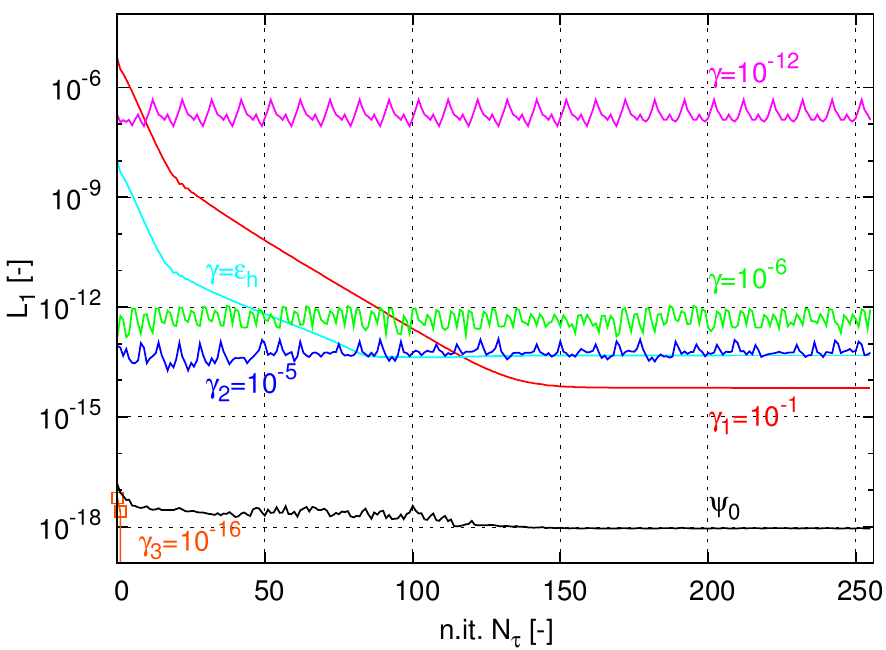}
 \caption{\small{The influence of $\gamma$ value in the mapping
                 function $\psig$ on convergence of the solution to \Eq{eq1} 
                 during re-initialization of the 1D regularized Heaviside function. 
                 $L_1$ norms are defined by \Eq{eq9}.}}
 \label{fig4}
 \end{figure}

In the beginning
of this study
we note that \Eq{eq28} 
could be considered as 
the condition 
which supplies
the maximal value 
of the constant 
$\gamma_{max} \approx 4 \eph$.
However,
to derive \Eq{eq19},
the more stringent condition 
$0 < \gamma \ll \Delta x$
is required;
convergence 
of \Eq{eq1} with 
$\gamma=\eph$ is
depicted in \Fig{fig4}.

Figure (\ref{fig4}) presents 
convergence of the solutions to \Eq{eq1}
in the case of re-initialization of the 1D regularized 
Heaviside function with different values of 
the constant $\gamma$.
Errors in solutions
to \Eq{eq1} are measured by 
the $L_1$ norm defined
by \Eq{eq9}.
The initial condition 
to \Eq{eq1} is given by \Eq{eq2} 
which is its stationary, 
analytical solution.
 \begin{figure}[h!] \nonumber
 \includegraphics[angle=-90]{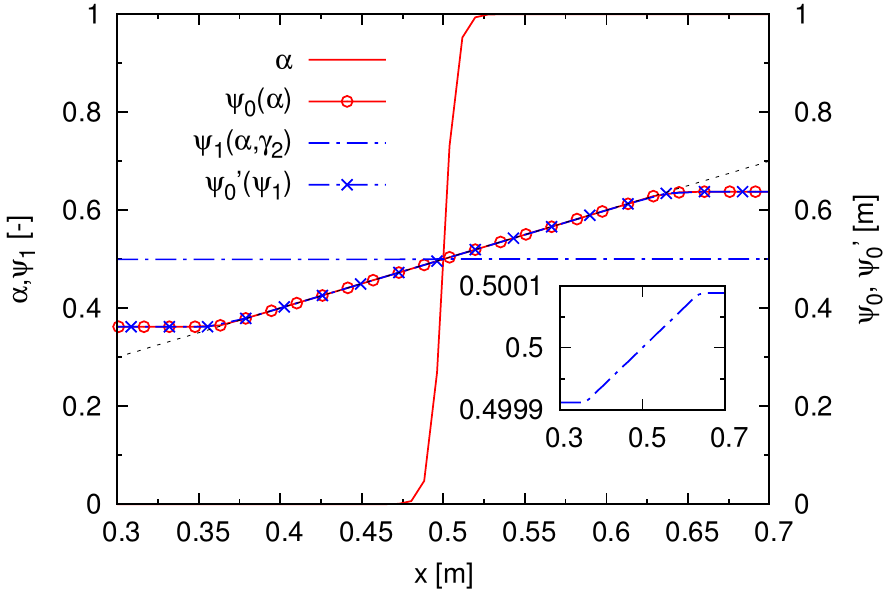}
 \includegraphics[angle=-90]{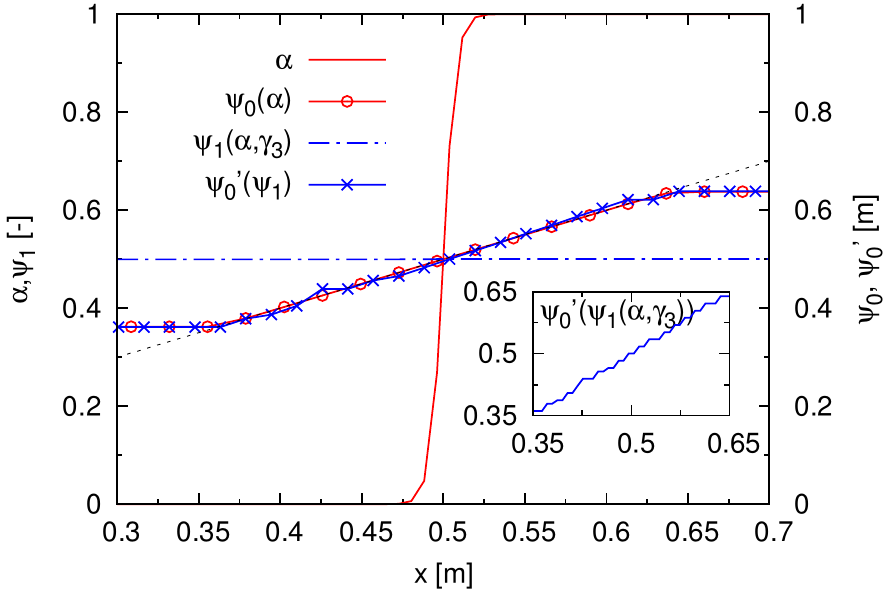}
 \caption{\small{The comparison of signed distance functions reconstructed
                 from $\psig$ with \Eq{eq19} after $N_\tau=256$ re-initialization steps
                 (a) $\gamma_2=10^{-5}$, (b) $\gamma_3=10^{-16}$
                 with $\psio$ obtained using \Eq{eq3}.
                 The lines with symbols are plotted every second or third point
                 to improve clarity in presentation of the results.}}
 \label{fig5}
 \end{figure}
$|\nabla \alpha|$
in \Eq{eq1} is discretized
with the mapping functions 
$\psig$ and $\psio$.
In order to obtain 
convergence when
$\gamma < 10^{-6}$ the 
size of time step 
$\Delta \tau=\eph$
was reduced 
to $\Delta \tau=\eph/2$.
In the present test case,
the interface $x_\Gamma=0.5\,m$ 
is localized exactly between
the two neighboring control volumes
$x_P,\, x_F$.

When $\gamma_1=0.1$,
one needs about 
$N_\tau \approx 140$ 
time steps
to achieve the 
stationary
solution. 
Similar to 
convergence with 
the exact signed distance function $\psio$, 
the most rapid convergence rate
and the smallest error 
is obtained when $\gamma_2=10^{-5}$.
For $\gamma_2 < 10^{-5}$,
the error of the solution to 
equation (\ref{eq1}) increases,
albeit remains on 
a constant level.
When $\gamma_3=10^{-16}$,
the solver needs only two iterations
and numerical zero 
is achieved
as shown in \Fig{fig4}.
We suspect 
the increment of the error level
is a numerical effect
associated with accuracy 
of the \textit{Fortran} compiler 
(all real variables 
in the Fastest solver 
are declared in double 
precision).

In figure \ref{fig5} we compare 
the signed distance functions 
reconstructed using \Eqs{eq3}{eq19} 
from the mapping functions $\psio$
and 
$\psigs$, $\psi_1 \lr \alpha,\gamma_3 \rr$
after $N_{\tau}=256$ re-initialization steps.
When $\gamma_3=10^{-16}$,
the reconstructed signed distance function 
$\psi_0' \lr \psi_1 \lr \alpha,\gamma_3 \rr\rr$ 
has a staircase shape 
due to the existence of 
numerical errors
in computation of high-order roots
in \Eq{eq8}, see \Fig{fig5}(b).
The latter result is 
empirical confirmation 
that the increment 
of $L_1$ norm levels,
observed in \Fig{fig4} when
$\gamma < \gamma_2=10^{-5}$,
has numerical origin. 
In figure \ref{fig5}(a), we  observe
that the signed distance function $\psizw$
reconstructed from $\psigs$ using equation
(\ref{eq19}) is identical with
$\psio$ obtained using equation (\ref{eq3}).

\subsection{Determination of the allowable interface width $\eph$}
\label{sec5.2}

In the previous sections 
we have chosen $\eph=\Delta x/2$
to set the support width of 
$\deltaa/\eph$,
as this value 
is also used in the literature 
\cite{olsson05,shukla10,tiwari13,balcazar2014,olsson07,mccaslin2014}.
However, at the beginning 
of Section \ref{sec5} 
it was mentioned
that other than
the signed distance function $\psio$, 
the Heaviside function $H \lr \psi_0 \rr$ 
is the stationary solution to \Eq{eq29} 
as well.
Subsequently, we will demonstrate
how this feature 
of \Eq{eq29} is preserved 
by the present numerical scheme 
which does not use
flux limiters or higher-order
essentially non-oscillatory 
 schemes
(see \ref{appB} for details).

The lack of flux limiters
in the present discretization
of the re-initialization equation 
and employment of the second-order flux 
limiters only during the advection step
may be an advantage, 
as the artificial deformations
of the interface may be avoided.
A deformation of the interface due to 
a \textit{minmod} flux limiter 
was observed in \cite{tiwari13}
during re-initialization 
of $\alpx$ with \Eq{eq4}.
The interface-grid lines alignment during advection 
is a known deficiency in compressive high-resolution 
schemes which use down-wind to maintain sharpness of the interface, 
and switch between higher and lower order differencing 
schemes to preserve its smoothness
\cite{ubbink99,waclawczyk06,waclawczyk08_2,waclawczyk08_3}.
Moreover, the numerical artifacts
described above cannot be accepted
during the reliable implementation 
of physical models based on
the variable relaxation velocity
$C \lr \bx,t\rr$, 
and the variable variance 
$\eph \lr \bx,t \rr$
in \Eq{eq1}.
A good example of 
the physical model 
requiring abovementioned
features of the numerical scheme, 
is the statistical model 
for the ensemble averaged description 
of interactions between the gas-liquid interface 
and turbulence
\cite{mwaclawczyk11,twaclawczyketal14,waclawczyk2015}. 

Since in \Eq{eq19} 
we have shown that $\psi_0 \approx \psi_0'$ 
for $0<\gamma \ll \Delta x$, 
in this section, for brevity,
we use only  $\psi_0$
for the discretization 
of $|\nabla \alpha|$
in \Eq{eq1}.
 \begin{figure}[h!] \nonumber
 \includegraphics[angle=-90]{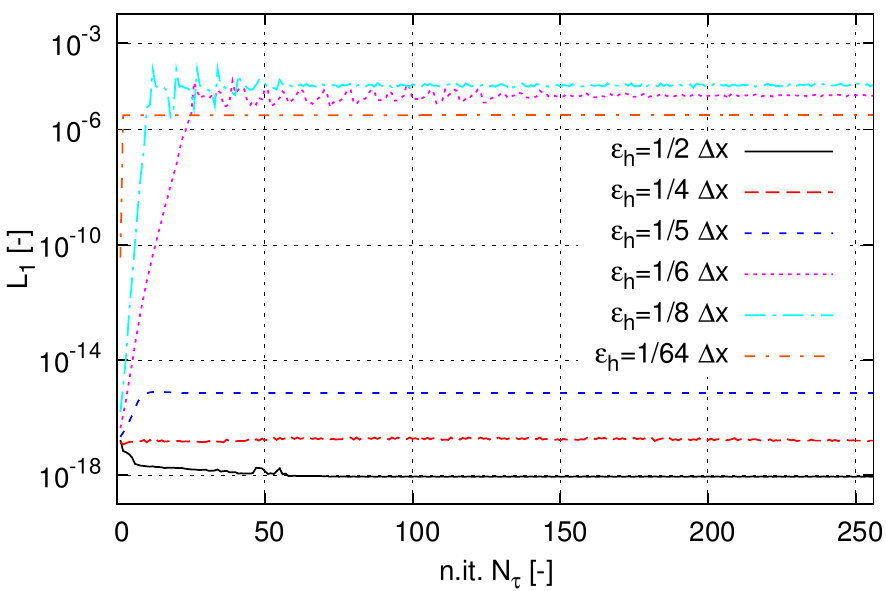}
 \caption{\small{The impact of the variable $\deltaa /\eph$
                 support width  $\eph = \Delta x /M$ where $M \le 16$, 
                 on convergence of the solution to the re-initialization 
                 equation (\ref{eq29}). Error is
                 measured by $L_1$ norm defined by \Eq{eq9}.}}
 \label{fig6}
 \end{figure}
In order to investigate
the influence of $\eph < \Delta x/2$ 
on convergence to the re-initialization equation,
we carry out the same test 
as in the previous section 
with $N_c=128$ 
and with the initial condition 
given by \Eq{eq2}. 
In the present case, 
the  position of the interface 
is set to $x_\Gamma=0.6\,m$,
and thus $\Gamma$ is not localized
exactly  between two neighboring
grid nodes $x_P,\,x_F$.  
Additionally,
$\eph=\Delta x/M$ where $M \ge 2$ is
an arbitrary integer number, 
the time step size 
is set to $\Delta \tau=\eph/2$.

In \Fig{fig6}, 
convergence of the 
solution to \Eq{eq29} 
with a varying width of the
interface $\eph=\Delta x/M$ 
can be observed.
The $L_1$ norm defined by \Eq{eq9} 
remains approximately 
constant for $M=2,4,5$ and 
increases rapidly for $M=6,8,64$.
This increment, however, does not lead to  
the divergence of the present numerical 
solutions, we observe
such behavior also  
when $M=16,32$ 
($L_1^M < 10^{-4}$).
 \begin{figure}[h!] \nonumber
 \includegraphics[angle=-90]{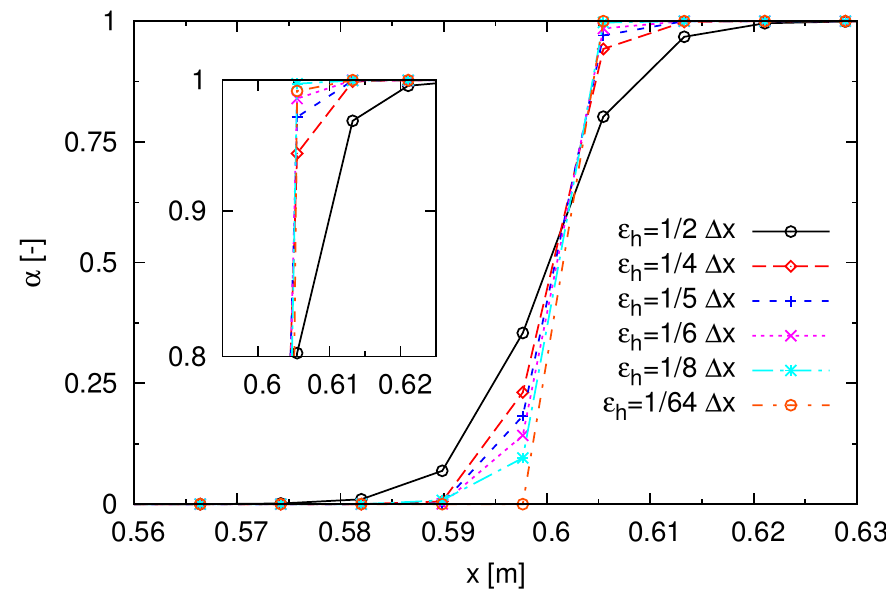}
 \includegraphics[angle=-90]{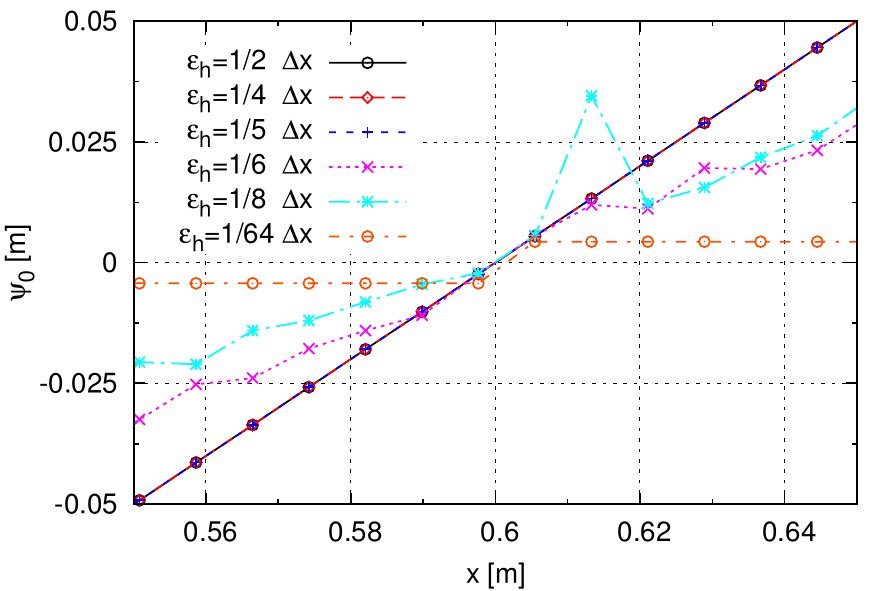}
 \caption{\small{The influence of 
                 the variable 
                 $\deltaa/\eph$ support width:
                 $\eph = \Delta x /M$ where $M \le 16$ 
                 on convergence to the solution to \Eq{eq29},
                 (a) $\alpps$, (b) $\psio$ 
                 after $N_\tau=256$ re-initialization steps,
                 the interface $\Gamma$ is
                 localized at $x_\Gamma=0.6\,m$.}}
 \label{fig7}
 \end{figure}
The observed increment 
in the $L_1$ norm
magnitude is related 
to the finite resolution
of the computational grid 
and the selected time step
size $\Delta \tau$.
 \begin{figure}[h!] \nonumber
 \includegraphics[angle=-90]{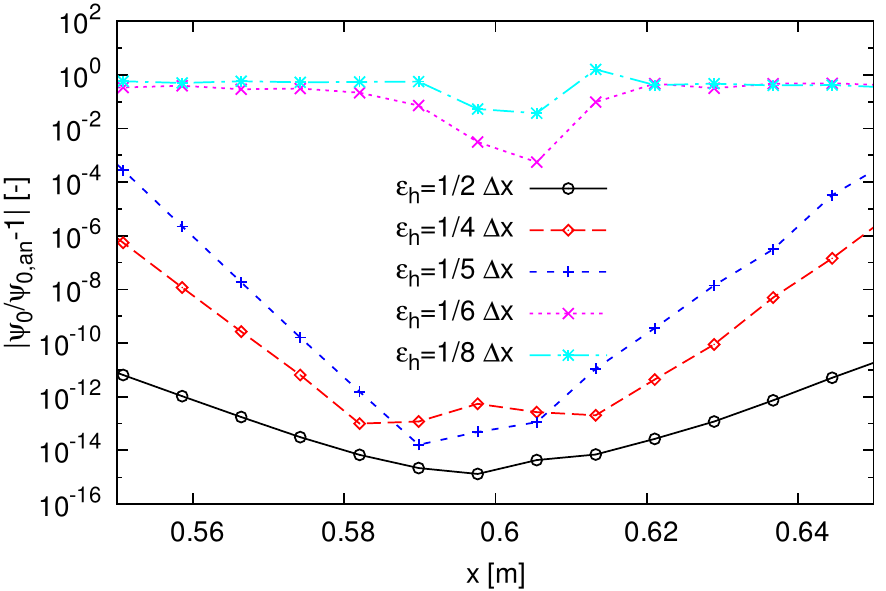}
 \caption{\small{The impact of the variable  
                 $\deltaa/\eph$ support width
                 $\eph = \Delta x /M$ where $M \le 8$ 
                 on the distribution of the error $|\psi_0/\psi_{0,an}-1|$, where
                 $\psi_{0,an}=x-x_\Gamma$ and $x_\Gamma = 0.6\,m$, 
                 $\psio$  is the signed distance function 
                 after $N_\tau=256$ 
                 re-initialization steps.}}
 \label{fig8}
 \end{figure}
We found that 
errors of the solutions to \Eq{eq1} remain 
smaller than the truncation error
($L_1^M < 10^{-16}$) when
$\Delta \tau = \Delta x/2^M$ and $M=2,\ldots,8$.
The explanation
of this fact  is straightforward
if one notices that in \Eq{eq29} 
for $\eph \to 0$ the value of 
$\nabla \deltaa=\lr 1-2\alpha \rr \nabla \alpha$
increases when $x \ne x_\Gamma$. 
For practical reasons, 
later in this
section we use 
$\Delta \tau=\eph/2$.
In such case,  
the errors increase with 
decreasing $\eph=\Delta x/M$ 
which follow from  oscillations 
of the signed distance function 
$\psi_0$, 
see figures \ref{fig7}(b) and \ref{fig8}.

In the present test case
it is found that 
if $\eph = \Delta x /128$ 
then $\deltaa=\alpha \lr 1-\alpha \rr = 0$ 
since in the given grid
$\forall x_i \le x_P : \alpha_i = 0$
and $\forall x_i \ge x_F : \alpha_i = 1$ where
$x_P,x_F$ are 
the two neighboring 
grid points
closest to the interface
position $x_\Gamma=0.6\,m$.
Since $\deltaa = 0$ 
in each point of
the domain 
only when $\alpha \lr \psi_0 \rr = H \lr \psi_0 \rr$,
the Heaviside function $H\lr \psi_0 \rr$
is the numerical
solution to \Eq{eq29} as well.

The results presented 
in \Fig{fig7}
show the support width 
of $\deltaa/\eph$ 
is restricted by
the accuracy of 
reconstruction
of the level-set functions
$\alpps$ and $\psio$.
In order to accurately
calculate gradient of $\alpps$,
one needs at least 
four points around 
the interface $\Gamma$
with correctly predicted
values of $\psio$.
In the 1D study
presented here, this condition is
satisfied when $M \le 4$
in $\eph=\Delta x/M$
(compare results in
figures \ref{fig7}-\ref{fig8}).

Now, we can compare features
of this new re-initialization equation 
with properties of the stationary solution
to the phase-filed equation
which was investigated by 
Sun and Beckermann in \cite{sun2007}.
The main difference lies in 
the fact that the stationary solution 
to the phase-field equation 
is always given by the regularized
Heaviside function represented
by the hyperbolic tangent profile;
see formulation of the 
phase field equation 
in \cite{sun2007}.
For this reason,
in \cite{sun2007}
it is recommend 
to use about $5-6$
grid points  to accurately
reconstruct $\alpx$.
In the present method when $\eph \to 0$,
the numerical solution to \Eq{eq29} is
given by the Heaviside
function $\alpps = H\lr \psi_0\rr$
because in such case
$\delta \lr H \rr = 0$ on the
given grid, 
and thus $\pdl{\alpha}{\tau} \equiv 0$
in \Eq{eq29}.
In  subsequent sections, 
we investigate
how the selection 
of the interface 
width $\eph$ 
affects
re-initialization 
of the interface $\Gamma$
and computation 
of its curvature $\kappa$.

\subsection{Interpretation of results}
\label{sec5.3}

Let us now shortly 
summarize ideas introduced
in the previous sections.
Up to now 
the two discretizations of
$|\nabla \alpha|$ in \Eq{eq1}
were introduced:
$|\nabla \alpha|=\zeta^2\deltaa^{1-\gamma}/\gamma|\nabla \psi_1|$
and 
$|\nabla \alpha|=\deltaa/\eph|\nabla \psi_0|$, 
see \Eq{eq22} and \Eq{eq26}
respectively. 
These two discretizations 
are equivalent
when $0< \gamma \ll \Delta x$,
as shown by equations 
(\ref{eq19}),
(\ref{eq28}) 
and in \Fig{fig5}.
The discretization 
of $|\nabla \alpha|$ 
with \Eq{eq26} is free 
from round-off errors
and is faster than the equivalent
discretization with \Eq{eq22};
in the latter case, 
higher-order roots of $\alpps$ 
must be calculated, 
see \Eq{eq8}.

With these facts,
we can now explain why
equation (\ref{eq1}) is solved
with the largest accuracy  
when the mapping
functions $\psio$
or $\psi_0' \lr \psigs \rr$ are 
used for approximation
of $|\nabla \alpha|$.
Since  $\psi_0 \approx \psi_0'$ is
the signed distance function 
$|\nabla \psio| \approx |\nabla \psi_0' \lr \psigs \rr|=1$,
the right-hand sides in
equations (\ref{eq29}) and (\ref{eq30}) 
are equal to zero when  
the initial condition 
to \Eq{eq1} is given by \Eq{eq2}.
This occurs regardless of
$\kappa$ values and 
explains the rapid convergence 
of the solution to \Eq{eq1}
to a steady state 
as it is depicted
in Figs.~\ref{fig4} and \ref{fig6}.
Therefore, the simplification
of \Eq{eq1} to \Eq{eq4}
put forward by Shukla et al. in \cite{shukla10} 
is justified only when 
$|\nabla \alpha|= \deltaa/\eph |\nabla \psi_0|$
where $\psi_0$ is the signed
distance function.
When $\psio \approx \psi_0' \lr \psigs \rr$
is not the signed distance function,
for example,
due to $\alpps$ deformations during advection, 
the right-hand sides in \Eqs{eq29}{eq30} 
are not equal to zero and 
the re-initialization 
process will begin.

In Section \ref{sec5.2} 
we have shown that 
the stationary solution to \Eq{eq29} 
is also given by the Heaviside 
function $H \lr \psi_0 \rr$.
This distinguishes 
the new re-initialization method
from re-initialization of the signed distance function 
put forward by Sussman et al. \cite{sussman94,sussman98}
and  the solution to 
the phase-field equation
investigated by
Sun and Beckermann \cite{sun2007}.
In the 1D case,
this feature of equation (\ref{eq29})
allows decreasing 
the interface width 
up to $\eph=\Delta x/4$
without large 
influence on the accuracy 
of reconstruction
of the corresponding
$\alpps$ and $\psio$ functions
as shown in \Figs{fig7}{fig8}.
The latter result suggests
that for $K$ dimensional
problems $\eph=\sqrt{K}\Delta x/4$,
where $K=1,2,3$. 
In \cite{olsson05} and 
in works that followed,
the interface width is set 
to $\eph \ge \Delta x /2$
and it is reported that solutions
of the re-initialization equation (\ref{eq1})
are not stable when $\eph < \Delta x/2$.
The new form of the re-initialization equation
(\ref{eq1}) given by equation (\ref{eq29})
provides both the explanation and the solution
to this problem.

Finally, we recall when $\eph \to 0$ 
then $\alpps \to H \lr \psi_0 \rr$, 
     $\deltaa / \eph \to \delta \lr \psi_0 \rr$
and $\pdl{\alpha}{\tau} \equiv 0$ in equation (\ref{eq29}).
In this case, the  advection equation 
of the conservative level-set function 
\blnm
\be
 \pd{\alpps}{t} + w_i \pd{\alpps}{x_i}=
 \frac{\deltaa}{\eph}\ls \pd{\psio}{t} + w_i \pd{\psio}{x_i} \rs  = 0,
\label{eq30a}
\ee
\elnm
becomes the advection equation 
of the phase indicator function  $H \lr \psi_0 \rr$
which is discretized 
in the VOF family of methods, 
$w_i$ in equation (\ref{eq30a})
is the $i-th$ component
of the interface velocity.
Equation (\ref{eq30a}) shows
that advection of the phase indicator
function $H \lr \psi_0 \rr$ is equivalent
to advection of the signed distance 
function $\psi_0$ localized within the
support of the Dirac's delta 
function $\delta \lr \psi_0 \rr$.
When $\eph \ne 0$, 
equation (\ref{eq30a})
describes advection of  $\alpps$
and advection of $\psio$  
which is localized within 
the support of $\deltaa / \eph$. 
Equation (\ref{eq30a}) is valid
in the whole domain of the solution unlike 
the transport equation of $\psi_0 \lr \bx,t \rr$
derived in \cite{sussman94}, valid only 
at the interface $\Gamma$ 
($w_i \ne 0$ when $x_i \in supp \ls \delta \lr \psi_0 \rr \rs$ 
and $w_i = 0$ elsewhere).

\section{Numerical experiments}
\label{sec6}
%
In the following sections,
we investigate 
the rates of convergence 
and
assess  
numerical errors during
re-initialization of 
both
the conservative level-set function $\alpps$
and the signed distance function $\psio$.
Re-initialization is
performed using  \Eq{eq29}  
which is equivalent to \Eq{eq1}
when $|\nabla \alpha|$ is 
discretized with the 
signed distance function 
$\psi_0$ or its approximation
$\psizw$, see \Eq{eq3} and \Eq{eq19}, 
respectively.
In particular,
we are interested in 
the case when the interface
$\Gamma$ is stationary,
$\eph = \sqrt{K} \Delta x/4 < \Delta x /2$
where $K=1,2,3$ is the
dimension of the problem,  
and the  number of
re-initialization
steps $N_\tau \gg 1$.
In order to present properties
of the new re-initialization method
in a broader context,
several advection test cases,
where equations (\ref{eq29}) and (\ref{eq30a})
are solved alongside,
are performed.
The numerical errors 
during  computation
of $\alpps$ derivatives
are measures of 
artificial deformations 
of the interface $\Gamma$
due to a re-initialization 
process.
In the next sections,
their identification
is our main concern.

\subsection{Re-initialization of stationary interfaces}
\label{sec6.1}

\subsubsection{Regularized Heaviside function}
\label{sec6.1.1}

In figures \ref{fig2} and \ref{fig3} 
the solutions to equation (\ref{eq1}) 
with the mapping function
$\psigf$ after $N_\tau=1$
re-initialization steps 
were illustrated,
therein the width of the
interface is set to $\eph=\Delta x/2$
and the time step
$\Delta \tau = D/C^2=\eph$.
In what follows, 
we discuss 
the results obtained
with the same numerical 
setup as in section (\ref{sec2})
but after $N_\tau=256$ 
re-initialization steps.
The results presented below
are obtained with $\psigf$ where $\gamma_1=0.1$,
and with the new mapping 
functions $\psio$ and 
$\psigs$ where $\gamma_2=10^{-5}$, see \Sec{sec5.1}.
 \begin{figure}[h!] \nonumber
 \includegraphics[angle=-90]{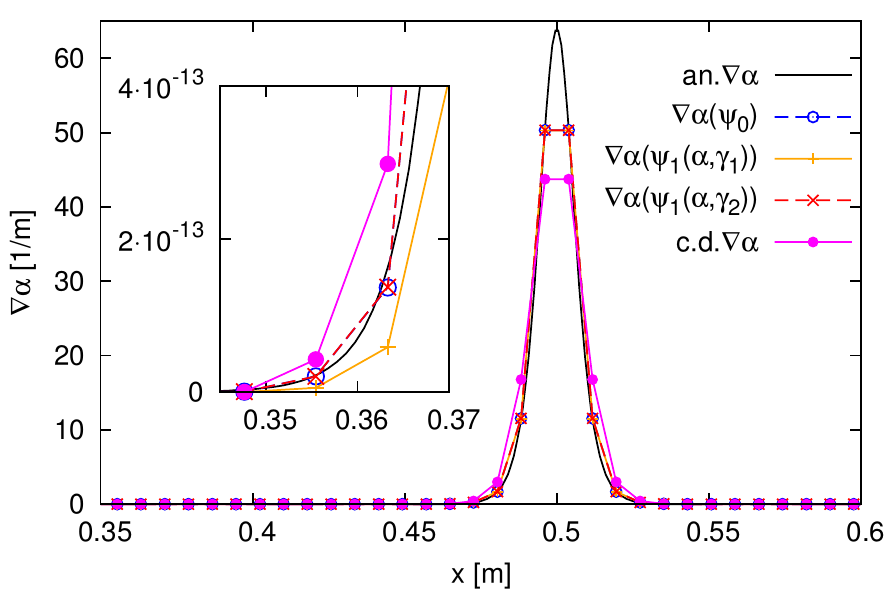}
 \caption{\small{The comparison of exact
                 $\nabla \alpha$ (black solid line)
                 with its numerical approximations  
                 after  $N_\tau=256$ 
                 re-initialization steps
                 of the 1D regularized Heaviside function.               
                 $|\nabla \alpha|$ in \Eq{eq1} is
                 discretized using the  mapping functions: 
                 $\psio$, $\psigf$ or $\psigs$ 
                 where $\gamma_1=0.1$, $\gamma_2=10^{-5}$
                 see \Eq{eq26} and \Eq{eq22}, respectively.
                 c.d.$\nabla \alpha$ denotes gradient 
                 of $\alpx$ computed using 
                 the central differencing 
                 scheme.}}
 \label{fig9}
 \end{figure}

\begin{figure}[h!] \nonumber
 \includegraphics[angle=-90]{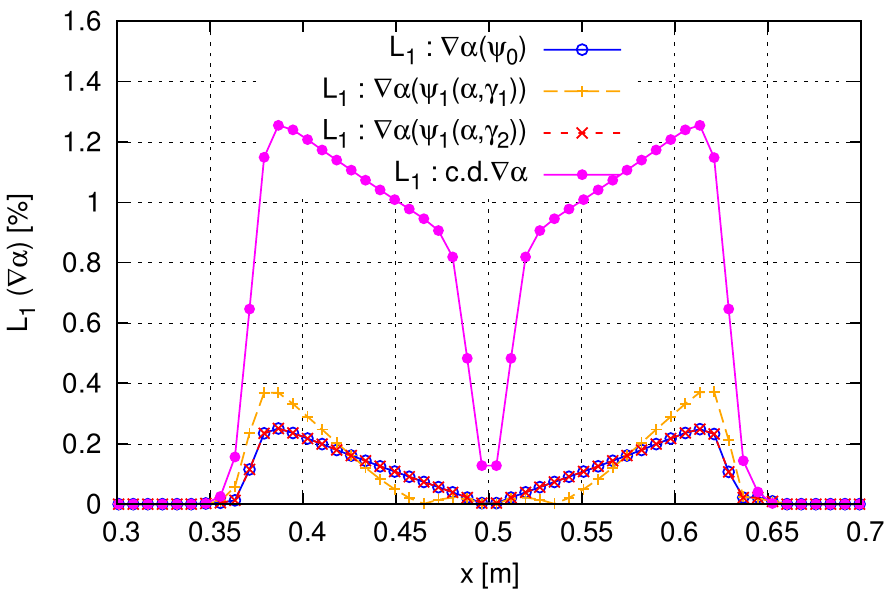}
 \caption{\small{Distributions of the  
                 $L_1 \lr \nabla \alpha \rr$ norms
                 defined by \Eq{eq10}  after  $N_\tau=256$
                 re-initialization steps
                 of the 1D regularized Heaviside function.
                 $L_1 \lr \nabla \alpha \rr$ is
                 calculated using the $\nabla \alpha$
                 discretizations  with the mapping functions:
                 $\psio$, $\psigf$ and $\psigs$ 
                 where $\gamma_1=0.1$, $\gamma_2=10^{-5}$,
                 c.d. $\nabla \alpha$ denotes 
                 the central differencing gradient $\alpx$
                 approximation.}}
 \label{fig10}
 \end{figure}

In figures \ref{fig9} and \ref{fig10},
approximations to the first component 
of $\nabla \alpha \lr \psi_0 \rr$
calculated using 
\Eq{eq22} and \Eq{eq26}
are presented.
In figure \ref{fig9}, it is observed
that both $\psio$ and $\psigs$
reconstruct the bell-like shape 
of the first-order analytical derivative of \Eq{eq2}
better than $\psigf$;
this result 
is expected in
the light of 
\Eq{eq19}.
The distributions of the
$L_1\lr \nabla \alpha \rr$ 
norms in \Fig{fig10} 
confirm these
observations.
We note that $\nabla \alpha$ 
is approximated with 
the smallest error 
in the neighborhood 
of the interface located 
at $x_\Gamma=0.5\,m$.

The distributions of the errors 
after $N_{\tau}=256$ re-initialization
steps in the case of  $\psio$ 
and $\psigs$
are  similar 
to the error distribution 
of  $\psigf$  after $N_{\tau}=1$ re-initialization steps
(compare \Fig{fig3} and \Fig{fig10}).
In the case of $\psigf$,
the $L_1\lr \nabla \alpha \rr$ norm 
(the numerical gradient 
 approximation within 
 \Eq{eq10})
is varying 
during 
the integration
of \Eq{eq1} in 
time $\tau$. 
This artificial deformation 
of  $\nabla \alpha$ in time
$\tau$ is the main reason for 
longer convergence of \Eq{eq1}
to the steady state
as shown in \Fig{fig4}.

Next,
we discuss  
the accuracy 
of computations 
of the 
second-order 
spatial 
derivatives
of the level-set
function
$\alpps$. 
In Section \ref{sec4},
the formulas 
for $\alpha_{,ij}$ 
using $\psig$ and $\psio$
were derived, see equations 
(\ref{eq25}) and
(\ref{eq27}), respectively.
To compute  
the second-order spatial derivatives 
of $\alpps$, one needs both 
the second and first order derivatives 
of the mapping functions $\psio$ or $\psig$.
 \begin{figure}[h!] \nonumber
 \includegraphics[angle=-90]{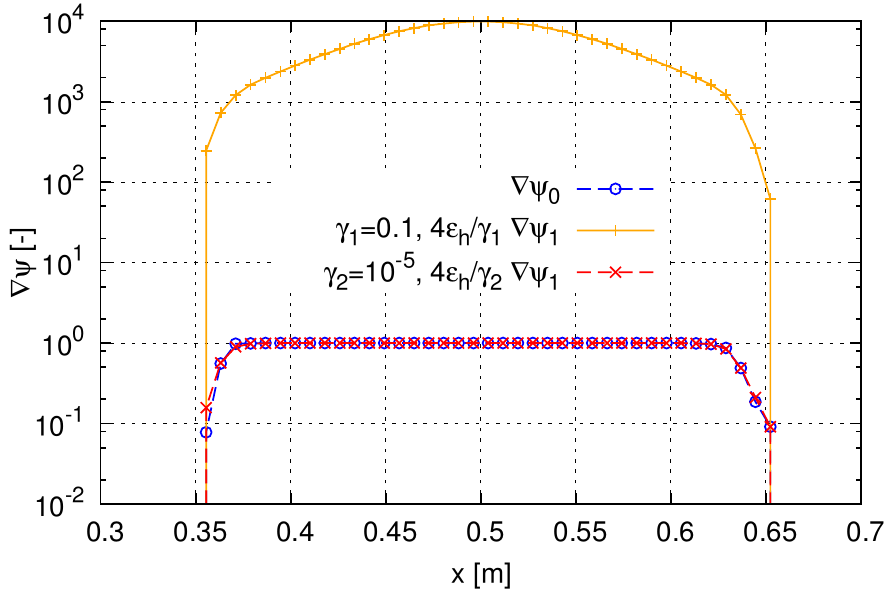}
 \includegraphics[angle=-90]{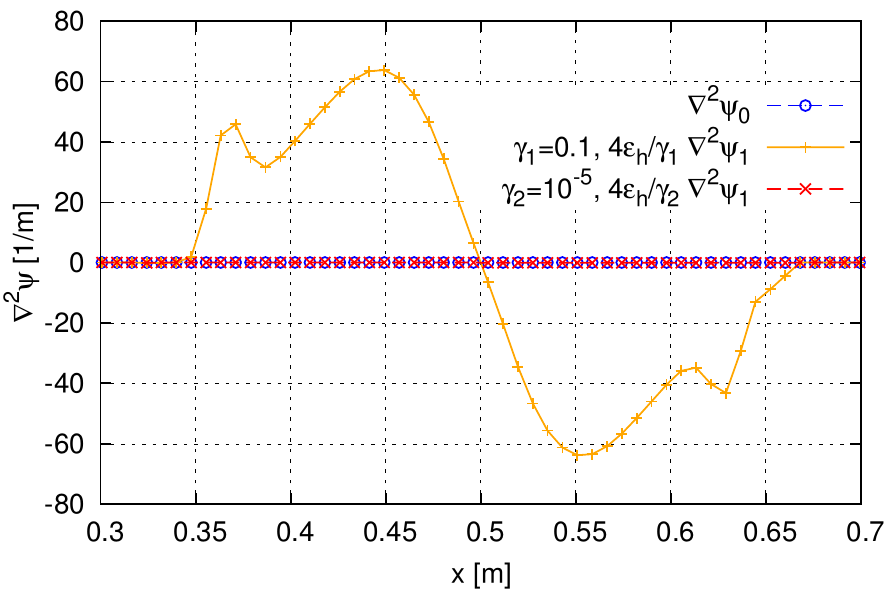}
 \caption{\small{The comparison
                 of (a) $\nabla \psi$, (b)  $\nabla^2\psi$ 
                 obtained after $N_\tau=256$ re-initialization 
                 steps of the 1D regularized Heaviside function.}}
 \label{fig11}
 \end{figure}
In our code they are calculated
using the discrete Gauss theorem
equivalent to the second-order accurate 
central difference gradient approximation, 
see \cite{peric02,schaefer06}.

In figure \ref{fig11}(a)
we compare the relation 
between $\nabla \psio$
and $4\eph/\gamma \nabla \psig$,
provided by \Eq{eq20}. 
One notes 
when $\psigf$
is used, 
the mapping 
function 
is not the signed 
distance function
since $|4\eph/\gamma_1 \nabla \psigf|\ne 1$.
On the other hand, 
gradient calculated
using $\psio$ and $\psi_0' \lr \psigs \rr$ 
gives 
$|\nabla \psi_0| \approx |\nabla \psi_0'| = 1$
inside the support of $\deltaa/\eph$;
with this, 
correctness 
of equations (\ref{eq19}) and (\ref{eq20})
is confirmed.

In figure \ref{fig11}(b),
the second-order 
spatial derivatives
of the mapping functions 
$\psio$ and $\psig$
are compared.
As $\psigf$ 
does not approximate 
the signed distance function,
its second-order derivative 
is not equal to zero everywhere
in the computational domain, 
see \Eq{eq21}.
Since $\psio$
is the signed distance
function    
and $\psipf$ 
is its approximation, 
see \Eq{eq19},
their second derivatives
are equal to zero
when $x \in supp \ls \deltaa/\eph \rs$.
Consequently, 
the accuracy 
of approximation
of the second-order derivatives 
using \Eq{eq25} and \Eq{eq27} 
is very similar
to accuracy 
achieved when 
the first-order derivatives
are computed with 
\Eq{eq22} and
\Eq{eq26}; 
compare results
in \Fig{fig10} and in \Fig{fig13}(a). 
 \begin{figure}[h!] \nonumber
 \includegraphics[angle=-90]{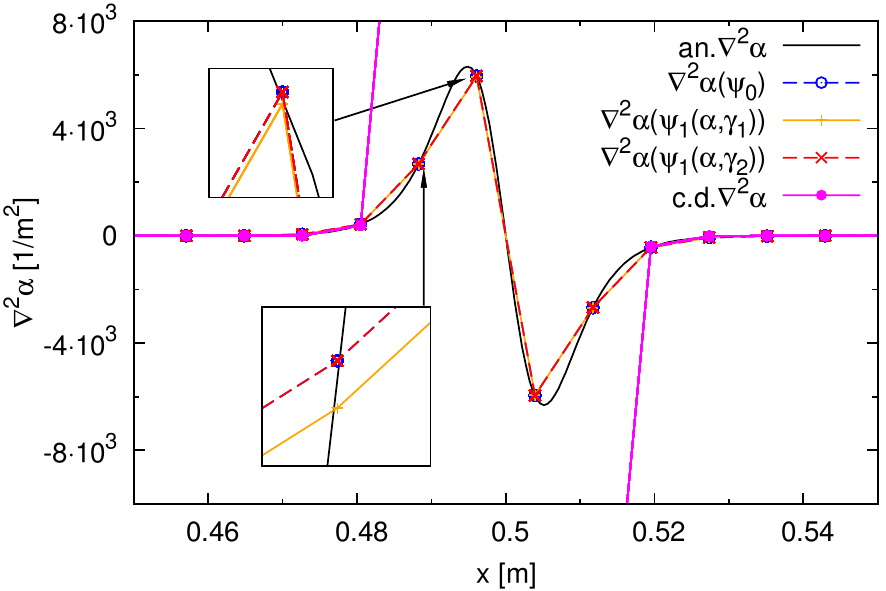}
 \caption{\small{The comparison of 
                 $\nabla^2\alpha\lr \psi \rr$
                 after $N_\tau=256$ re-initialization steps
                 of the 1D regularized Heaviside function
                 obtained with \Eq{eq25} or \Eq{eq27}
                 with the second-order derivative calculated 
                 analytically from \Eq{eq2} (solid-black line)
                 and the  central difference (c.d $\nabla^2\alpha$)
                 approximation.  
                 $\psi$ is one of the mapping functions: 
                 $\psix$, $\psigf$ or $\psigs$.}}
 \label{fig12}
 \end{figure}
 \begin{figure}[h!] \nonumber
 \includegraphics[angle=-90]{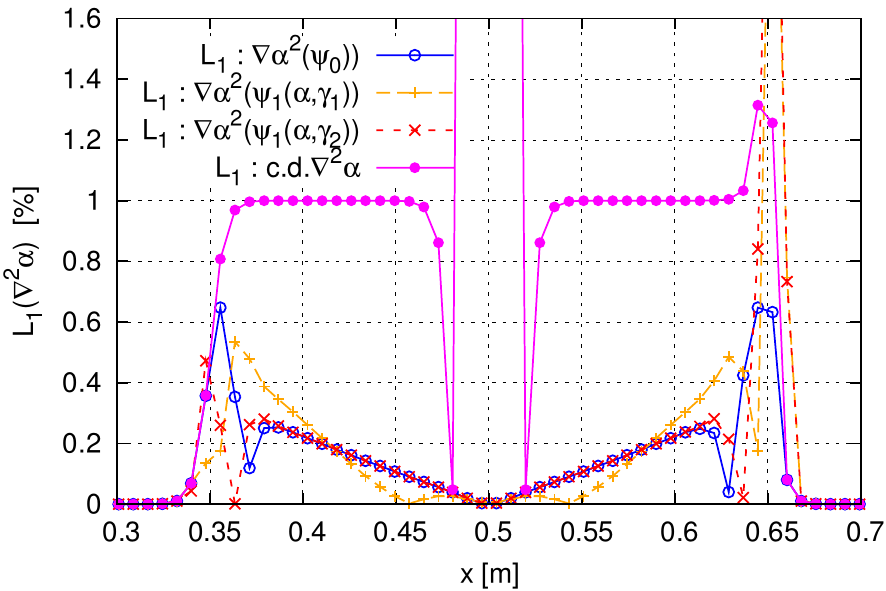}
 \includegraphics[angle=-90]{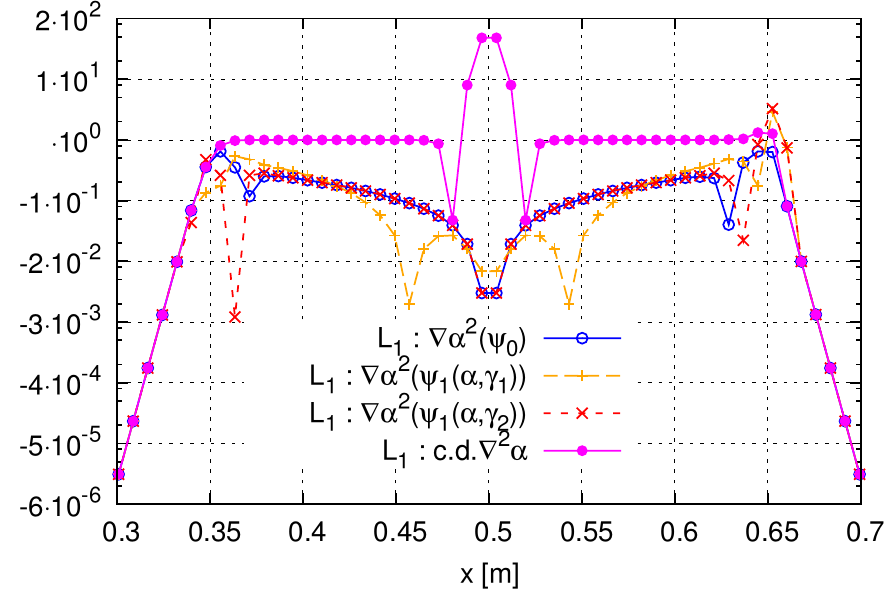}
\caption{\small{The comparison of 
                $L_1 \lr \nabla^2\alpha \lr \psi \rr \rr$ 
                norms after $N_\tau=256$ re-initialization steps
                of the 1D regularized Heaviside function.
                $\psi$ is one of the mapping functions: 
                $\psigf$, $\psigs$ or $\psio$,
                $L_1$ norm is defined by \Eq{eq10}. 
                c.d. $\nabla^2\alpha$
                denotes the second-order derivative computed
                using the central differencing method, 
                figures (a) and (b) present the same
                results but the Y-axis scales  
                are different.}}
 \label{fig13}
 \end{figure}

In figures 
\ref{fig12}-\ref{fig13},
the non-zero terms in 
$\nabla^2 \alpha$,
calculated with \Eq{eq25} 
and \Eq{eq27},  
are compared
with 
analytically calculated
second-order derivative 
of \Eq{eq2} and its 
central difference 
approximation.
The first observation in \Fig{fig12}
is a very good reconstruction of the
second-order
derivative in the case of all three
mapping functions; 
see also 
\Fig{fig13} where the distribution
of $L_1 \lr \nabla^2 \alpha \rr$ error
defined by \Eq{eq10} 
is depicted.
As expected, 
the second order 
derivatives of $\alpp$
calculated with $\psigs$ and $\psio$
are closest to each other and to the
analytical solution. 
In \Fig{fig13},
one
notices 
when 
$\psigs$ and $\psio$
are used, 
the differences 
between these two 
approximations 
are visible
only in points where 
the jumps of $\psig$
with $\epsilon=0$
in \Fig{fig1}
are present
(about $x=0.35\,m$ and 
$x=0.65\,m$).
Thus, 
oscillations
observed
in the $L_1 \lr \nabla^2\alpha \rr$ norm
can be attributed 
to the truncation errors
due to a floating point underflow.

In the light of 
the 1D 
re-initialization
studies presented
above 
we conclude 
that the
most accurate 
discretization
of $|\nabla \alpha|$ 
in equation (\ref{eq1})
is achieved 
with the mapping 
function
$\psio$, see
\Eq{eq26}.
This discretization 
is also the most 
natural one 
since the mapping 
between $\psio$ 
and $\alpps$
is well defined
by \Eq{eq3}.
We recall that
after proper selection
of the mapping function,
\Eq{eq1} is equivalent
to \Eqs{eq29}{eq30}
which are 
the conservative
and non-conservative
form 
of the re-initialization equation 
of the signed distance function $\psio$ and 
the conservative level-set function
$\alpps$. 
Since equation (\ref{eq1}) 
is solved accurately
when $x \in supp \ls \deltaa/\eph \rs$,
there 
is no need for 
an introduction
of additional 
techniques 
in the present 
solution 
procedure,
which
reconstruct 
the signed distance 
function $\psi_0$
in the vicinity
of the interface
$\Gamma$
.

\paragraph{Influence of initial conditions and the interface width on the convergence rate}
\label{sec6.1.1.1}

%
In most of the previous examples,
the interface $\Gamma$ was localized 
exactly in-between neighboring control volumes 
$x_P,\,x_F$; 
the width of the $\deltaa/\eph$ support 
was set to $\eph=\Delta x/2$
and  equation (\ref{eq1}) was initialized 
with its own analytical solution
given by equation (\ref{eq2}).
 \begin{figure}[h!] \nonumber
 \includegraphics[angle=-90]{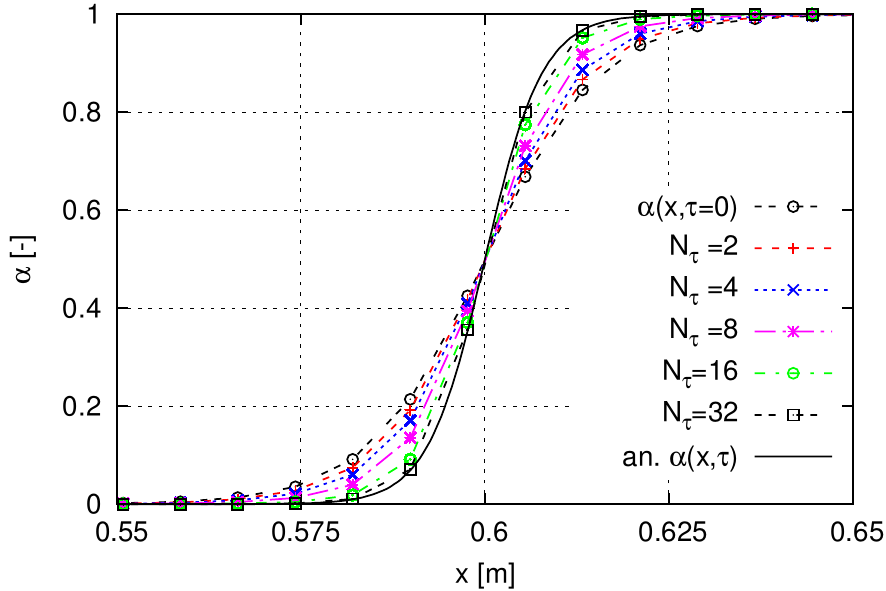}
 \includegraphics[angle=-90]{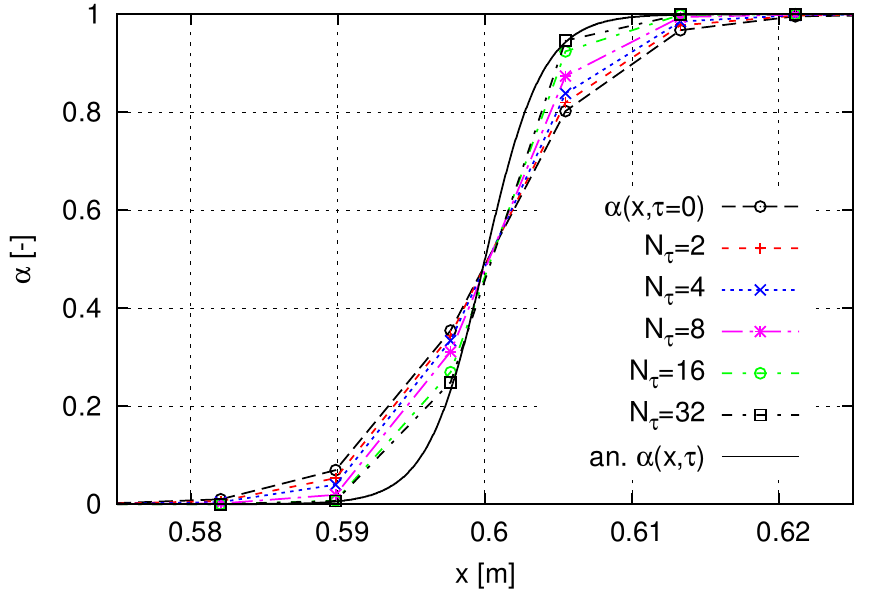}
 \caption{\small{Convergence of the level-set function $\alpps$ 
                 towards its analytical solution, 
                 $\alpha \lr \bx,\tau=0 \rr$
                 is set using $\epsilon_{h,0}=2 \eph$ where
                 (a) $\eph=\Delta x/2$, 
                 (b) $\eph=\Delta x/4$,
                 the time step size $\Delta \tau=\eph/2$.
                 $|\nabla \alpha|$ in \Eq{eq1} 
                 is discretized with $\psio$.}}
 \label{fig14}
 \end{figure}
However, during advection of 
the level-set functions $\alpps$ and $\psio$, 
the re-initialization equation (\ref{eq29})
must handle more general
initial conditions.
 \begin{figure}[h!] \nonumber
 \includegraphics[angle=-90]{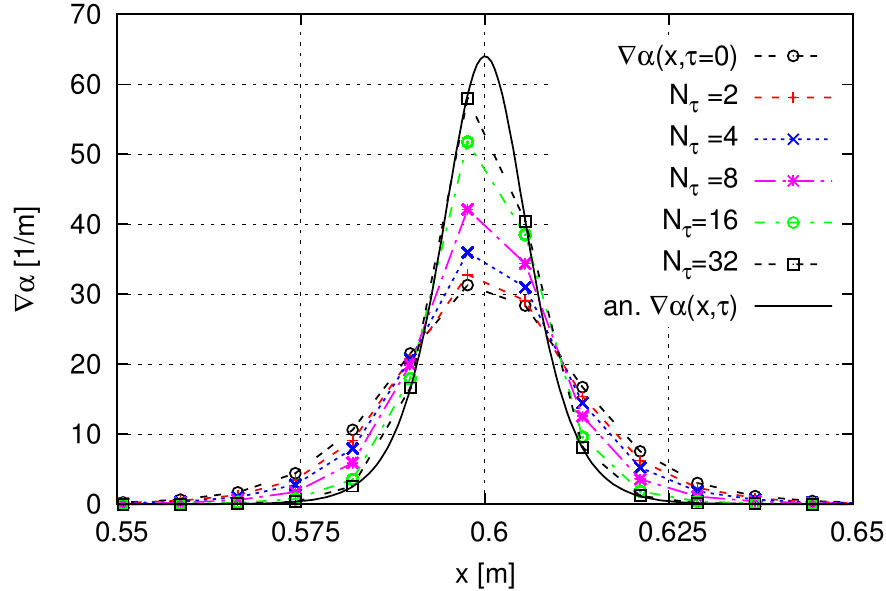}
 \includegraphics[angle=-90]{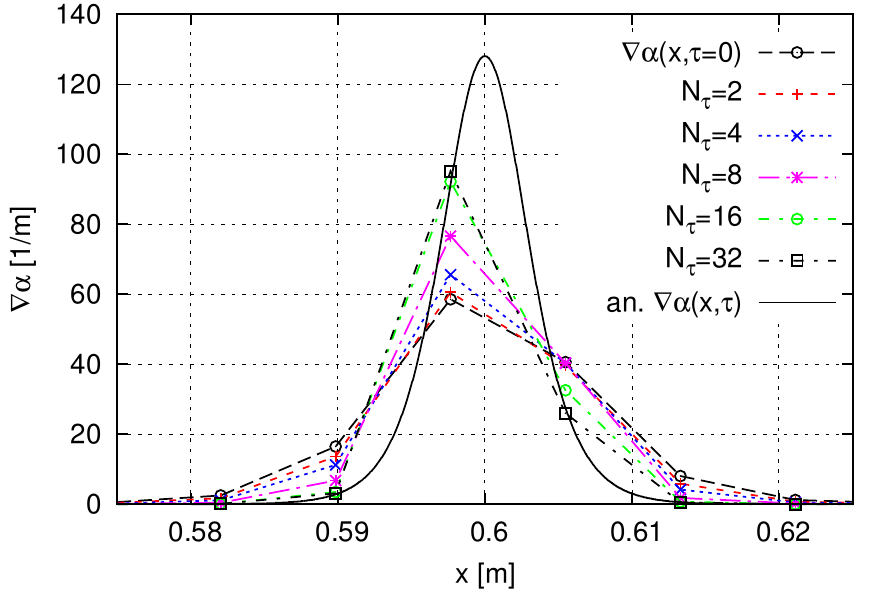}
 \caption{\small{Convergence of $\nabla \alpps$  
                 towards its analytical solution,  
                 $\alpha \lr \bx,\tau=0 \rr$
                 is set using $\epsilon_{h,0}=2 \eph$ where:
                 (a) $\eph=\Delta x/2$, 
                 (b) $\eph=\Delta x/4$,
                 the time step size $\Delta \tau=\eph/2$.
                 $|\nabla \alpha|$ in \Eq{eq1} 
                 is discretized with $\psio$.}}
 \label{fig15}
 \end{figure}
 \begin{figure}[h!] \nonumber
 \includegraphics[angle=-90]{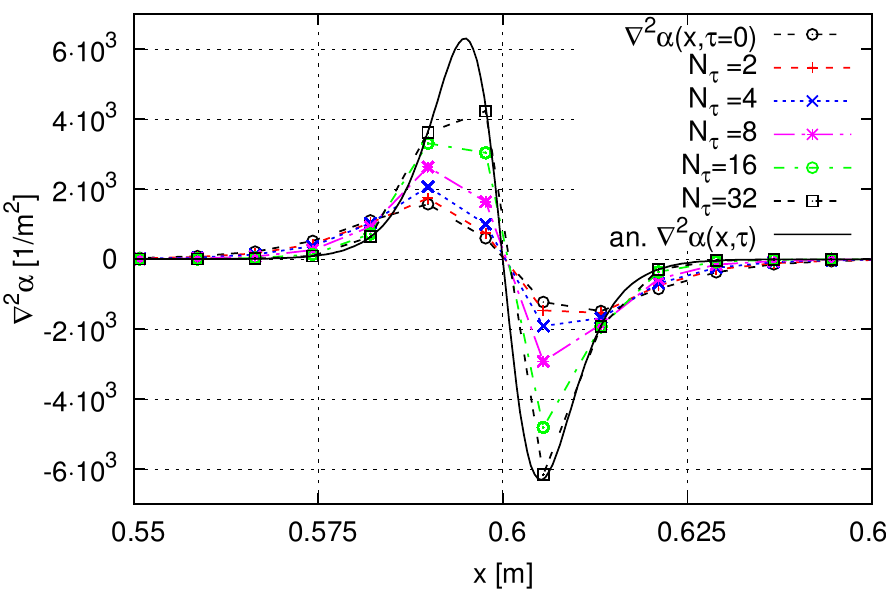}
 \includegraphics[angle=-90]{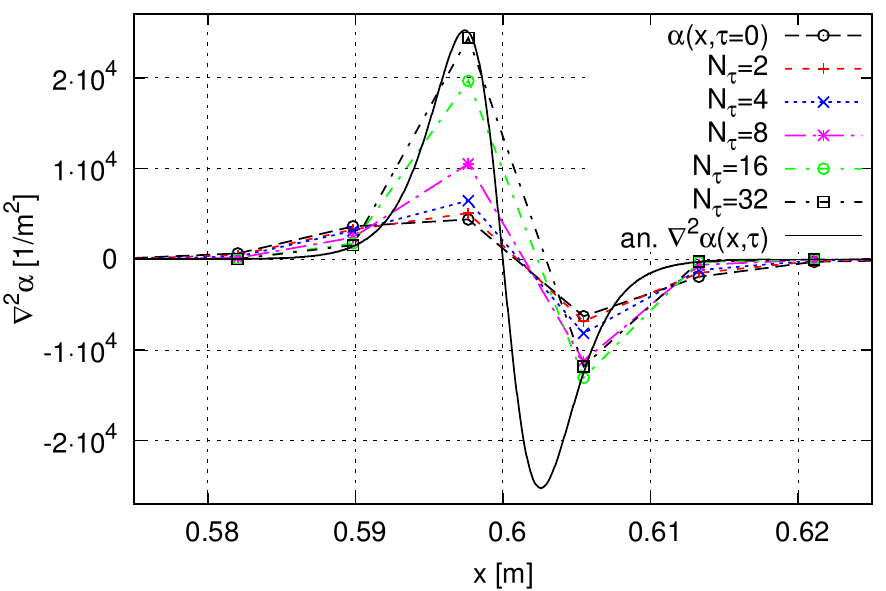}
 \caption{\small{Convergence of $\nabla^2 \alpps$  
                 towards its analytical solution,  
                 $\alpha \lr \bx,\tau=0 \rr$
                 is set using $\epsilon_{h,0}=2 \eph$ where:
                 (a) $\eph=\Delta x/2$, 
                 (b) $\eph=\Delta x/4$,
                 the time step size $\Delta \tau=\eph/2$.
                 $|\nabla \alpha|$ in \Eq{eq1} 
                 is discretized with $\psio$.}}
 \label{fig16}
 \end{figure}
For this reason,
in what follows we 
study influence of 
the arbitrary interface location 
(here $x_\Gamma=0.6\,m$), 
and the support width $\eph$
on the convergence rate 
of $\alpps$ and its spatial 
derivatives during re-initialization
of the 1D regularized Heaviside function.

In this study,
the initial condition to \Eq{eq1} 
is not given by its analytical solution;
at $\tau=0$ we set $\epsilon_{h,0}=2 \eph$
and we consider two widths 
of the interface:
$\eph=\Delta x/2$ 
and $\eph=\Delta x/4$.
The discretization of the 1D computational domain 
is the same as described in \Sec{sec6.1.1}.
For brevity, 
only the results obtained with $\psio$ 
used in discretization of  $|\nabla \alpha|$
in \Eq{eq1} are presented
in this section.

\begin{figure}[h!] \nonumber
  \includegraphics[angle=-90]{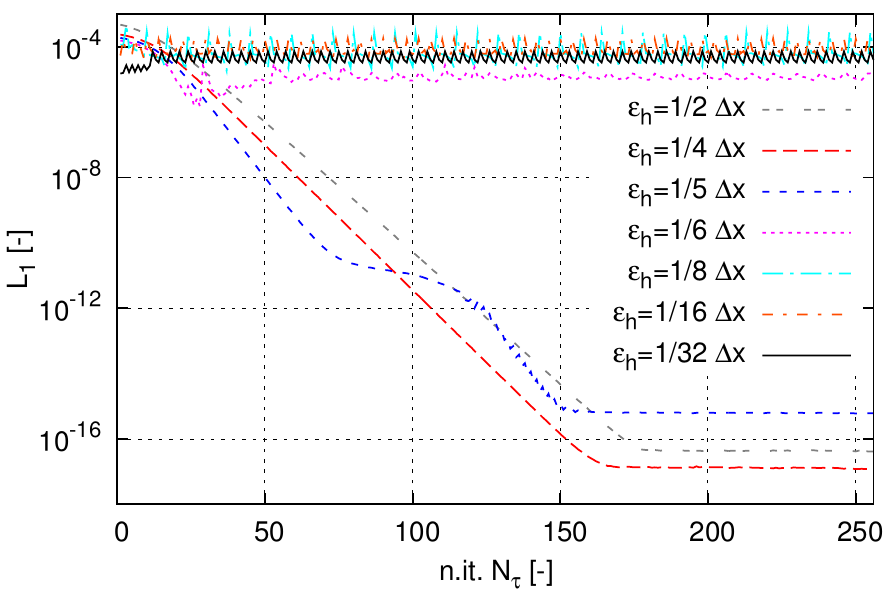}
  \caption{\small{The $L_1$ norms defined by \Eq{eq9}
                  obtained during re-initialization of 
                  the 1D regularized Heaviside function $\alpps$
                  with \Eq{eq29} and the variable interface width $\eph$.
                  The interface is located at $x_\Gamma=0.6\,m$, 
                  at $\tau=0$ $\epsilon_{h,0}=2 \eph$, 
                  $\Delta \tau =\eph/2$,
                  $\eph=\Delta x/M$ and $M=2,\ldots,32$.}}
 \label{fig17}
 \end{figure}

In figures \ref{fig14}-\ref{fig16} 
it is observed that
the numerical solution to \Eq{eq1} 
converges towards its analytical 
counterpart independent from 
the selected final support
width of $\deltaa/\eph$.
We emphasize that re-initialization
of $\alpps$ with $\epsilon_{h,0}=2\eph$
is also possible when $\eph=\Delta x/M$. 
When $M > 4$ the solution of \Eq{eq29} 
is convergent but the error level in 
the representation of $\alpps$ and $\psio$ 
grows, see \Fig{fig17}.
As it is discussed 
in Section \ref{sec5.2},
this occurs due to finite spatial 
and temporal resolutions.

Since the interface $\Gamma$ 
is not localized exactly in-between 
neighboring control volumes $x_P,\,x_F$, 
obtained solutions are  not symmetrical
as depicted in \Figs{fig14}{fig16}.
In spite of this,
the shapes of the level-set function 
$\alpps$ and its first and second-order 
spatial derivatives are correctly 
reconstructed.
We note
that in the limit of $N_\tau \to \infty$,
the present numerical solution to 
equation (\ref{eq29}) tends to 
its stationary analytical solution 
given by equation (\ref{eq2})
and its first and second
order spatial derivatives.
The discretizations 
of equation (\ref{eq1})
presented in the extant literature
do not guarantee 
convergence towards
its stationary analytical solution,
artificial deformations 
of the interface $\Gamma$ 
that emerge when $N_\tau \gg 1$
are the result 
of this inconsistency.
%

\subsubsection{Circular interface}
\label{sec6.1.2}

In the case
of the 1D regularized Heaviside function 
re-initialized in Section \ref{sec6.1.1}, 
equations (\ref{eq1}) and  (\ref{eq4})
are equivalent  since 
$\kappa = -\nabla \cdot \bng \equiv 0$.
In 2D or 3D cases  
that are  discussed next,
$\kappa \ne 0$. 
Choosing the mapping function 
as the signed
distance function 
$\psio$ allows us
to write equation (\ref{eq1}) 
in the form of equations (\ref{eq29}) or (\ref{eq30})
and solve one of these equations
without neglecting the influence
of $\kappa$ on 
the re-initialization process.
When equation (\ref{eq1}) is solved
correctly, equation (\ref{eq3}) that 
defines the signed distance function $\psio$ 
holds;  we use this fact in 
the present solution procedure.

In what follows, we assess 
influence of the selected
mapping function and 
the interface width 
($\eph=\Delta x/2$ or $\eph=\sqrt{2}\Delta x/4$)
on the accuracy of approximations
of $\alpps$ derivatives, 
see \Eqs{eq22}{eq26} 
or \Eqs{eq25}{eq27},
respectively.  
These derivatives 
are used in calculation 
of the interface  
curvature $\kappa$ according
to the formula 
\blnm
\begin{eqnarray}
 \begin{split}
   \kappa= & \lr \alpha_1^2 \alpha_{22} 
               + \alpha_2^2 \alpha_{11} 
               + \alpha_1^2 \alpha_{33} 
               + \alpha_3^2 \alpha_{11} 
               + \alpha_2^2 \alpha_{33} 
               + \alpha_3^2 \alpha_{22}   \right. \\ & \left.
               - 2 \alpha_1 \alpha_2 \alpha_{12} 
               - 2 \alpha_1 \alpha_3 \alpha_{13} 
               - 2 \alpha_2 \alpha_3 \alpha_{23} \rr / |\nabla \alpha|^3 \\
         = &  \lc \psi_{0,1}^2 \ls \psi_{0,22} + \frac{\psi_{0,2}^2}{\eph} \lr 1-2 \alpha \rr \rs 
               +   \psi_{0,2}^2 \ls \psi_{0,11} + \frac{\psi_{0,1}^2}{\eph} \lr 1-2 \alpha \rr \rs \right. \\ & \left.
               +   \psi_{0,1}^2 \ls \psi_{0,33} + \frac{\psi_{0,3}^2}{\eph} \lr 1-2 \alpha \rr \rs
               +   \psi_{0,3}^2 \ls \psi_{0,11} + \frac{\psi_{0,1}^2}{\eph} \lr 1-2 \alpha \rr \rs \right. \\ &  \left.
               +   \psi_{0,2}^2 \ls \psi_{0,33} + \frac{\psi_{0,3}^2}{\eph} \lr 1-2 \alpha \rr \rs
               +   \psi_{0,3}^2 \ls \psi_{0,22} + \frac{\psi_{0,2}^2}{\eph} \lr 1-2 \alpha \rr \rs \right. \\ & \left.
               - 2 \psi_{0,1} \psi_{0,2} \ls \psi_{0,12} + \frac{\psi_{0,1} \psi_{0,2}}{\eph} \lr 1-2\alpha \rr \rs   \right. \\ & \left.
               - 2 \psi_{0,1} \psi_{0,3} \ls \psi_{0,13} + \frac{\psi_{0,1} \psi_{0,3}}{\eph} \lr 1-2\alpha \rr \rs   \right. \\ & \left.
               - 2 \psi_{0,2} \psi_{0,3} \ls \psi_{0,23} + \frac{\psi_{0,2} \psi_{0,3}}{\eph} \lr 1-2\alpha \rr \rs   \rc
               / |\nabla \psi_0|^3 ,
  \label{eq31k}
 \end{split}
\end{eqnarray}
\elnm
which is written for the
mapping function $\psio$ and
is valid when 
$x_i \in supp \ls \deltaa/\eph \rs$.
Unlike in the standard 
approach which uses only 
the signed distance function $\psi_0$,
equation (\ref{eq31k}) contains
terms that contribute to $\kappa$ 
exclusively away from 
the interface $\Gamma$. 
These terms are multiplied by 
factor $\lr 1-2\alpha \rr$
and they vanish at $\Gamma$,
i.e., 
when $\alpha \lr \psi_0 = 0 \rr =1/2$.
At the interface $\Gamma$ 
equation (\ref{eq31k}) reduces 
to $\kappa$ definition 
given in \cite{osher03}.

In the following sections 
we investigate 
the convergence rate 
of the circular interface curvature 
on five gradually refined meshes 
$m_i= 2^{4+i} \times 2^{4+i}$
where $i=1,\dots,5$. 
The initial condition to \Eq{eq1}
is given by \Eq{eq2} where
\blnm
 \be
  \psi_0 \lr \bx,\tau=0 \rr =  \ls \sum_{i=1}^2 \lr x_i-x_{0,i} \rr^2 \rs^{1/2} - R,
 \label{eq31}
 \ee
\elnm
$\lr x_{0,1},x_{0,2} \rr =\lr 0.5\,m, 0.5\,m \rr$ 
denotes the center of the circle  
with the radius $R=0.2\,m$. 
In this test case,
the computational domain 
is quadratic box 
$\Omega = <0,1>\times<0,1>\,m^2$,
and the number of grid nodes
depends on the size of 
the grid $m_i$.
%

%
\paragraph{Convergence of the re-initialization equation}
\label{sec6.1.2.1}

\begin{figure}[ht!] \nonumber
 \begin{centering}
 \begin{minipage}{.5\textwidth}
  \centering
  \includegraphics[width=.666\textwidth,height=1.\textwidth, angle=-90]{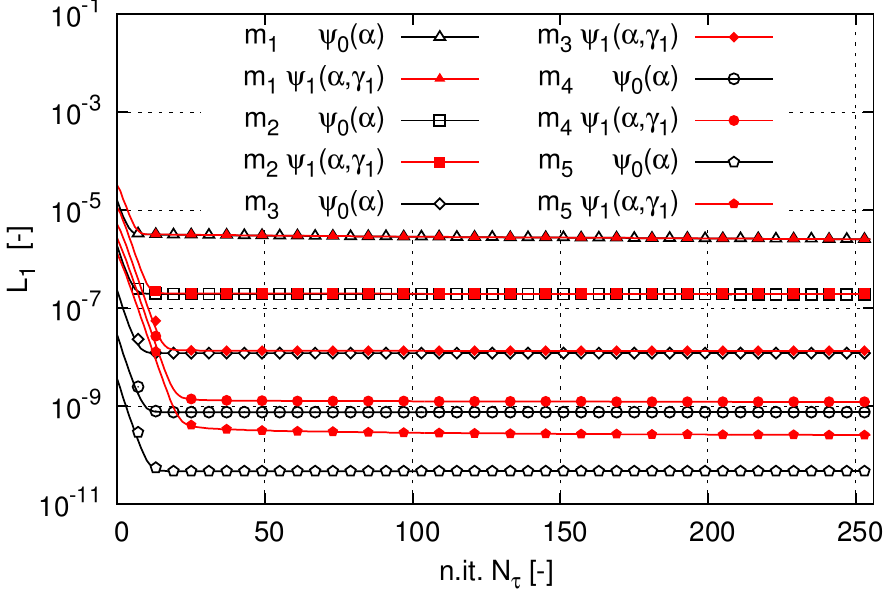}
 \end{minipage}%
 \begin{minipage}{.5\textwidth}
  \centering
  \includegraphics[width=.666\textwidth,height=1.\textwidth,angle=-90]{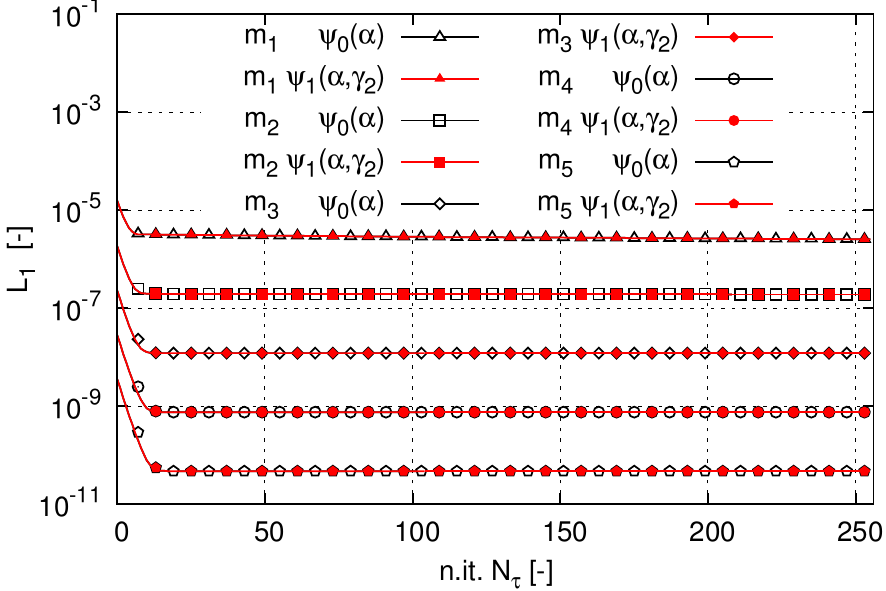}
 \end{minipage}
 \begin{minipage}{.5\textwidth}
  \centering
  \includegraphics[width=.666\textwidth,height=1.\textwidth,angle=-90]{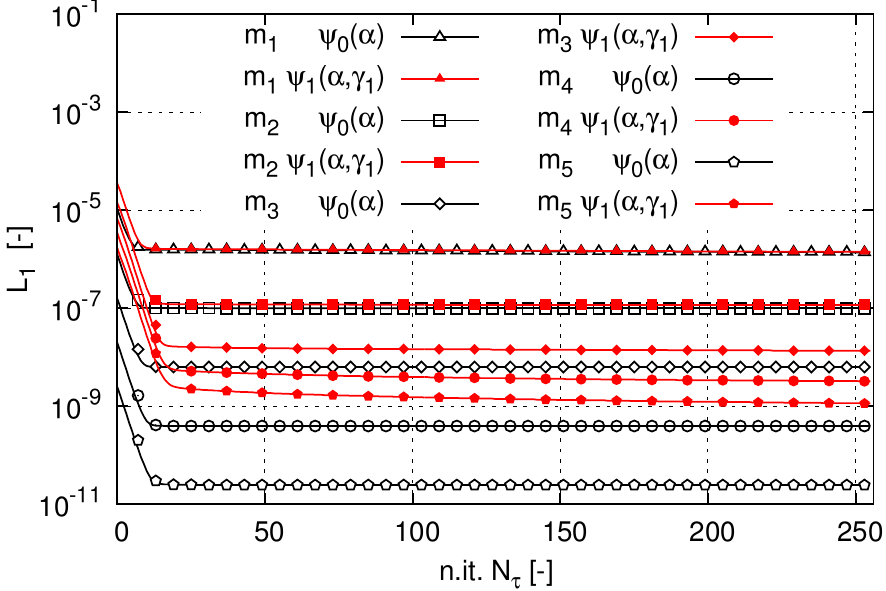}
 \end{minipage}%
 \begin{minipage}{.5\textwidth}
  \centering
  \includegraphics[width=.666\textwidth,height=1.\textwidth,angle=-90]{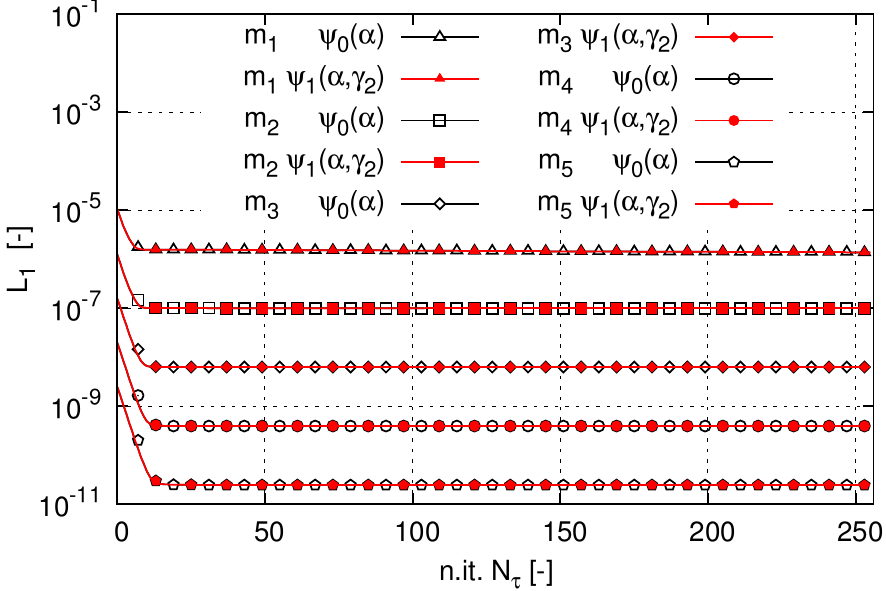}
 \end{minipage}
 \end{centering}
 \caption{\small{Convergence of the solution to \Eq{eq1} 
                 during re-initialization of the 2D circular interface,
                 $L_1$ norms defined by \Eq{eq9} 
                 are obtained for
                 the mapping functions:
                 $\psigf$, $\psigs$ and $\psio$,
                 the interface width
                 is set to $\eph=\Delta x/2$ (top)
                 and $\eph=\Delta x \sqrt{2}/4$ (bottom), 
                 $\Delta \tau=\eph$,
                 symbols correspond to
                 every sixth or every twelfth 
                 iteration in time $\tau$.}}

 \label{fig18}
 \end{figure}

In what follows, 
the convergence rate 
to the solution to \Eq{eq1}
during $N_\tau=256$ 
re-initialization steps
on grids $m_i$, $i=1,\dots,5$
is presented.
We compare the results
obtained with two $\deltaa/\eph$ support
widths: $\eph= \Delta x/2$,
$\eph=  \sqrt{2}\Delta x/4$,
where the
time step size is
$\Delta \tau=C/D^2=\eph$.
Unlike in the 1D case, 
the convergence rates and $L_1$ norms 
on gradually refined grids
are practically the same 
when $\psigs$ and $\psio$  
mapping functions are used 
(compare results in figures
 \ref{fig18}(b)(d)
 and in figure \ref{fig4}).
These results 
again confirm 
the correctness 
of the relation
given by \Eq{eq19}.
In the case of $\psio$, $\psigs$
and for both interface widths $\eph$
the stationary solution 
is achieved after about $N_{\tau} \approx 10$
iterations in time $\tau$.
The integration of \Eq{eq1} with 
the mapping function $\psigf$ requires
about $N_{\tau} \approx 25$ iterations
to achieve the steady state solution, 
see figures \ref{fig18}(a)(c).
In this latter case,
the convergence rate is lower 
and the error level is higher
when compared with
the results obtained 
with 
$\psio$ and $\psigs$. 

Lack of immediate convergence 
of the numerical solution to 
a steady state
(as in the 1D case depicted in \Fig{fig4}) 
is explained
by additional numerical errors 
which are introduced to the solution 
of \Eq{eq1} during discretization
process in the 2D case.
The two sources of 
the numerical errors
can be identified: 
the second-order
discretization of the fluxes
on the RHS of \Eq{eq1},
and the computation 
of $\psio$  
gradients and normals 
to the interface $\Gamma$
with the central-differencing
scheme 
which is 
known to be
mathematically 
exact approximation
of the spatial derivatives
only in
the 1D case,
see \cite{peric02}.

During 
numerical
experiments 
it was found 
that
differences 
between  rates
of convergence
and their 
levels for the
$\psio$, $\psigs$ 
and $\psigf$
mapping functions 
are more pronounced
when
$\Delta x/4 \le \eph < \Delta x \sqrt{2}/4$ 
and
$\Delta \tau <  D/C^2$.
For the sake of brevity, 
we subsequently use the
time step size
$\Delta \tau=D/C^2=\eph$
for the two different 
interface widths $\eph=\Delta x/2$
and $\eph=\Delta x \sqrt{K}/4$ where $K=2,3$ 
for 2D or 3D problems, 
respectively.

\begin{figure}[h!] \nonumber
 \begin{centering}
 \begin{minipage}{.75\textwidth}
  \centering
  \includegraphics[width=.6428\textwidth,height=1.\textwidth,angle=-90]{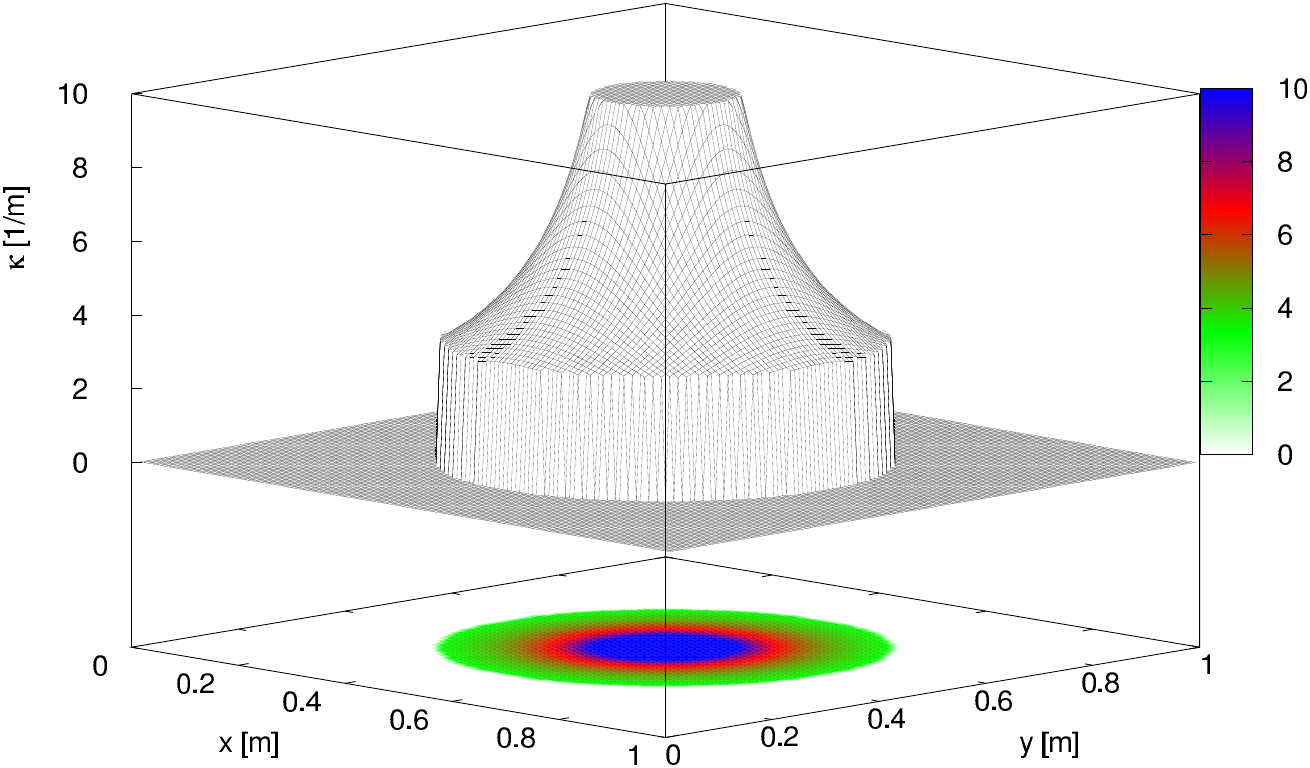}
 \end{minipage}
 \end{centering}
 \caption{\small{The curvature field $\kappa_3'$ 
                 after $N_\tau=256$ re-initialization
                 steps of the 2D circular interface 
                 obtained on the grid $m_3$ 
                 with the mapping function $\psio$.
                 The interface width
                 is set to $\eph=\Delta x/2$.}}
 \label{fig19}
 \end{figure}

\begin{figure}[ht!] \nonumber
 \begin{centering}
 \begin{minipage}{.33\textwidth}
  \centering
  \includegraphics[width=1.\textwidth,height=1.\textwidth,angle=-90]{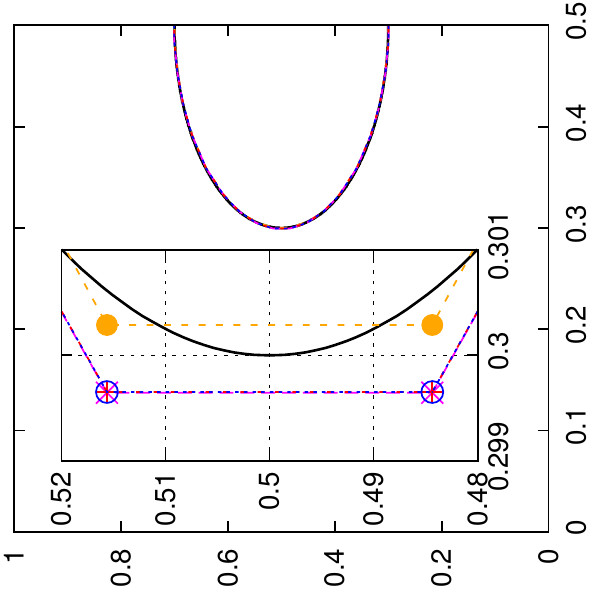}
 \end{minipage}%
 \begin{minipage}{.33\textwidth}
  \centering
  \includegraphics[width=1.\textwidth,height=1.\textwidth,angle=-90]{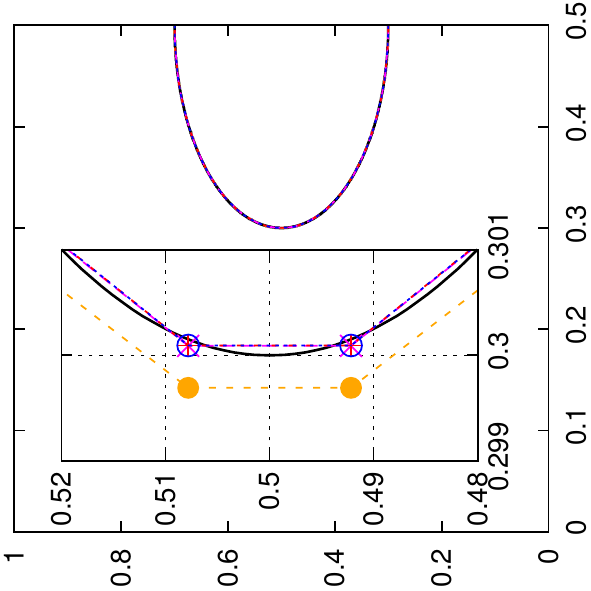}
 \end{minipage}%
 \begin{minipage}{.33\textwidth}
  \centering
  \includegraphics[width=1.\textwidth,height=1.\textwidth,angle=-90]{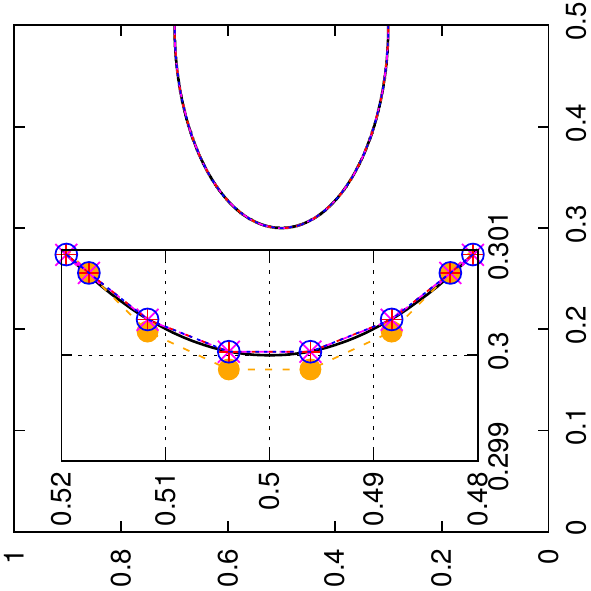}
 \end{minipage}
 \begin{minipage}{.33\textwidth}
  \centering
  \includegraphics[width=1.\textwidth,height=1.\textwidth,angle=-90]{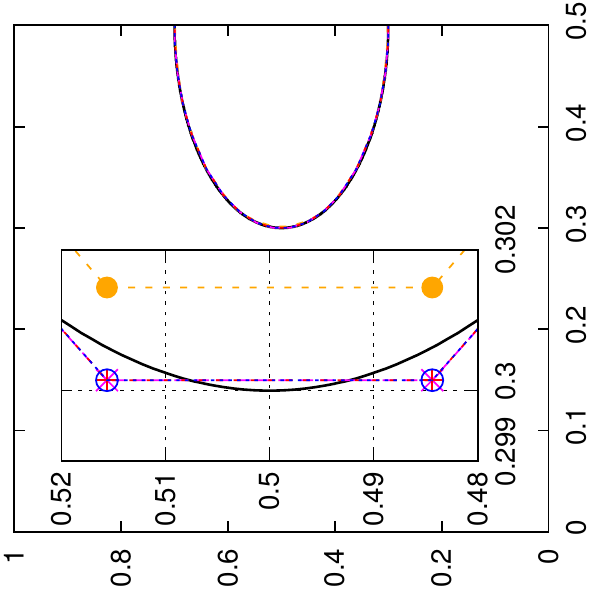}
 \end{minipage}%
 \begin{minipage}{.33\textwidth}
  \centering
  \includegraphics[width=1.\textwidth,height=1.\textwidth,angle=-90]{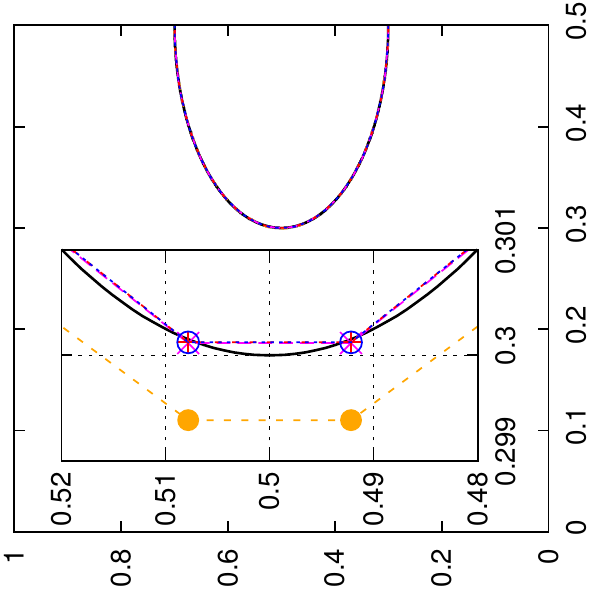}
 \end{minipage}%
 \begin{minipage}{.33\textwidth}
  \centering
  \includegraphics[width=1.\textwidth,height=1.\textwidth,angle=-90]{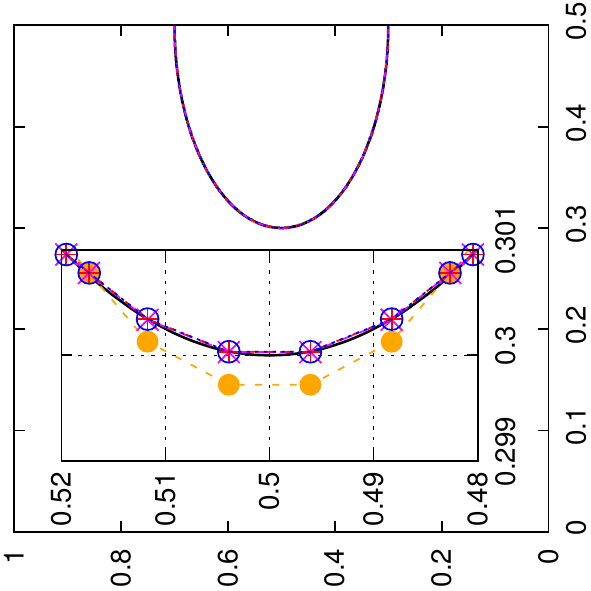}
 \end{minipage}
 \end{centering}
 \caption{\small{Convergence of different interface representations 
                 towards the analytical interface (black solid line) 
                 after $N_\tau=256$ re-initialization steps.
                 The interface width 
                 is set to $\eph=\Delta x/2$ (top) and 
                 $\eph=\sqrt{2} \Delta x/4$ (bottom).                 
                 The interface $\Gamma$ is captured using: 
                 the regularized Heaviside function $\alpp$ (solid orange dots),
                 the signed distance functions 
                 $\psizw \lr \psig \rr$ where $\gamma_1=0.1$ (red crosses), 
                 $\gamma_2=10^{-5}$ (magenta stars),
                 and $\psio$ (void dark blue dots) on grids: 
                  $m_1$, $m_2$, $m_3$ from left to right.}}
 \label{fig20}
 \end{figure}

\paragraph{Computation of the circular interface curvature}  
\label{sec6.1.2.2}
%
Next, we compute 
a numerical approximation $\kappa'$ 
of the exact curvature $\kappa$ 
using equation (\ref{eq31k}),
and we investigate 
its convergence rate
on five gradually 
refined grids after $N_\tau=256$ 
re-initialization steps
when $\eph=\Delta x/2$
or  $\eph=\sqrt{2}\Delta x/4$.
Since $\alpp$ is 
the level-set function, 
$\kappa'$ is calculated not only
at the interface 
$\alpha \lr \bx_\Gamma,t\rr=0.5$
but also
at $\alpha \lr \bx^1,t\rr=0.05$
and at $\alpha \lr \bx^2,t\rr=0.95$, 
see \Fig{fig19}. 
The exact curvature $\kappa_i$
of a circle on grid $m_i$ 
at the interface $\bx_{i,\Gamma}$ 
is constant and  equal 
to $\kappa_i=1/\mathcal{R}_i=1/\lr |\bx_{i,\Gamma}-\bx_0|\rr$.
$\mathcal{R}_i$ is the numerical approximation
to $R$ and is determined separately
for $\alpps$, 
    $\psizw \lr \psi_1 \rr$ 
and $\psio$
on each grid $m_i$,
$i=1,\dots,5$, 
see \Fig{fig20}. 
In the cases 
$\alpha \lr \bx^1_i,t \rr$
and $\alpha \lr \bx_i^2,t \rr$ 
curvatures are defined by
$\kappa_1=1/(\mathcal{R}_i+r_i^1)$
and $\kappa_2=1/(\mathcal{R}_i+r_i^2)$
where
\blnm
\bea
\label{eq32}
   &&r_i^1 = \epsilon_{h,i} \ln{\ls \frac{\alpha \lr \bx_i^1,t \rr+\epsilon}{1-\alpha \lr \bx_i^1,t \rr+\epsilon}\rs},\\ 
   &&r_i^2 = \epsilon_{h,i} \ln{\ls \frac{\alpha \lr \bx_i^2,t \rr+\epsilon}{1-\alpha \lr \bx_i^2,t \rr+\epsilon}\rs},
\label{eq33}
\eea
\elnm
$\epsilon_{h,i}$ 
depends on the 
grid size and
$\epsilon=5\cdot 10^{-16}$
is a small constant. 

We note when the mapping function 
is used in discretization of \Eq{eq1},
two representations of the interface $\Gamma$ 
do exist.
These two representations are given by: 
$\alpha \lr \bx_\Gamma,t\rr=1/2$
and 
$\psizw \lr 1/2 \rr =0$ or 
$\psi_0 \lr 1/2 \rr=0$,
and they are equivalent 
when $\eph \to 0$,
see \Fig{fig20}.
For this reason 
the $L_1^\kappa$ error 
is calculated as follows, first, 
the interface iso-lines are computed for each grid $m_i$:
\begin{itemize}
\item[1.]{ $\bx_i^\alpha$ from $\alpha \lr \bx_\Gamma,t\rr=1/2$,}
\item[2.]{$\bx_i^{\psi_0'}$      from $\psizw \lr 1/2 \rr =0$, see \Eq{eq19},}
\item[3.]{$\bx_i^{\psi_0}$   from $\psi_0 \lr 1/2 \rr= 0$}, see \Eq{eq3}.
\end{itemize}
Later, when we refer to sets 
of points  representing the interface iso-lines 
we use the notation 
$\bx_i^\omega$ where 
$\omega=\alpha,\,\psizw$ or $\psi_0$.
Next, the value
of the curvature 
$\kappa_i=1/\mathcal{R}_i=1/\max{\lr |\bx_i^\omega - \bx_0| \rr}$ 
is computed,
using equations 
(\ref{eq32}) and (\ref{eq33})
away from the interface $\Gamma$.
Now, from a given 
$\kappa_i'$ field
that is a numerical approximation
of the curvature $\kappa_i=1/\mathcal{R}_i$,
the iso-contours 
given by the sets of points $\bx_i^{\omega,\kappa'}$ 
are determined on each grid $m_i$
and for each $\omega$.
In figure \ref{fig19},
approximation of 
the curvature field $\kappa_3'$
obtained on the grid $m_3$ with
$\eph=\Delta x/2$ is presented.

For each grid $m_i$
the $\bx_i^\omega$ and $\bx_i^{\omega,\kappa'}$
iso-lines are divided into $N_i^s$ sample points;
hence, the error of the 
interface curvature approximation 
is defined as 
\blnm
\be
 L_{1,i}^\kappa = \frac{1}{N_i^s} \sum_{l=1}^{N_i^s} | \br_{i,l}^{\kappa'} - \br_{i,l}^\Gamma|,
\label{eq34}
\ee
\elnm
where $\br_{i,l}^{\kappa'}= |\bx_{i,l}^{\omega,\kappa'}-\bx_0|$ 
and   $\br_{i,l}^\Gamma= |\bx_{i,l}^\omega-\bx_0|$,
$\bx_0$ is the center of the circular interface.
Such formula for the error
accounts only for the error 
in the curvature approximation $\kappa_i'$, 
and distinguishes it from the error in 
the interface $\Gamma$
position approximation.

In tables \ref{tab1} and \ref{tab2},
results of the curvature $\kappa_i'$ 
convergence study on five grids $m_i$
and after $N_\tau=256$
re-initialization
steps are given,
for the two interface widths
$\eph=\Delta x/2$
and $\eph=\Delta x \sqrt{2}/4$,
respectively.
%
%
\begin{table}[!ht]
\centering 
\begin{tabular*}{0.95\textwidth}{@{}  l c c  c c  c c @{}} 
\cline{1-7}
 m.f. & \multicolumn{2}{c}{$\psigf$} & \multicolumn{2}{c}{$\psigs$} & \multicolumn{2}{c}{$\psix$} \\ [0.5ex] 
\cline{1-7} 
 $\Gamma$ & $\alpha$ & $\psip$ & $\alpha$ & $\psip$ & $\alpha$ & $\psio$ \\ [0.5ex] 
\cline{1-7} 
$m_1$ & 4.5127e-3  & 3.7965e-3  & 4.0082e-3  &  2.9365e-3   &  4.0082e-3  &  2.9365e-3   \\ 
$m_2$ & 1.8758e-3  & 6.9689e-4  & 1.9332e-3  &  2.1592e-4   &  1.9333e-3  &  2.1591e-4    \\
$m_3$ & 9.6959e-4  & 7.2554e-4  & 8.9694e-4  &  1.4100e-4   &  8.9634e-4  &  1.4258e-4   \\
$m_4$ & 8.1422e-4  & 8.0377e-4  & 2.6027e-4  &  4.1669e-5   &  2.6027e-4  &  4.1669e-5   \\  
$m_5$ & 3.2637e-2  & 3.2472e-2  & 1.4827e-4  &  1.0275e-5   &  1.4827e-4  &  1.0275e-5   \\ [0.5ex] 
\cline{1-7} 
$m_1$ & 5.3185e-3  & 5.2789e-3  & 4.9513e-3  &  4.9016e-3   &  4.9513e-3  &  4.9016e-3   \\ 
$m_2$ & 5.0811e-4  & 5.8833e-4  & 3.0678e-4  &  1.5043e-4   &  3.0676e-4  &  1.5044e-4   \\
$m_3$ & 5.2778e-4  & 5.8279e-4  & 1.9617e-4  &  1.3111e-4   &  1.9617e-4  &  1.3111e-4   \\ 
$m_4$ & 7.7801e-4  & 7.3043e-4  & 1.0271e-4  &  3.7678e-5   &  1.0271e-4  &  3.7678e-5   \\  
$m_5$ & 2.9551e-2  & 2.9521e-2  & 4.0428e-5  &  1.0021e-5   &  4.0428e-5  &  1.0021e-5   \\ [0.5ex] 
\cline{1-7} 
$m_1$ & 5.9374e-3  & 6.3264e-3  & 5.5758e-3  &  5.8275e-3   &  5.5758e-3   & 5.8275e-3   \\ 
$m_2$ & 1.7941e-3  & 8.8314e-4  & 1.2051e-3  &  1.3091e-4   &  1.2051e-3  &  1.3091e-4   \\
$m_3$ & 1.1891e-3  & 7.3411e-4  & 6.6088e-4  &  1.0361e-4   &  6.6330e-4  &  1.0460e-4   \\ 
$m_4$ & 1.0401e-3  & 7.8492e-4  & 4.6755e-4  &  3.6848e-5   &  4.6755e-4  &  3.6848e-5   \\ 
$m_5$ & 2.6359e-2  & 2.6527e-2  & 2.3224e-4  &  9.6314e-6   &  2.3224e-4  &  9.6314e-6   \\ [0.5ex]   
\cline{1-7}  
\end{tabular*}
\caption{Convergence 
         of the curvature 
         $\kappa_i'$ at $\alpx = \lr 0.05,0.5,0.95 \rr$
         iso-lines (from top to bottom)
         after $N_\tau=256$ re-initialization steps, 
         the interface width $\eph=\Delta x/2$.
         Errors are defined using $L_{1,i}^\kappa$ norm
         given by \Eq{eq34}. 
         The results are calculated on five
         grids $m_i$ 
         with three mapping
         functions (m.f.)
         $\psigf$,
         $\psigs$,  
         $\psix$
         and two interface representations 
         $\Gamma \lr \alpps \rr$ 
         and $\Gamma \lr \psio \approx \psip \rr$.}
\label{tab1} 
\end{table}
%
%
\begin{table}[!ht]
\centering 
\begin{tabular*}{0.95\textwidth}{@{}  l c c  c c  c c @{}} 
\cline{1-7}
 m.f. & \multicolumn{2}{c}{$\psigf$} & \multicolumn{2}{c}{$\psigs$} & \multicolumn{2}{c}{$\psix$} \\ [0.5ex] 
\cline{1-7} 
 $\Gamma$ & $\alpha$ & $\psip$ & $\alpha$ & $\psip$ & $\alpha$ & $\psio$ \\ [0.5ex] 
\cline{1-7} 
$m_1$ &  3.3388e-3  &  2.2473e-3  &  3.7087e-3  &  4.0166e-4  &  3.7087e-3  &  4.0166e-4  \\ 
$m_2$ &  1.7094e-3  &  1.2542e-3  &  2.5111e-3  &  5.2983e-4  &  2.5111e-3  &  5.2979e-4  \\
$m_3$ &  1.5087e-3  &  1.7981e-3  &  1.2191e-3  &  1.6811e-4  &  1.2191e-3  &  1.6811e-4  \\
$m_4$ &  2.3557e-3  &  2.4407e-3  &  3.9579e-4  &  4.0151e-5  &  3.9579e-4  &  4.0151e-5  \\
$m_5$ &  3.1094e-3  &  3.0644e-3  &  1.8679e-4  &  1.0546e-5  &  1.8679e-4  &  1.0546e-5  \\ [0.5ex] 
\cline{1-7} 
\cline{1-7} 
$m_1$ & 3.5091e-3  &  2.6268e-3  & 2.3302e-3  & 1.4589e-3 & 2.3302e-3  & 1.4589e-3  \\ 
$m_2$ & 8.0637e-4  &  1.0710e-3  & 7.4319e-4  & 3.5706e-4 & 7.4319e-4  & 3.5708e-4  \\
$m_3$ & 1.4803e-3  &  1.6132e-3  & 3.2301e-4  & 1.4815e-4 & 3.2301e-4  & 1.4815e-4  \\ 
$m_4$ & 2.8712e-3  &  2.6359e-3  & 1.8123e-4  & 3.8763e-5 & 1.8123e-4  & 3.8763e-5  \\  
$m_5$ & 2.5595e-3  &  2.5462e-3 &  7.2441e-5  & 1.0676e-5 & 7.2441e-5  & 1.0676e-5  \\ [0.5ex] 
\cline{1-7} 
$m_1$ &  6.8061e-3  & 3.8563e-3  &  5.4077e-3  & 2.9623e-3  &  5.4077e-3  & 2.9623e-3  \\ 
$m_2$ &  2.6807e-3  & 1.3214e-3  &  1.7051e-3  & 2.7312e-4  &  1.7051e-3  & 2.7313e-4  \\
$m_3$ &  1.8291e-3  & 1.5433e-3  &  7.0173e-4  & 1.3536e-4  &  7.0173e-4  & 1.3536e-4  \\ 
$m_4$ &  2.5678e-3  & 2.3315e-3  &  6.8413e-4  & 3.7086e-5  &  6.8413e-4  & 3.7086e-5  \\ 
$m_5$ &  3.0363e-3  &  2.8712e-3  &  3.1737e-4  &  1.0543e-5  &  3.1737e-4  &  1.0543e-5  \\ [0.5ex] 
\cline{1-7}  
\end{tabular*}
\caption{Convergence 
         of the curvature 
         $\kappa_i'$ at $\alpx = \lr 0.05,0.5,0.95 \rr$
         iso-lines (from top to bottom) 
         after $N_\tau=256$ re-initialization steps, 
         the interface width $\eph=\Delta x \sqrt{2}/4$.
         Errors are defined using $L_{1,i}^\kappa$ norm
         given by \Eq{eq34}. 
         The results are calculated on five
         grids $m_i$ 
         with three mapping
         functions (m.f.) 
         $\psigf$,
         $\psigs$,  
         $\psix$
         and two interface representations 
         $\Gamma \lr \alpps \rr$ 
         and $\Gamma \lr \psio \approx \psip \rr$.}
\label{tab2} 
\end{table}
\begin{table}[ht]
\centering 
\begin{tabular}{@{}  l c c  c c  c c @{}} 
\cline{1-7}
 m.f. & \multicolumn{2}{c}{$\psigf$} & \multicolumn{2}{c}{$\psigs$} & \multicolumn{2}{c}{$\psix$} \\ [0.5ex] 
\cline{1-7} 
 $\Gamma$ & $\alpha$ & $\psi_0'(\psi_1)$ & $\alpha$ & $\psi_0'(\psi_1)$ & $\alpha$ & $\psi_0$ \\ [0.5ex] 
\cline{1-7} 
  $m_1$   & 5.2562e-3  & 5.1339e-3  & 4.8451e-3  &  4.5552e-3   &  4.8451e-3  &  4.5552e-3  \\ 
  $m_2$   & 1.3927e-3  & 7.2278e-4  & 1.1484e-3  &  1.6575e-4   &  1.1484e-3  &  1.6576e-4  \\
  $m_3$   & 8.9547e-4  & 6.8082e-4  & 5.8466e-4  &  1.2524e-4   &  5.8531e-4  &  1.2613e-4  \\ 
  $m_4$   & 8.7741e-4  & 7.7305e-4  & 2.7684e-4  &  3.8732e-5   &  2.7684e-4  &  3.8732e-5  \\ 
  $m_5$   & 2.9515e-2  & 2.9506e-2  & 1.4031e-4  &  9.9756e-6   &  1.4031e-4  &  9.9756e-6  \\ [0.5ex] 
\cline{1-7}  
\end{tabular}
\caption{Convergence 
         of the curvature $\kappa_i'$ in the narrow band
         of $0.05 \le \alpha \lr \bx_\Gamma,t \rr \le 0.95$
         the interface width  $\eph=\Delta x/2$. 
         The table contains arithmetical mean of the
         values in appropriate columns and rows 
         of \Tab{tab1}.}
\label{tab3} 
\end{table}
\begin{table}[ht]
\centering 
\begin{tabular}{@{}  l c c  c c  c c @{}} 
\cline{1-7}
 m.f.     & \multicolumn{2}{c}{$\psigf$} & \multicolumn{2}{c}{$\psigs$} & \multicolumn{2}{c}{$\psix$} \\ [0.5ex] 
\cline{1-7} 
 $\Gamma$ & $\alpha$ & $\psi_0'(\psi_1)$ & $\alpha$ & $\psi_0'(\psi_1)$ & $\alpha$ & $\psi_0$ \\ [0.5ex] 
\cline{1-7} 
 $m_1$    &  4.5513e-3  &  2.9101e-3  & 3.8155e-3  &  1.6076e-3  & 3.8155e-3  &  1.6076e-3  \\
 $m_2$    &  1.7322e-3  &  1.2156e-3  & 1.6531e-3  &  3.8667e-4  & 1.6531e-3  &  3.8666e-4  \\
 $m_3$    &  1.6061e-3  &  1.6515e-3  & 7.4794e-4  &  1.5054e-4  & 7.4794e-4  &  1.5054e-4  \\
 $m_4$    &  2.5982e-3  &  2.4694e-3  & 4.2038e-4  &  3.8666e-5  & 4.2038e-4  &  3.8666e-5 \\
 $m_5$    &  2.9017e-3  &  2.8272e-3  & 1.9220e-4  &  1.0588e-5  & 1.9220e-4  &  1.0588e-5   \\ [0.5ex] 
\cline{1-7}  
\end{tabular}
\caption{Convergence 
         of the curvature $\kappa_i'$ in the narrow band
         of $0.05 \le \alpha \lr \bx_\Gamma,t \rr \le 0.95$
         the interface width   $\eph= \sqrt{2}\Delta x/4$. 
         The table contains arithmetical mean of the
         values in appropriate columns and rows 
         of \Tab{tab2}.}
\label{tab4} 
\end{table}

\begin{figure}[h!] \nonumber
 \includegraphics[angle=-90]{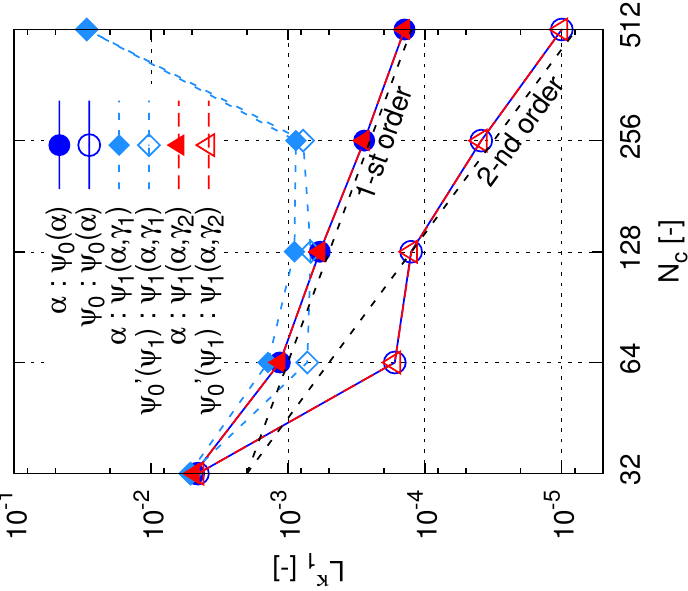}
 \includegraphics[angle=-90]{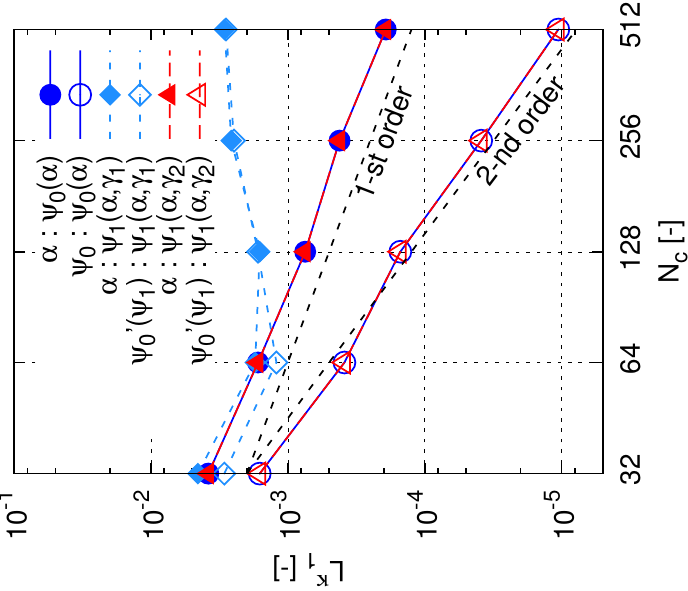}
 \caption{\small{Convergence of the 2D circular interface curvature $\kappa_i'$
                 in the narrow band of the level-set function 
                 $0.05 \le \alpx \le 0.95$ with 
                 (a) $\eph=\Delta x/2$  and
                 (b) $\eph=\Delta x \sqrt{2}/4$ 
                 after $N_\tau = 256$ 
                 re-initialization steps.
                 In the legend $\alpha : \psi_0 \lr \alpha \rr$,
                 the symbol on the left denotes the interface
                 representation, the symbol on the right denotes
                 the mapping function used to compute 
                 derivatives in \Eq{eq31k}.
                 $N_c$ is the number of grid points in $x,y$
                 directions.}}
 \label{fig21}
 \end{figure}

In tables \ref{tab3} and \ref{tab4},
the errors from tables \ref{tab1}-\ref{tab2}
averaged in the narrow band $(0.05,0.5,0.95)$ 
are presented. 
The results  
in tables \ref{tab3} and \ref{tab4}
are depicted in \Fig{fig21}, 
and they can be interpreted as 
the convergence rate 
of the circular interface
curvature in the narrow band 
$0.05 \le \alpx \le 0.95$
with two different
$\deltaa/\eph$ support widths
$\eph=\Delta x/2$ 
or $\eph=\Delta x \sqrt{2}/4$.

In figure \ref{fig21},
one observes
that the
second order 
convergence rate
of the curvature 
is achieved
for the
mapping functions $\psigs$
and $\psio$ 
when the signed
distance function 
interface
representations 
given by $\psio \approx \psizw \lr \psigs \rr$ 
are chosen.
The convergence 
rates of the 
interface 
curvature 
differ
dependent on whether 
the interface
is represented by 
$\alpps$ or $\psio$.
This latter observation
is explained by the fact
that to reconstruct the jump
at $\Gamma$ 
the level-set
function $\alpp$ 
requires the 
constant number 
of grid points 
(from five to three 
dependent on selected $\eph$,
see \Fig{fig7}) 
regardless of 
the grid resolution
used in 
simulation.
At the same time,
the accuracy of the representation
of the signed distance function
$\psio$, or its approximation
$\psi_0' \lr \psigs \rr$, 
increases proportionally to
the number of grid points $N_c$.

In the case of the 
mapping function
$\psigf$, convergence 
of the interface curvature 
is not achieved
although equation (\ref{eq1}) 
is solved 
to a fully convergent state 
on all grids $m_i$, $i=1,\dots,5$,
see \Fig{fig18}.
This result may be explained
by the existence of 
a non-zero second order derivative 
of $\psigf$ 
function 
from the
interface $\bx_\Gamma$
that impairs 
the calculation
of the 
interface curvature, 
see \Fig{fig11}(b). 

For both interface widths
$\eph=\Delta x/2$
or $\eph=\sqrt{2}\Delta x/4$, 
the second order convergence rate
of $\kappa_i'$ is obtained
for $\psio$ and $\psigs$,
see \Fig{fig21}. 
This confirms
discussion 
regarding the allowable 
width of the $\deltaa/\eph$
support presented in Section \ref{sec5.2}.
The two-three grid points 
(see \Fig{fig7} and \Fig{fig14}) 
required to reconstruct the interface 
curvature $\kappa$ with 
second-order accuracy
are also needed in the construction 
of the flux limiters 
during a solution of the advection 
equation (\ref{eq30a});
the two-three grid points 
is a typical resolution
of the VOF interface capturing methods 
\cite{ubbink99,waclawczyk06,waclawczyk08_2, waclawczyk08_3}.

\subsubsection{Wavy interface}
\label{sec6.1.3}

McCaslin and Desjardins \cite{mccaslin2014} 
noticed 
the feedback mechanism 
between 
erroneous 
normals $\bng$
and the solution 
of the re-initialization equation (\ref{eq1}) 
when  $N_\tau \gg 1$.
This numerical phenomenon 
occurs due to errors introduced
during discretization 
of equation (\ref{eq1}) 
and leads to artificial deformations 
of the level-set function $\alpp$ 
increasing with time $\tau$.
According to  \cite{mccaslin2014},
this defect in 
the re-initialization procedure 
is particularly noticeable 
in regions where the
interface $\Gamma$ 
is stationary.
In order
to reduce 
this erroneous feedback between
$\bng$ and $\alpp$ 
in  \cite{mccaslin2014}  
it is proposed to vary 
the amount of re-initialization
spatially and to localize 
the solution of \Eq{eq1}
only in the regions where
the interface is advected 
or $|\bng \cdot \bu| \gg 1$.
With this method, 
an additional function 
$\beta \lr \bx,t\rr$, 
variable in time
and space is introduced.
The diffusive 
and compressive
fluxes on the RHS of equation (\ref{eq1})
are  multiplied by $\beta \lr \bx,t\rr$
to localize 
re-initialization
depending on the local
flow conditions.
The function 
$\beta \lr \bx,t \rr$
has to satisfy
the condition
$\bng \cdot \nabla \beta = 0$ 
used to
introduce
an additional 
equation for 
$\beta \lr \bx,t \rr$.
We recognize 
that
in the method
put forward by
the present paper
$\delta \lr \alpha \rr$
given by \Eq{eq24}
(see also \Eqs{eq29}{eq30}), 
plays a role 
similar to the function
$\beta \lr \bx,t\rr$
in \cite{mccaslin2014}.
However, in the present work,
we do not need to vary the number 
of re-initialization steps
in time and space relative to
local flow conditions.

In this section,
we carry out
re-initialization of 
the 3D wavy interface
similar to the test case
proposed in  \cite{mccaslin2014}.
The initial condition to \Eq{eq1}
is given by  \Eq{eq2} where the 
signed distance function 
\blnm
\be
 \psix = 1/2-y+A_0 \sin{\lr 4\pi x\rr \sin{\lr 4\pi z\rr}},
\label{eq35}
\ee
\elnm
defines the wavy interface,
$A_0=0.03125\,m$ 
on all grids $m_i$.
Substitution of \Eq{eq35} into
\Eq{eq2} allows us to compute
the exact curvature $\kappa$ 
of the interface $\Gamma$.
Convergence 
of $\kappa_i'$, 
which is 
a numerical 
approximation 
of the exact interface 
curvature $\kappa_i$, 
is studied
on the four gradually
refined meshes: 
$m_i=2^{4+i} \times 2^{4+i} \times 2^{4+i}$
where $i=2,\dots,4$
and  
a mesh 
$m_{3.5}= 192 \times 192 \times 192$
with the number of grid nodes
in-between $m_3$, $m_4$.
The  
computational
domain  is
a cubic box
$\Omega = <\!-0.5,0.5\!> 
  \times <\!0,1\!> 
  \times <\!-0.5,0.5\!>\,m^3$
discretized
using a uniform 
grid nodes 
distribution,
the time step size
is set to 
$\Delta \tau =D/C^2=\eph$.

Using experience  
gained in the 
previous 2D test cases 
(see \Sec{sec6.1.2}),
two mapping functions
$\psio$ and $\psigs$ 
are used for 
discretization
of $|\nabla \alpha|$
in \Eq{eq1}.
Similar to  
the previous example,
we focus on convergence 
of the interface 
curvature in 
the narrow band 
$0.05 \le \alpha \le 0.95$.
Additionally, 
as in Section \ref{sec6.2}, 
the influence of the $\deltaa/\eph$ support width
$\eph=\Delta x/2$ 
or $\eph=\sqrt{3} \Delta x/4$
on convergence rate of $\kappa_i'$ 
is studied.
To asses effects  of
the $\alpha-\bng$ coupling
on the interface deformations,
we consider three test cases
in which $\bn_\Gamma$ is constant
or variable in time $\tau$:
\begin{itemize}
\item[$T_1\,:$]{$\bng=const$, $\nabla \alpha$  is calculated before first iteration,
                $N_\tau=1024$, $\eph=\Delta x/2$.}
\item[$T_2\,:$]{$\bng \ne const$, 
                $\nabla \alpha$  is calculated before each  iteration, $N_\tau=1024$, $\eph=\Delta x/2$.}
\item[$T_3\,:$]{as the case $T_2$ but with $\eph=\sqrt{3}\Delta x/4$.}

\end{itemize}
%

\paragraph{Computation of numerical errors in the wavy interface curvature}
\label{sec6.1.3.1}

%
In the beginning,
we note that 
the curvature 
$\kappa$ derived from
the analytical solution given
by \Eq{eq2} and \Eq{eq35}
is exact only at the interface 
$\Gamma$.
The comparison with
$\kappa'$ that is calculated
away from the interface $\Gamma$
may lead to incorrect 
interpretations of the 
convergence results 
in the narrow band 
of the level-set function 
$0.05 \le \alpha \le 0.95$.
This is due to 
the 3D wavy interface curvature
$\kappa_i \ne const$,
in the direction normal 
and tangential to the interface $\Gamma$.
For this reason, computation
of the errors in the numerical 
approximation $\kappa_i'$ 
of the exact curvature $\kappa_i$,
is performed two ways.
In the first approach, 
the error of the numerical solution
is determined directly from the difference
between 
$\kappa_i$  
and 
$\kappa_i'$, 
which are 
both computed 
in the centers of the control
volumes $N_i^p$ localized 
in the narrow band
$0.05 \le \alpha \le 0.95$ 
on the mesh $m_i$.
The $L_1^\kappa$, $L_2^\kappa$ and $L_\infty^\kappa$ 
norms on the mesh $m_i$ 
are defined by following
formulas  
\blnm
\be
 L_{1,i}^\kappa = \frac{1}{N_i^p} \sum_{l=1}^{N_i^p} | \kappa_l - \kappa_l'  |,
\label{eq36}
\ee
\elnm
\blnm
\be
 L_{2,i}^\kappa =  \frac{1}{N_i^p} \ls \sum_{l=1}^{N_i^p} \lr \kappa_l - \kappa_l'  \rr^2 \rs^{1/2},
\label{eq37}
\ee
\elnm
\blnm
\be
  L_{\infty,i}^\kappa = \max{\lr | \kappa_l - \kappa_l'  | \rr}, \quad l=1,\dots,N_i^p.
 \label{eq38}
\ee
\elnm

In the second approach,
the convergence rate
of $\kappa_i'$ is determined 
dependent on the interface 
representation:
by the conservative 
level-set function
$\Gamma \lr \alpha \lr \bx_i,t \rr \rr$
or by the signed distance function
$\Gamma \lr \psio \rr$.
Since in the present case
the curvature $\kappa_i$
is variable in space and the relation between
$\kappa_i'$ and $\bx_{i,\Gamma}$ is not known, 
it is very difficult to apply 
the procedure presented in Section \ref{sec6.1.2}
for computation of $\kappa_i'$.
For this reason, we use a simplified 
approach using available post-processing
tools.
First, we calculate 
$\kappa_i$ and 
$\kappa_i'$
in the centers 
of the control volumes,
then the norm  
$L_i^\kappa = \kappa_i - \kappa_i'$
is evaluated.
Afterwards, 
$L_i^\kappa$
is interpolated 
to the interface 
represented by the conservative level-set function
$\Gamma \lr \alpha \lr \bx_i,t \rr \rr$
or the signed distance function
$\Gamma \lr \psio \rr$
obtained as iso-surfaces
in the post-processing
software.
The maximal 
value of $|L_i^\kappa \lr \Gamma \lr \alpha \rr \rr|$ 
or $|L_i^\kappa \lr \Gamma \lr \psi_0 \rr \rr|$ 
on grid $m_i$
allows us to estimate
$L_{\infty,i}^\kappa \lr \alpha \rr$
and
$L_{\infty,i}^\kappa \lr \psi_0 \rr$
norms at $\Gamma \lr \alpha \rr$
and $\Gamma \lr \psi_0 \rr$, 
respectively.
%

\paragraph{Convergence of the re-initialization equation}
\label{sec6.1.3.2}

\begin{figure}[h!] \nonumber
 \begin{centering}
 \begin{minipage}{.5\textwidth}
  \centering
  \includegraphics[width=.666\textwidth,height=1.\textwidth,angle=-90]{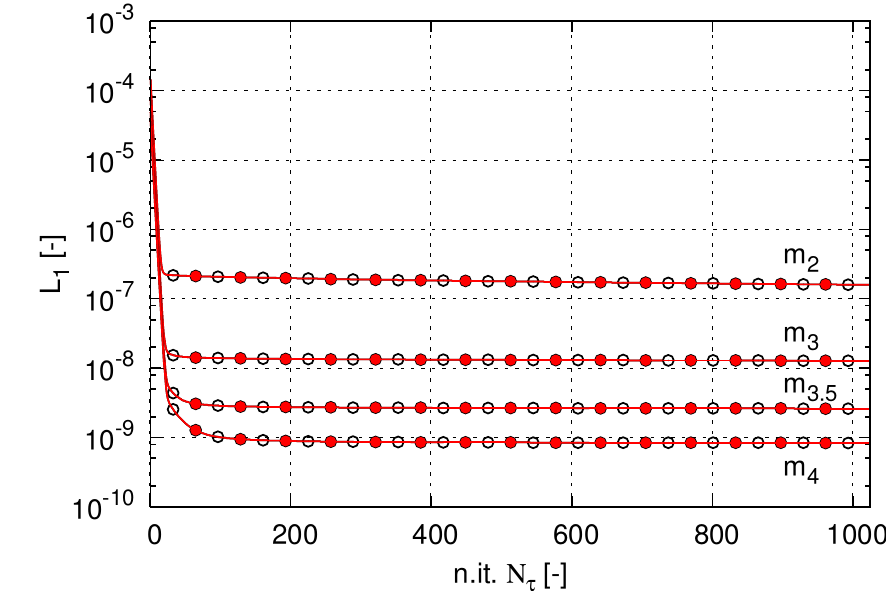}
 \end{minipage}
 \begin{minipage}{.5\textwidth}
  \centering
  \includegraphics[width=.666\textwidth,height=1.\textwidth,angle=-90]{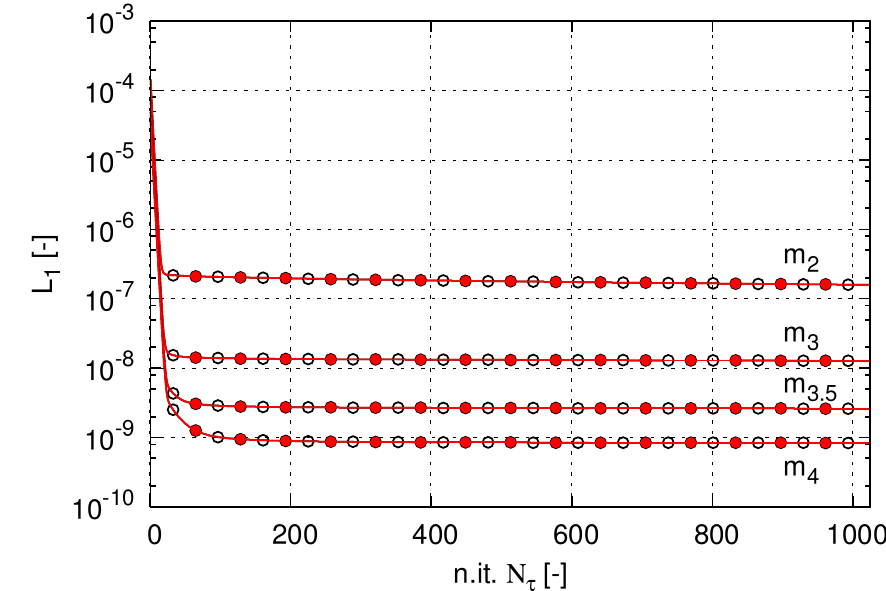}
 \end{minipage}
 \hspace*{2.5cm}
 \begin{minipage}{.5\textwidth}
  \centering
  \includegraphics[width=.666\textwidth,height=1.\textwidth,angle=-90]{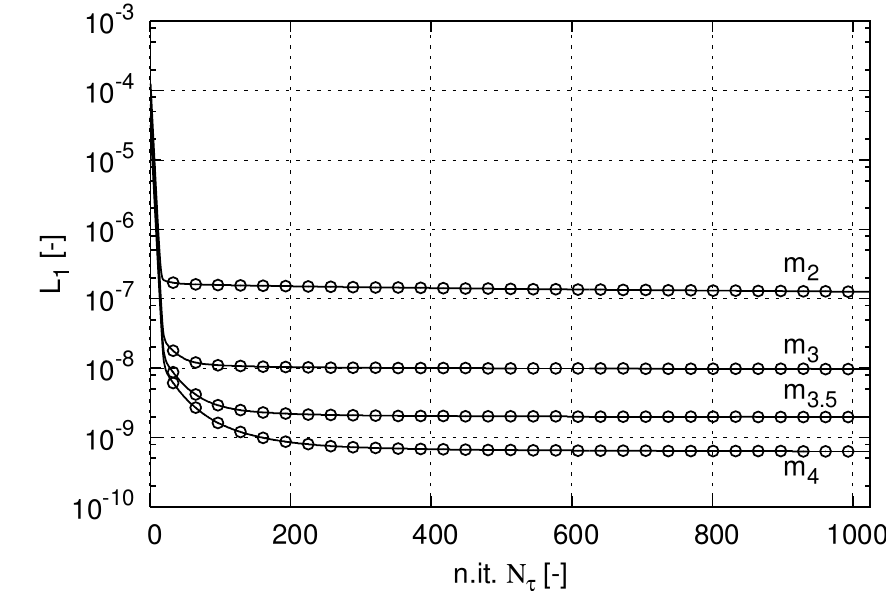}
 \end{minipage}
 \end{centering}
 \caption{\small{Convergence of the solution to re-initialization 
                 \Eq{eq1}
                 in the test cases
                 (a) $T_1$, (b) $T_2$, (c) $T_3$.
                 $|\nabla \alpha|$ is
                 discretized using the mapping functions 
                 $\psio$ (black void dots)
                 $\psi_0' \lr \psigs \rr$ (red solid dots), 
                 the results were obtained
                 on four gradually refined grids $m_i$, 
                 $L_1$ norm is defined by \Eq{eq9}.}}
 \label{fig22}
 \end{figure}

In figure \ref{fig22} 
convergence of the solution
to \Eq{eq1} in the test cases 
$T_1$, $T_2$ and $T_3$ 
with two mapping functions 
$\psio$ and $\psigs$
is presented.
We observe that 
to obtain the stationary solution
(the constant convergence level) 
about $N_\tau \approx 100$ 
iterations are needed in the cases $T_1$, $T_2$
and about $N_\tau \approx 200$ iterations 
in the test case $T_3$ (here only  $\psio$ is used).

We note that at 
least second-order 
convergence rate 
of $L_1$ norm
(defined by \Eq{eq9})
in time $\tau$
may be deduced from
\Fig{fig22} in the
all three test
cases.
Surprisingly, 
the rate of convergence
of the solution to \Eq{eq1} 
is exactly the same in the
test cases $T_1$ and $T_2$, 
see \Fig{fig22}(a)(b).
Since the interface curvature
$\kappa$ is implicitly included
in the re-initialization equation (\ref{eq1}) 
(see \Eqs{eq29}{eq30}),
previously discussed numerical 
errors introduced 
by calculation of gradients in 2D-3D 
space have a larger impact on the
convergence rate 
than the errors in  computation
of $\kappa$.

The convergence rate in the test
case $T_3$ presented in \Fig{fig22}(c) 
is lower than in previously
discussed tests $T_1$, $T_2$.
This result is to some degree in 
opposition 
to the convergence studies presented in figures 
\ref{fig6} and \ref{fig18}.
However,  the wavy interface curvature   
is variable not only in the direction normal 
but also in the direction tangential to 
the interface $\Gamma$.
For this reason
the reduced number of grid points 
when the interface width is set to
$\eph=\sqrt{3}\Delta x/4 < \Delta x/2$ 
may lead to slower convergence.

\paragraph{Errors in computations of the wavy interface curvature}
\label{sec6.1.3.3}

During numerical experiment
$T_1$ it was found that 
values of the errors 
defined by equations
(\ref{eq36})-(\ref{eq38})
remain constant
and equal to error
after $N_\tau=1$
re-initialization
steps,
see table \ref{tab5}.
Similar to the convergence rates 
presented in \Fig{fig22},
the values of the norms defined by
\Eqs{eq36}{eq38}  in table \ref{tab5} 
are identical for the mapping 
functions $\psigs$ and $\psio$.
This result is consistent
with data presented in tables 
\ref{tab1}-\ref{tab2}
and is expected in the light of 
equation (\ref{eq19}).
\begin{table}[ht]
\centering 
\begin{tabular}{@{}  l  c c c  c c c @{}} 
\cline{1-7}
 m.f.    & \multicolumn{3}{c}{$\psigs$} & \multicolumn{3}{c}{$\psio$} \\ [0.5ex] 
\cline{1-7} 
 L       & $L_1^\kappa$ & $L_2^\kappa$ & $L_\infty^\kappa$ & $L_1^\kappa$ & $L_2^\kappa$ & $L_\infty^\kappa$ \\ [0.5ex] 
\cline{1-7} 
$m_2$    & 4.6636e-2 & 5.7389e-2 & 1.2425e-1  &  4.6636e-2   &  5.7389e-2  &  1.2425e-1  \\ 
$m_3$    & 1.1639e-2 & 1.4403e-2 & 3.1545e-2  &  1.1639e-2   &  1.4403e-2  &  3.1545e-2  \\ 
$m_{3.5}$ & 5.1776e-3 & 6.4156e-2 & 1.4061e-2  &  5.1776e-3   &  6.4156e-3  &  1.4061e-2  \\
$m_4$    & 2.9106e-3 & 3.6092e-3 & 7.9169e-3  &  2.9106e-3   &  3.6092e-3  &  7.9169e-3  \\ [0.5ex] 
\cline{1-7}  
\end{tabular}
\caption{The $L_1^\kappa$,  $L_2^\kappa$ and $L_\infty^\kappa$ 
         norms obtained using \Eqs{eq36}{eq38} 
         in the test case $T_1$ after $N_\tau=1024$ 
         re-initialization steps.}
\label{tab5} 
\end{table}
\begin{table}[ht]
\centering 
\begin{tabular}{@{}  l  c c c @{}} 
\cline{1-4}
 m.f.    & \multicolumn{3}{c}{$\psio$}  \\ [0.5ex] 
\cline{1-4} 
 L       & $L_1^\kappa$ & $L_2^\kappa$ & $L_\infty^\kappa$ \\ [0.5ex] 
\cline{1-4} 
$m_2$    &  4.6414e-2   &  5.7396e-2   &  1.2425e-1  \\ 
$m_3$    &  1.1657e-2   &  1.4478e-2   &  3.1545e-2  \\ 
$m_{3.5}$ &  5.2041e-3   &  6.4583e-3   &  1.4060e-2  \\
$m_4$    &  2.9263e-3   &  3.6317e-3   &  7.9169e-3  \\ [0.5ex] 
\cline{1-4}  
\end{tabular}
\caption{The $L_1^\kappa$,  $L_2^\kappa$ and $L_\infty^\kappa$ 
         norms obtained using \Eqs{eq36}{eq38} 
         in the test case $T_3$ after $N_\tau=1$ 
         re-initialization steps.}
\label{tab6} 
\end{table}

\begin{figure}[ht!] \nonumber
 \begin{centering}
 \begin{minipage}{.5\textwidth}
  \centering
  \includegraphics[width=.666\textwidth,height=1.\textwidth,angle=-90]{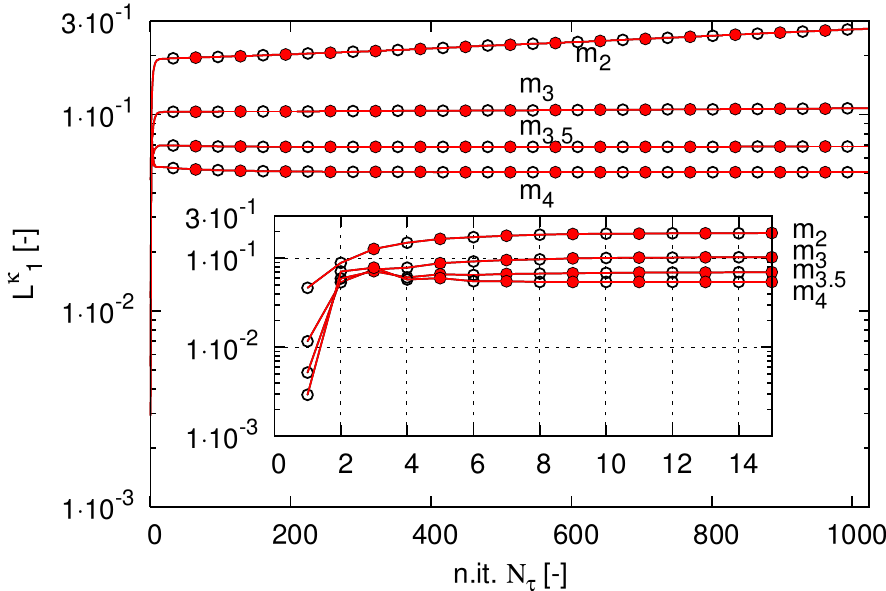}
 \end{minipage}%
 \begin{minipage}{.5\textwidth}
  \centering
  \includegraphics[width=.666\textwidth,height=1.\textwidth,angle=-90]{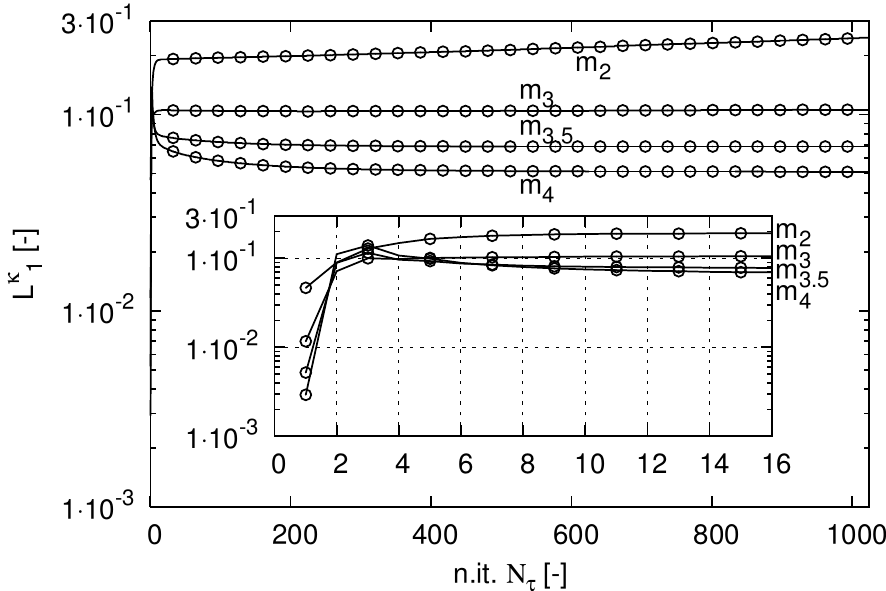}
 \end{minipage}
 \begin{minipage}{.5\textwidth}
  \centering
  \includegraphics[width=.666\textwidth,height=1.\textwidth,angle=-90]{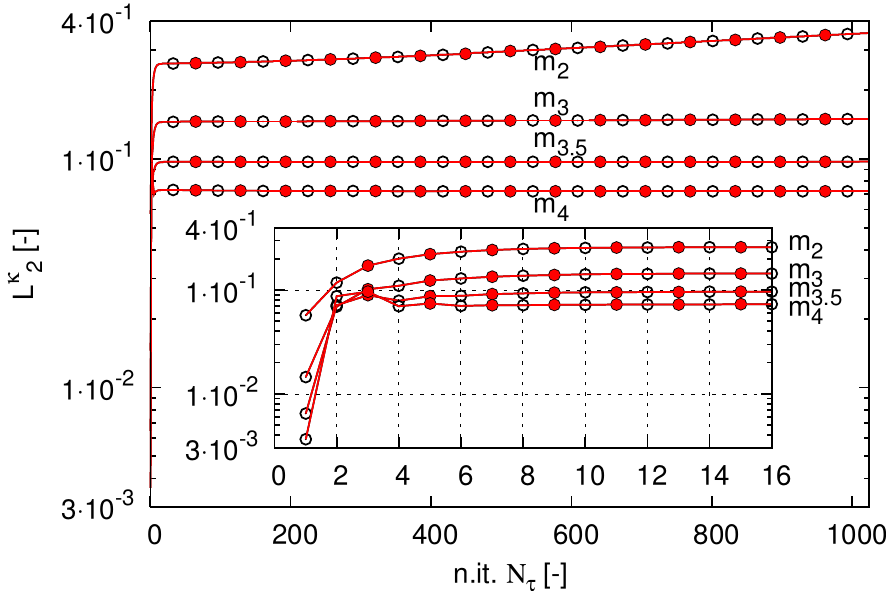}
 \end{minipage}%
 \begin{minipage}{.5\textwidth}
  \centering
  \includegraphics[width=.666\textwidth,height=1.\textwidth,angle=-90]{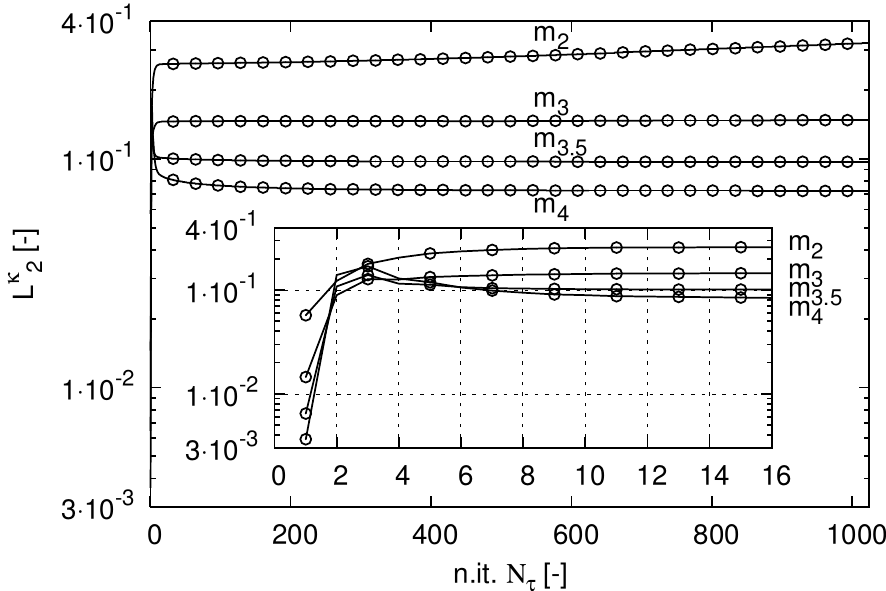}
 \end{minipage}
 \begin{minipage}{.5\textwidth}
  \centering
  \includegraphics[width=.666\textwidth,height=1.\textwidth,angle=-90]{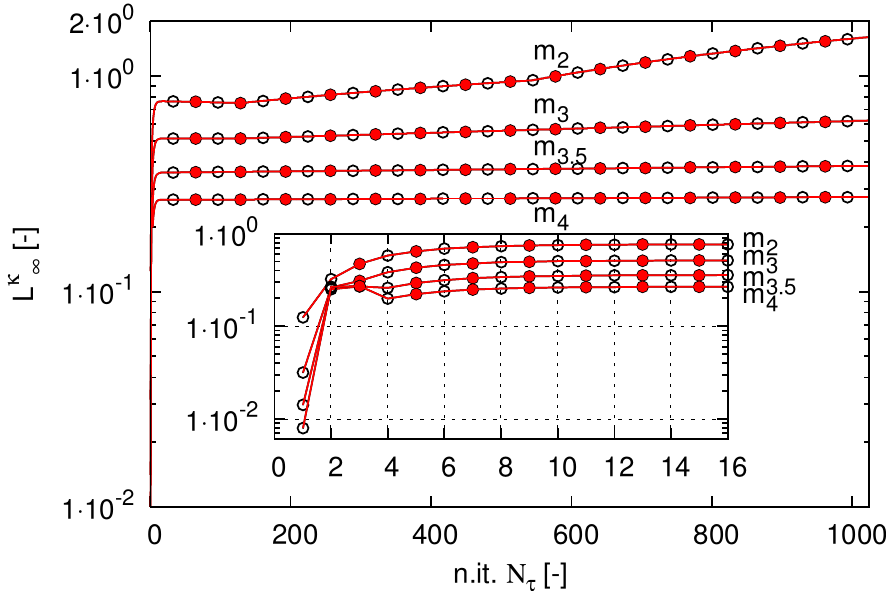}
 \end{minipage}%
 \begin{minipage}{.5\textwidth}
  \centering
  \includegraphics[width=.666\textwidth,height=1.\textwidth,angle=-90]{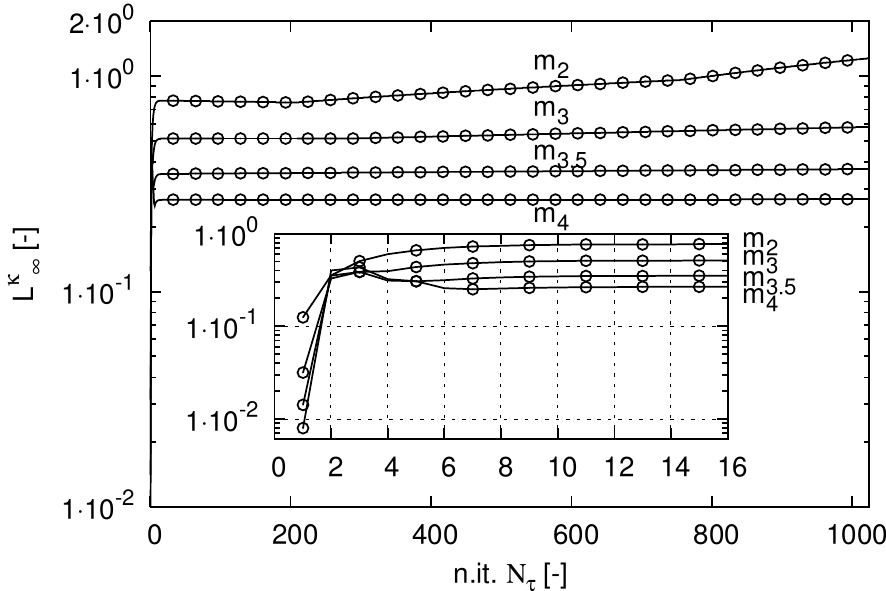}
 \end{minipage}
 \end{centering}
 \caption{\small{Convergence of $L_1^\kappa$, $L_2^\kappa$
                 and $L_\infty^\kappa$ norms during $N_\tau=1024$
                 re-initialization steps of the 3D wavy interface            
                 with \Eq{eq1} in the test cases: $T_2$ (left) and
                 $T_3$ (right).
                 $|\nabla \alpha|$ is discretized
                 using mapping functions $\psio$ (black void dots) 
                 and $\psizw \lr \psigs \rr$ (red solid dots).}}
 \label{fig23}
 \end{figure}

In figure \ref{fig23},
the evolution
of the errors
recorded
during test 
cases $T_1$, $T_2$ and $T_3$ 
are depicted.
We note that after 
the first re-initialization step,
levels of the all errors are
equal to the values 
obtained in the test case
$T_1$ 
(compare results
in \Fig{fig23} 
at $N_\tau=1$
with values 
in tables 
\ref{tab5}-\ref{tab6}).
When $N_\tau > 1$ values of 
the errors grow,
but after about $N_\tau \approx 6$
iterations they reach an almost
constant level which remain
until the end of re-initialization.
Such behavior is observed
in the three finest meshes
$m_2,\,m_{3.5},\,m_4$.
The solution in the mesh $m_2$
is somewhat unresolved, since 
in this case we
have only two grid nodes
per wave amplitude
(in \cite{mccaslin2014} three nodes were 
used).
\begin{table}[ht]
\centering 
\begin{tabular}{@{}  l  c c c  c c c @{}} 
\cline{1-7}
  m.f.   & \multicolumn{3}{c}{$\psigs$} & \multicolumn{3}{c}{$\psio$} \\ [0.5ex] 
\cline{1-7} 
   L     & $L_1^\kappa$ & $L_2^\kappa$ & $L_\infty^\kappa$ & $L_1^\kappa$ & $L_2^\kappa$ & $L_\infty^\kappa$ \\ [0.5ex] 
\cline{1-7} 
$m_2$    & 2.7461e-1 & 3.5492e-1 & 1.51936    &  2.7461e-1 & 3.5492e-1 & 1.51936    \\ 
$m_3$    & 1.0791e-1 & 1.4944e-1 & 6.2188e-1  &  1.0791e-1 & 1.4944e-1 & 6.2188e-1   \\ 
$m_{3.5}$ & 6.8839e-2 & 9.7209e-2 & 3.8353e-1  &  6.8839e-2 & 9.7209e-2 & 3.8353e-1  \\
$m_4$    & 5.0941e-2 & 7.2103e-2 & 2.7458e-1  &  5.0941e-2 & 7.2103e-2 & 2.7458e-1  \\ [0.5ex] 
\cline{1-7}  
\end{tabular}
\caption{The $L_1^\kappa$,  $L_2^\kappa$ and $L_\infty^\kappa$ 
         norms  obtained using \Eqs{eq36}{eq38} 
         in the test case $T_2$ after $N_\tau=1024$ 
         re-initialization steps.}
\label{tab7} 
\end{table}
\begin{table}[ht]
\centering 
\begin{tabular}{@{}  l  c c c @{}} 
\cline{1-4}
 m.f.     & \multicolumn{3}{c}{$\psio$}  \\ [0.5ex] 
\cline{1-4} 
  L    & $L_1^\kappa$ & $L_2^\kappa$ & $L_\infty^\kappa$ \\ [0.5ex] 
\cline{1-4} 
$m_2$    &  2.4721e-1 & 3.2123e-1 & 1.210778   \\ 
$m_3$    &  1.0624e-1 & 1.4785e-1 & 5.8091e-1   \\ 
$m_{3.5}$ &  6.8784e-2 & 9.7189e-2 & 3.7048e-1  \\
$m_4$    &  5.1173e-2 & 7.2384e-2 & 2.6857e-1  \\ [0.5ex] 
\cline{1-4}  
\end{tabular}
\caption{The $L_1^\kappa$,  $L_2^\kappa$ and $L_\infty^\kappa$ 
         norms  obtained using \Eqs{eq36}{eq38} 
         in the test case $T_3$ after $N_\tau=1024$ 
         re-initialization steps.}
\label{tab8} 
\end{table}
The end values 
of the errors
obtained in the test
cases $T_2$, $T_3$ are given 
in tables 
\ref{tab7}-\ref{tab8},
respectively. 

\begin{figure}[ht!] \nonumber
 \begin{centering}
 \begin{minipage}{.333\textwidth}
  \centering
  \includegraphics[width=1.166\textwidth,height=1.\textwidth,angle=-90]{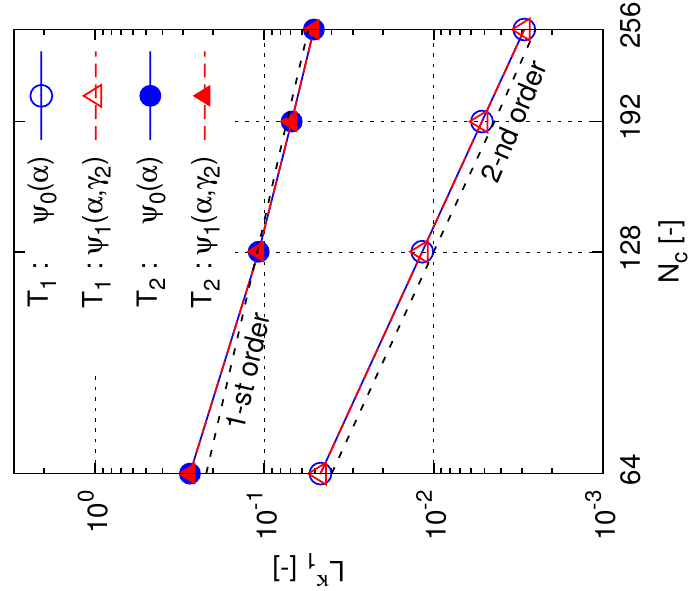}
 \end{minipage}%
 \begin{minipage}{.333\textwidth}
  \centering
  \includegraphics[width=1.166\textwidth,height=1.\textwidth,angle=-90]{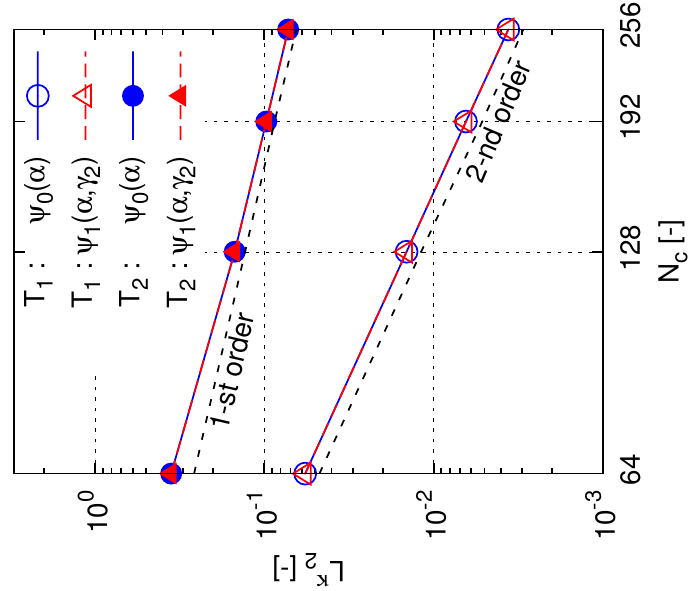}
 \end{minipage}%
 \begin{minipage}{.333\textwidth}
  \centering
  \includegraphics[width=1.166\textwidth,height=1.\textwidth,angle=-90]{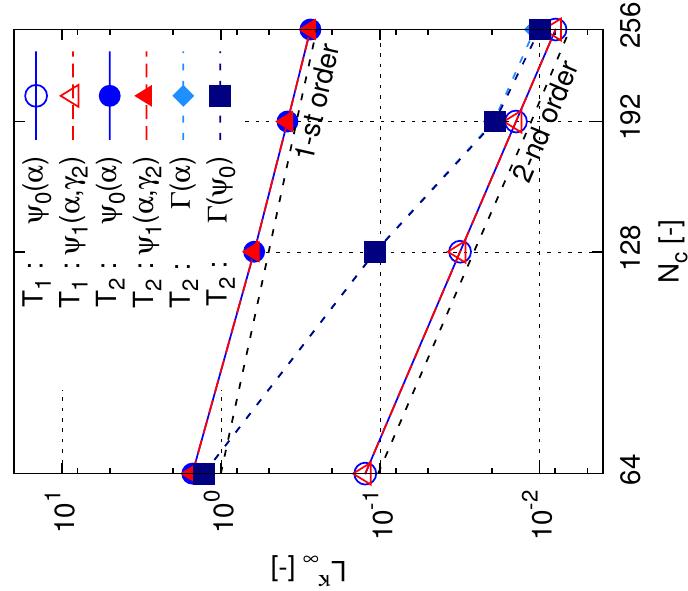}
 \end{minipage}
 \begin{minipage}{.333\textwidth}
  \centering
  \includegraphics[width=1.166\textwidth,height=1.\textwidth,angle=-90]{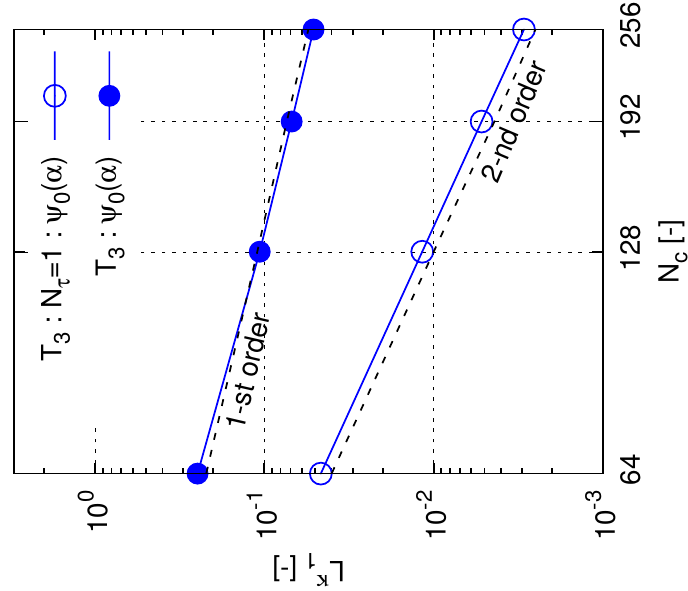}
 \end{minipage}%
 \begin{minipage}{.333\textwidth}
  \centering
  \includegraphics[width=1.166\textwidth,height=1.\textwidth,angle=-90]{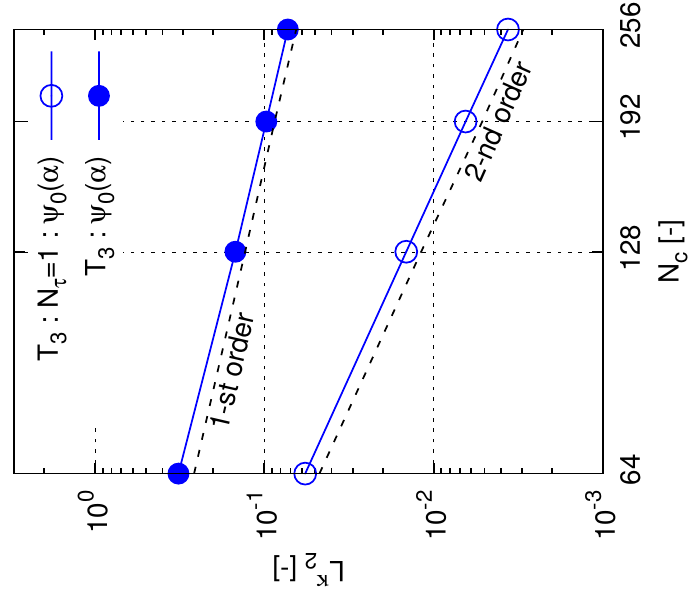}
 \end{minipage}%
 \begin{minipage}{.333\textwidth}
  \centering
  \includegraphics[width=1.166\textwidth,height=1.\textwidth,angle=-90]{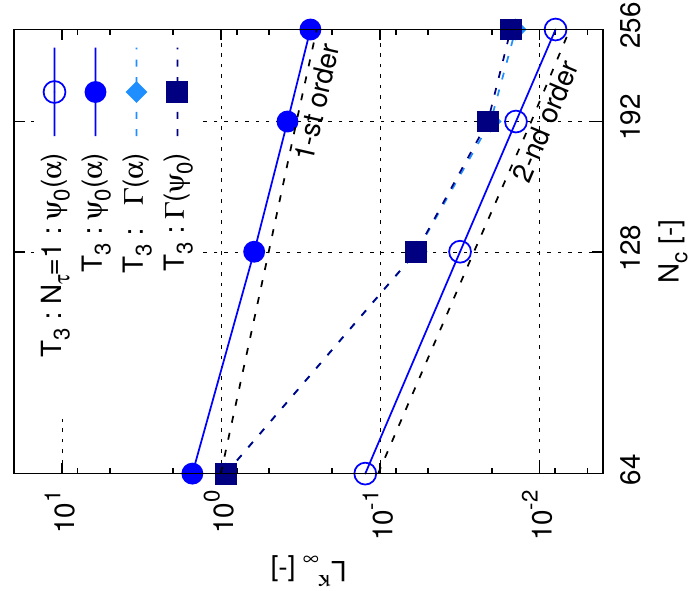}
 \end{minipage}
 \end{centering}
 \caption{\small{The comparison of
                 convergence rates for 
                 $L_1^\kappa$, $L_2^\kappa$
                 and $L_\infty^\kappa$ norms 
                 after $N_\tau=1$  (void symbols)
                 and $N_\tau=1024$ (solid symbols)
                 re-initialization steps
                 in the test cases $T_1$, $T_2$ (top) and $T_3$ (bottom). 
                 $|\nabla \alpha|$ in \Eq{eq1}
                 is discretized
                 using the mapping functions $\psio$ (dots, diamonds, squares) and
                 $\psizw \lr \psigs \rr$ (triangles).}}
 \label{fig24}
 \end{figure}

In figure \ref{fig24},
the results given in tables
\ref{tab5}-\ref{tab8}
are depicted.
We observe that at $N_\tau=1$
(and in the test case $T_1$)
errors defined by 
equations (\ref{eq36})-(\ref{eq38}) 
show the second order
convergence rate. 
We emphasize once again
that in the test case $T_1$,
the values of all errors 
remain constant 
during all 
$N_\tau=1024$  steps.
Hence, for $\bng=const$
in time $\tau$
the second-order
convergence rate 
is also achieved in the
narrow band of $\alpps$.
In the test cases $T_2$ and $T_3$,
the first-order
convergence rate 
of the interface curvature
is detected in the narrow
band of the conservative
level set function 
$0.05 \le \alpps \le 0.95$.
 \begin{figure}[ht!] \nonumber
 \begin{centering}
 \begin{minipage}{.49\textwidth}
  \centering
  \includegraphics[width=0.9\textwidth,height=.5\textwidth,angle=0]{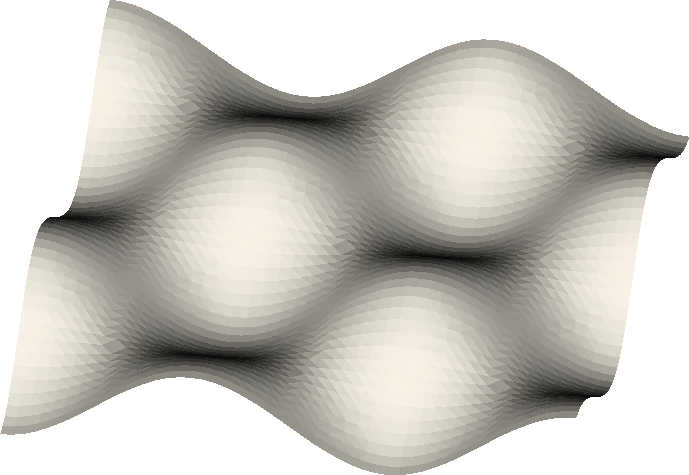}
 \end{minipage}%
 \begin{minipage}{.49\textwidth}
  \centering
  \includegraphics[width=0.8\textwidth,height=.65\textwidth,angle=0]{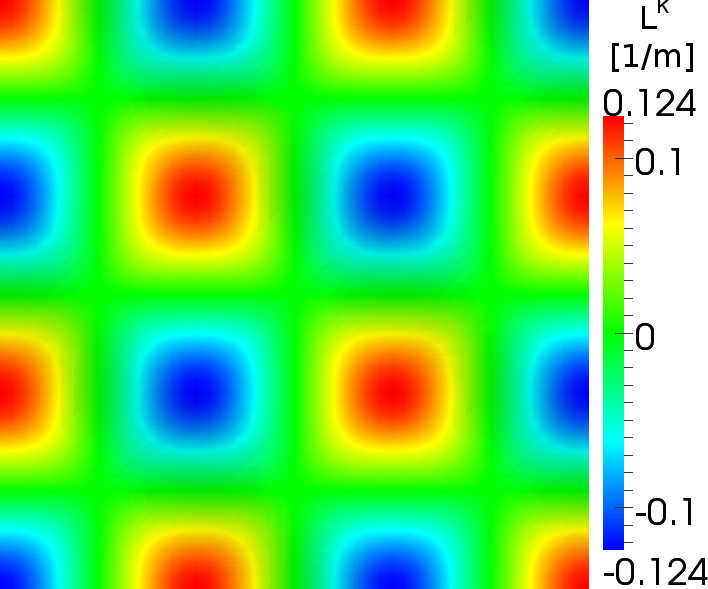}
 \end{minipage}
 \begin{minipage}{.49\textwidth}
  \centering
  \includegraphics[width=0.9\textwidth,height=.5\textwidth,angle=0]{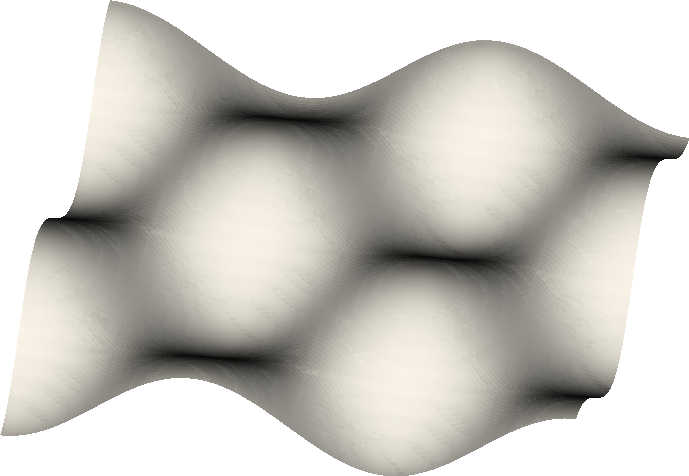}
 \end{minipage}%
 \begin{minipage}{.49\textwidth}
  \centering
  \includegraphics[width=0.8\textwidth,height=.65\textwidth,angle=0]{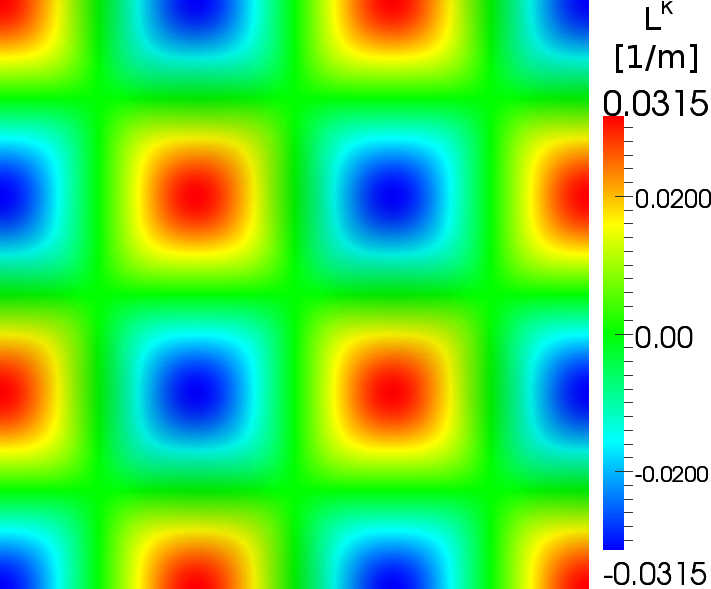}
 \end{minipage}
 \begin{minipage}{.49\textwidth}
  \centering
  \includegraphics[width=0.9\textwidth,height=.5\textwidth,angle=0]{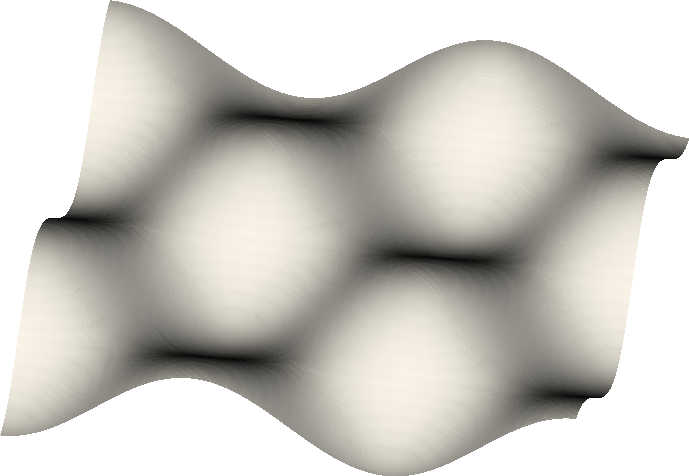}
 \end{minipage}%
 \begin{minipage}{.49\textwidth}
  \centering
  \includegraphics[width=0.8\textwidth,height=.65\textwidth,angle=0]{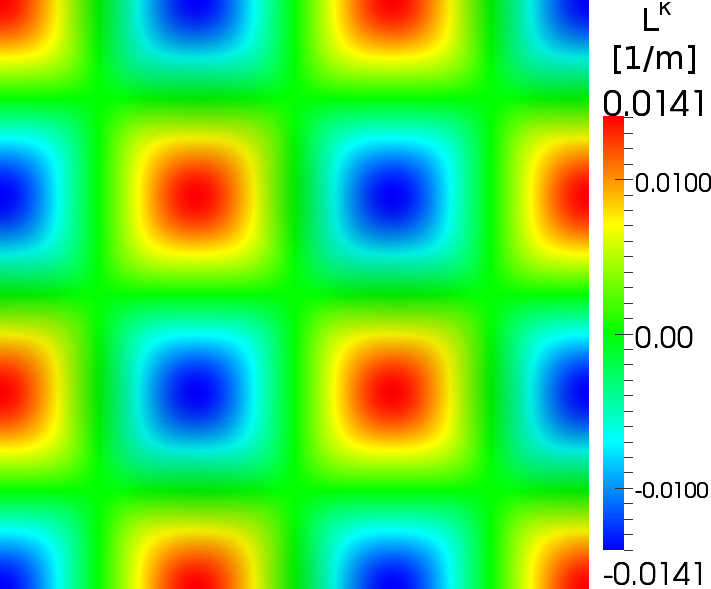}
 \end{minipage}
 \begin{minipage}{.49\textwidth}
  \centering
  \includegraphics[width=0.9\textwidth,height=.5\textwidth,angle=0]{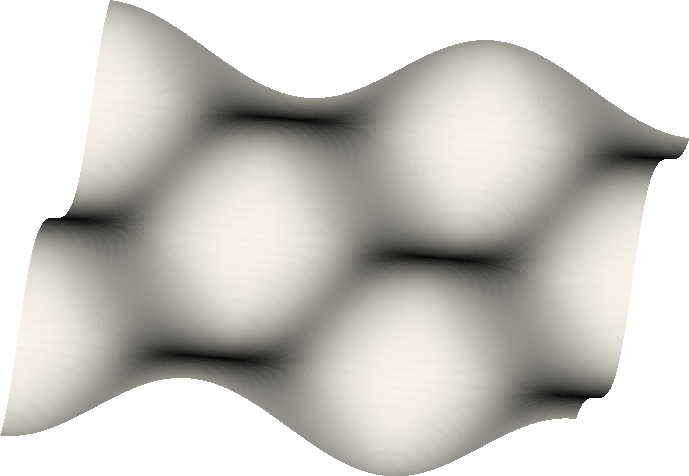}
 \end{minipage}%
 \begin{minipage}{.49\textwidth}
  \centering
  \includegraphics[width=0.8\textwidth,height=.65\textwidth,angle=0]{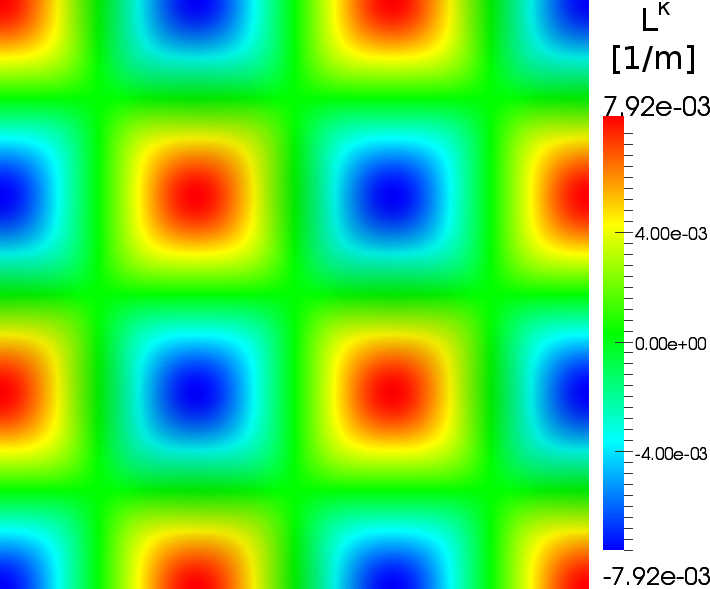}
 \end{minipage}
 \end{centering}
 \caption{\small{The comparison of
                 $\Gamma \lr \psi_0 \rr$ (left) 
                 and
                 $L_i^\kappa=\kappa_i-\kappa_i'$ 
                 norms (right)
                 in the test case $T_1$
                 on four grids $m_1$, $m_2$,
                 $m_{3.5}$ and $m_4$ from top 
                 to bottom.
                 Re-initialization
                 of the 3D wavy interface
                 is performed with 
                 the mapping function $\psio$.
                 The error $L_i^\kappa$ is
                 interpolated to the interface
                 $\Gamma \lr \psi_0 \rr$,
                 the interface width
                 is set to $\eph=\Delta x/2$.}}
 \label{fig25}
 \end{figure}
 \begin{figure}[ht!] \nonumber
 \begin{centering}
 \begin{minipage}{.49\textwidth}
  \centering
  \includegraphics[width=0.9\textwidth,height=.5\textwidth,angle=0]{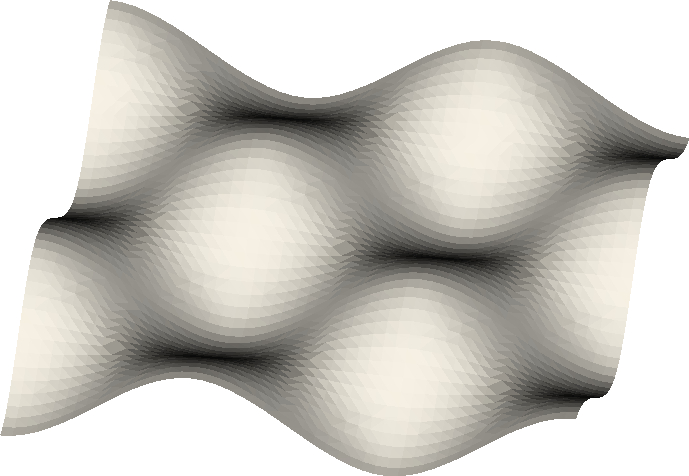}
 \end{minipage}%
 \begin{minipage}{.49\textwidth}
  \centering
  \includegraphics[width=0.8\textwidth,height=.65\textwidth,angle=0]{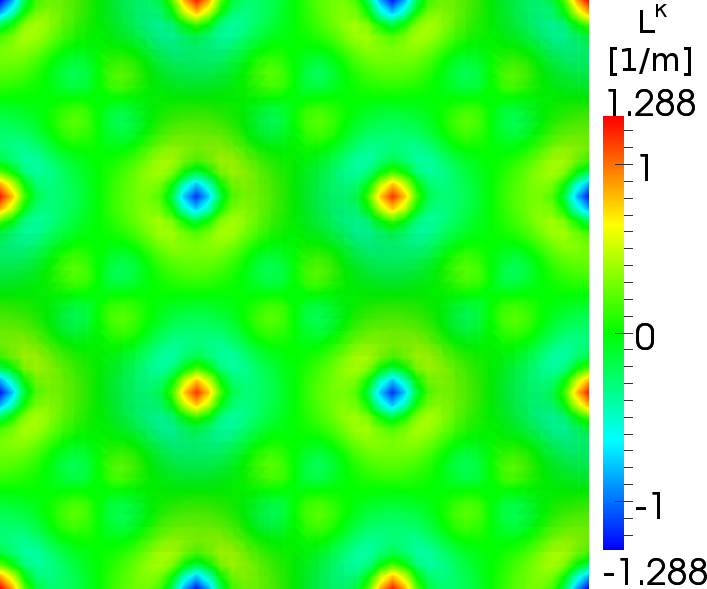}
 \end{minipage}
 \begin{minipage}{.49\textwidth}
  \centering
  \includegraphics[width=0.9\textwidth,height=.5\textwidth,angle=0]{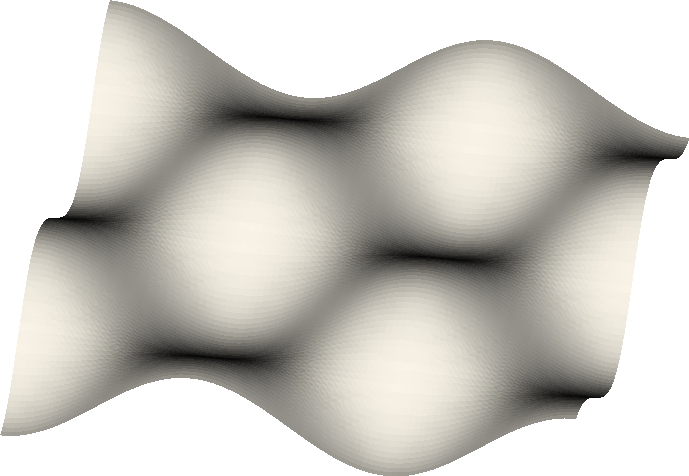}
 \end{minipage}%
 \begin{minipage}{.49\textwidth}
  \centering
  \includegraphics[width=0.8\textwidth,height=.65\textwidth,angle=0]{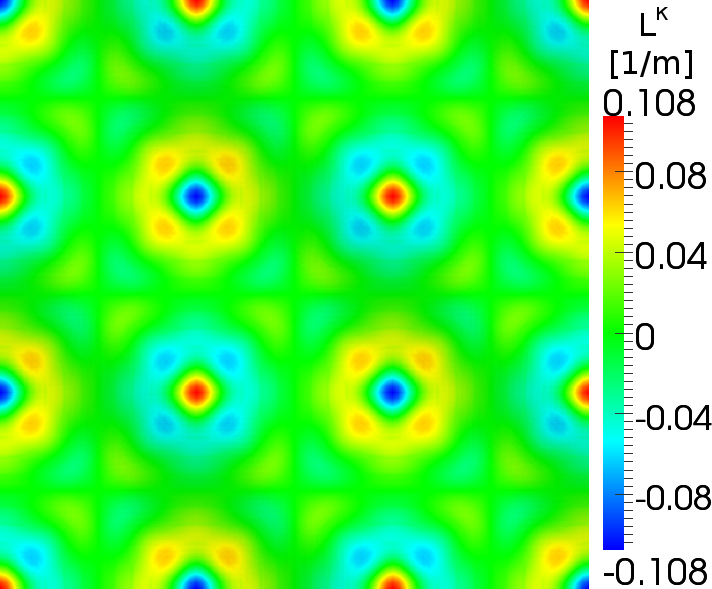}
 \end{minipage}
 \begin{minipage}{.49\textwidth}
  \centering
  \includegraphics[width=0.9\textwidth,height=.5\textwidth,angle=0]{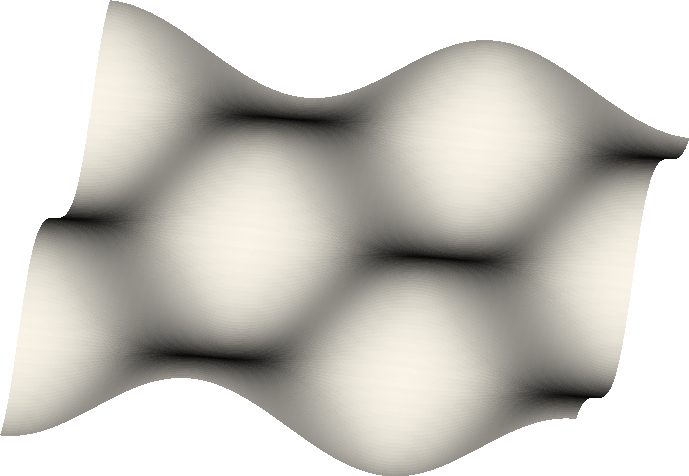}
 \end{minipage}%
 \begin{minipage}{.49\textwidth}
  \centering
  \includegraphics[width=0.8\textwidth,height=.65\textwidth,angle=0]{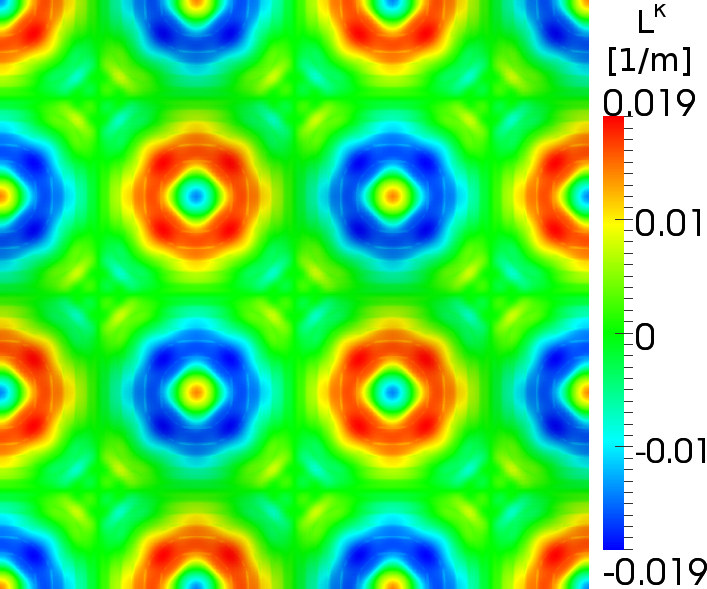}
 \end{minipage}
 \begin{minipage}{.49\textwidth}
  \centering
  \includegraphics[width=0.9\textwidth,height=.5\textwidth,angle=0]{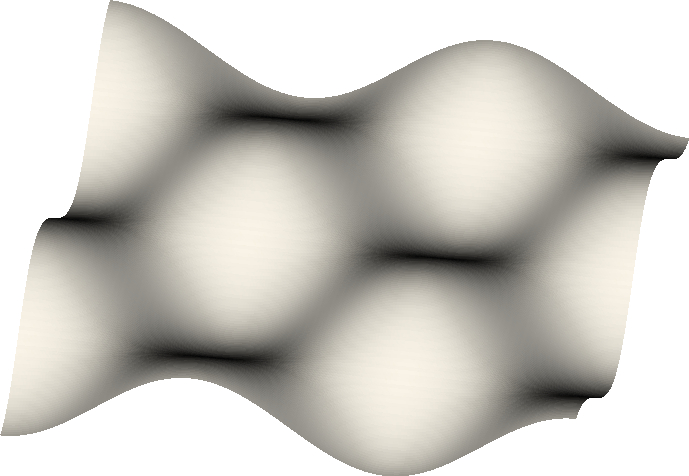}
 \end{minipage}%
 \begin{minipage}{.49\textwidth}
  \centering
  \includegraphics[width=0.8\textwidth,height=.65\textwidth,angle=0]{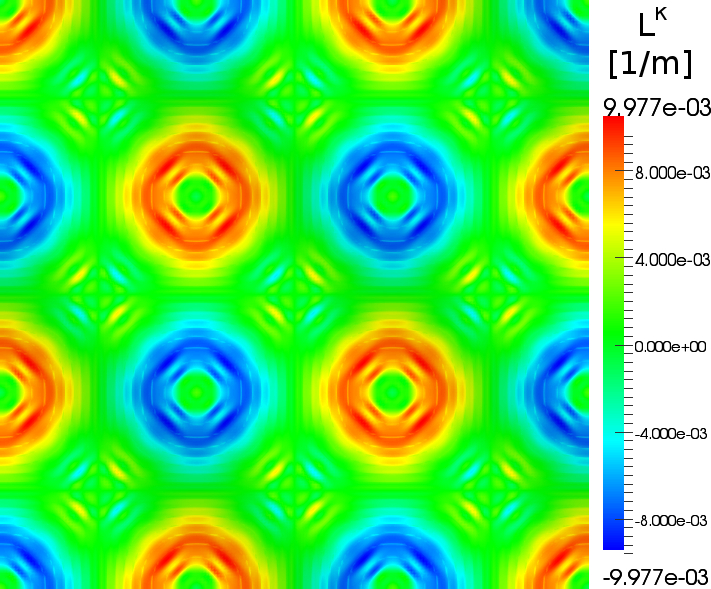}
 \end{minipage}
 \end{centering}
 \caption{\small{The comparison of
                 $\Gamma \lr \psi_0 \rr$ (left)
                 and
                 $L_i^\kappa=\kappa_i-\kappa_i'$ 
                 norms (right)
                 in the test case $T_2$
                 on four grids $m_1$, $m_2$,
                 $m_{3.5}$ and $m_4$ from top 
                 to bottom.
                 Re-initialization
                 of the 3D wavy interface
                 is performed with 
                 the mapping function $\psio$.
                 The error $L_i^\kappa$ is
                 interpolated to the interface
                 $\Gamma \lr \psi_0 \rr$,
                 the interface width
                 is set to $\eph=\Delta x/2$}}
 \label{fig26}
 \end{figure}
 \begin{figure}[ht!] \nonumber
 \begin{centering}
 \begin{minipage}{.49\textwidth}
  \centering
  \includegraphics[width=0.9\textwidth,height=.5\textwidth,angle=0]{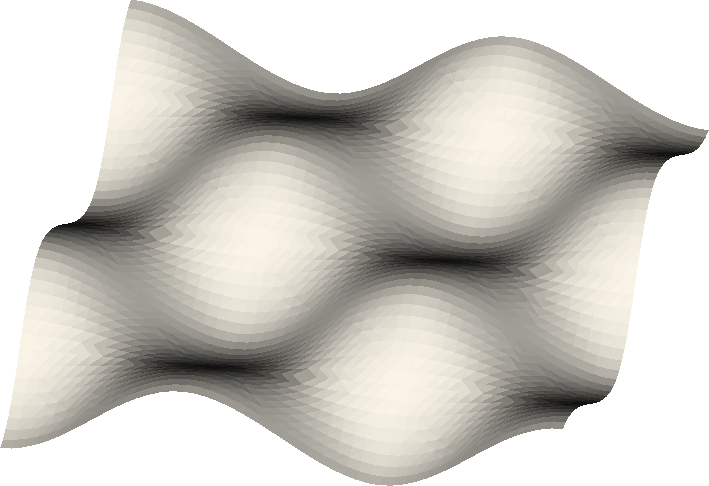}
 \end{minipage}%
 \begin{minipage}{.49\textwidth}
  \centering
  \includegraphics[width=0.8\textwidth,height=.65\textwidth,angle=0]{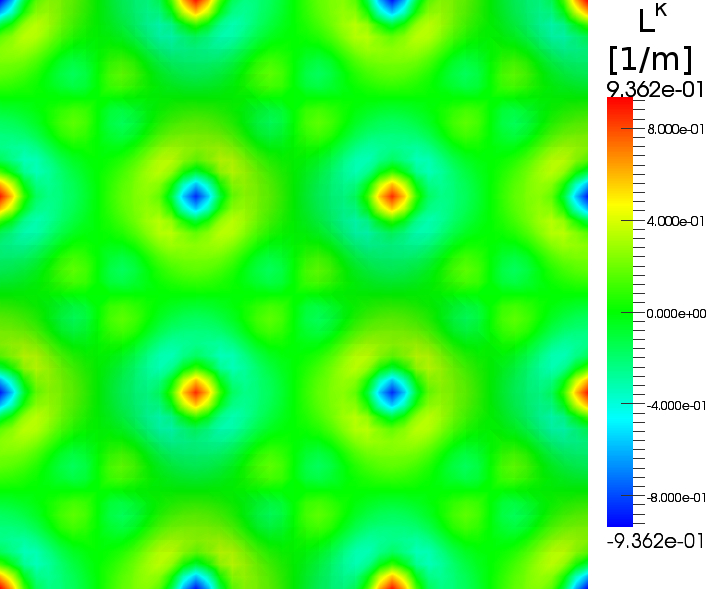}
 \end{minipage}
 \begin{minipage}{.49\textwidth}
  \centering
  \includegraphics[width=0.9\textwidth,height=.5\textwidth,angle=0]{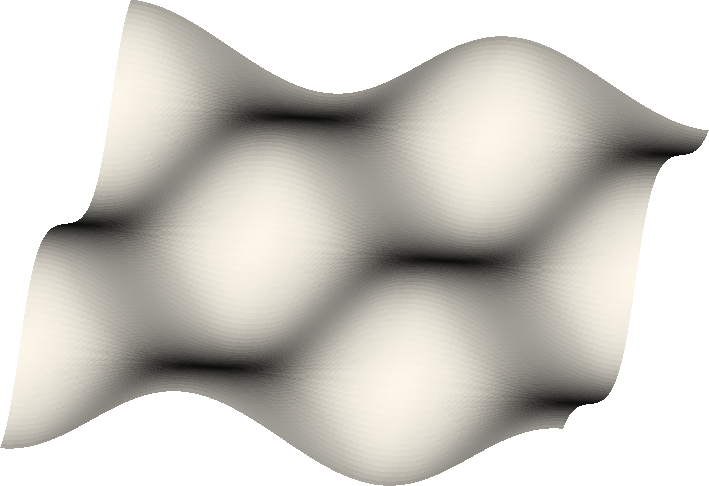}
 \end{minipage}%
 \begin{minipage}{.49\textwidth}
  \centering
  \includegraphics[width=0.8\textwidth,height=.65\textwidth,angle=0]{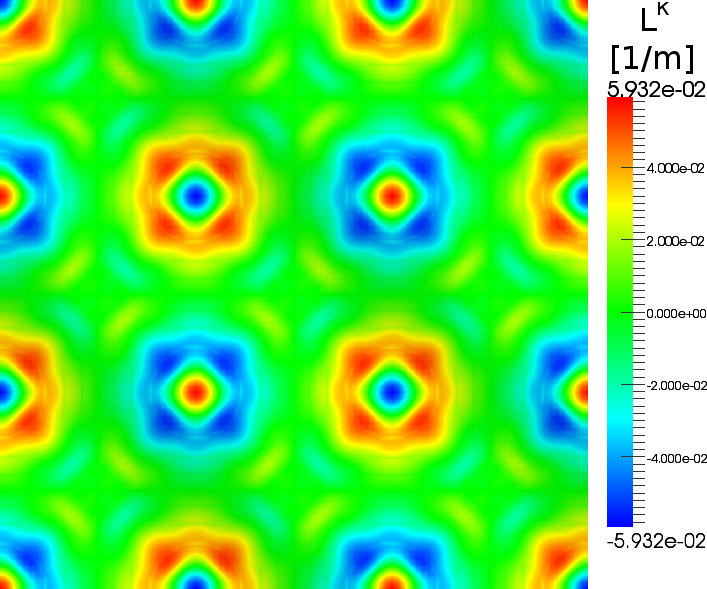}
 \end{minipage}
 \begin{minipage}{.49\textwidth}
  \centering
  \includegraphics[width=0.9\textwidth,height=.5\textwidth,angle=0]{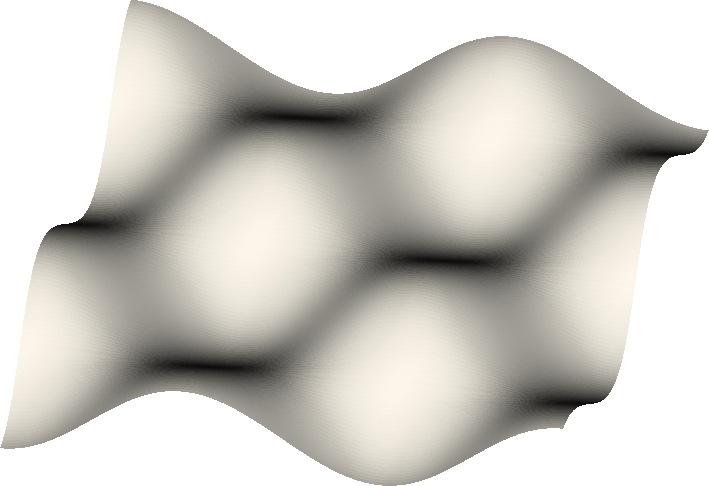}
 \end{minipage}%
 \begin{minipage}{.49\textwidth}
  \centering
  \includegraphics[width=0.8\textwidth,height=.65\textwidth,angle=0]{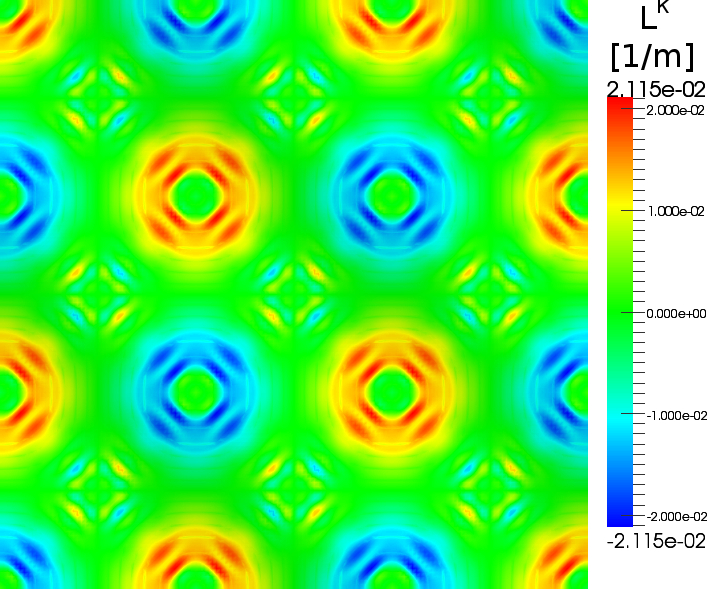}
 \end{minipage}
 \begin{minipage}{.49\textwidth}
  \centering
  \includegraphics[width=0.9\textwidth,height=.5\textwidth,angle=0]{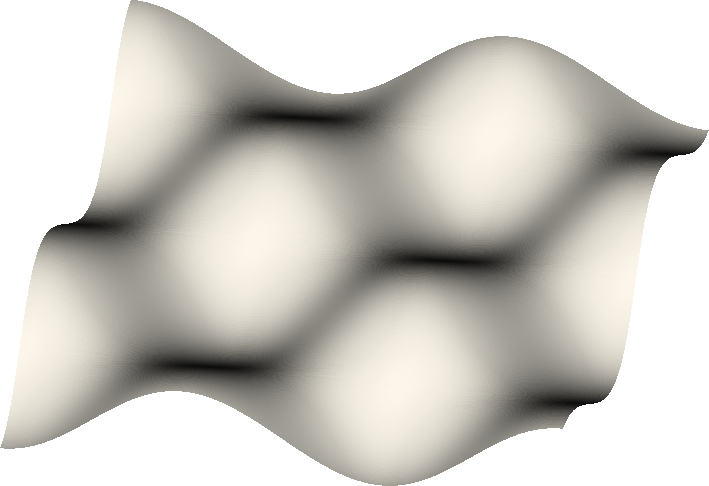}
 \end{minipage}%
 \begin{minipage}{.49\textwidth}
  \centering
  \includegraphics[width=0.8\textwidth,height=.65\textwidth,angle=0]{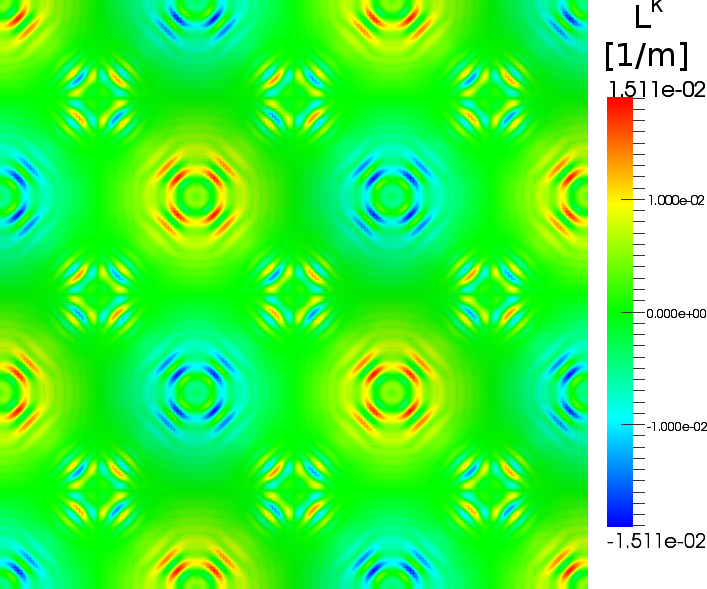}
 \end{minipage}
 \end{centering}
 \caption{\small{The comparison of
                 $\Gamma \lr \psi_0 \rr$ (left)
                 and
                 $L_i^\kappa=\kappa_i-\kappa_i'$ 
                 norms (right)
                 in the test case $T_3$
                 on four grids $m_1$, $m_2$,
                 $m_{3.5}$ and $m_4$ from top 
                 to bottom.
                 Re-initialization
                 of the 3D wavy interface
                 is performed with 
                 the mapping function $\psio$.
                 The error $L_i^\kappa$ is
                 interpolated to the interface
                 $\Gamma \lr \psi_0 \rr$,  
                 the interface width
                 is set to $\eph=\sqrt{3} \Delta x/4$.}}
 \label{fig27}
 \end{figure}

As discussed previously,
there is yet another way in which 
the errors in $\kappa'$ may be
estimated.
In tables \ref{tab9}-\ref{tab10},
values of the
$L_\infty^\kappa=max(|L^\kappa|=|\kappa-\kappa'|)$ 
norms obtained 
with the mapping 
functions $\psigs$, $\psio$
in the test cases $T_2$, $T_3$ 
are given.
Here, almost no difference 
is  detected between the convergence rates  
obtained using $\Gamma \lr \alpps \rr$ or $\Gamma \lr \psio \rr$ 
interface representations,
also depicted in figure \ref{fig24}(c)(f). 
The values in tables \ref{tab9}-\ref{tab10}
show the second-order
convergence rate of $\kappa_i'$ 
is detected if we take into account
only the errors  at the interface 
$\Gamma \lr \alpha \rr$ or $\Gamma \lr \psi_0 \rr$
(compare with $L_\infty^\kappa$ norms
presented in 
figures \ref{fig25}-\ref{fig27}).
In the case $T_3$, 
the convergence rate of $\kappa_i'$
on the finest grid $m_4$ 
is somewhat lower;
this may be explained by
the reduced number of grid points due to 
the smaller interface width 
$\eph=\sqrt{3}\Delta x/4 < \Delta x/2$ used
to resolve $\alpps$ and $\psio$.

Finally, 
in figures \ref{fig25}-\ref{fig27}
we can compare 
shapes of 
the reconstructed interfaces 
obtained in the test cases $T_i$ where $i=1,2,3$.
The differences 
in the interface $\Gamma$ shapes 
presented in these figures 
are barely recognizable.
\begin{table}[ht]
\centering 
\begin{tabular}{@{}  l  c c  c c @{}} 
\cline{1-5}
 m.f.    & \multicolumn{2}{c}{$\psigs$} & \multicolumn{2}{c}{$\psio$} \\ [0.5ex] 
\cline{1-5} 
 $\Gamma$    & $\alpha$ & $\psi_0'(\psi_1)$ & $\alpha$ & $\psi_0$ \\ [0.5ex] 
\cline{1-5} 
$m_2$    &  1.28881    &  1.28783    &  1.28763    &  1.28772    \\ 
$m_3$    &  1.0813e-1  &  1.0771e-1  &  1.0816e-1  &  1.0767e-1  \\
$m_{3.5}$ &  1.9419e-2  &  1.9067e-2  &  1.9413e-2  &  1.9061e-2  \\ 
$m_4$    &  1.0741e-2  &  9.9776e-3  &  1.0721e-2  &  9.9771e-3  \\[0.5ex]  
\cline{1-5}  
\end{tabular}
\caption{Convergence 
         of the interface curvature $\kappa_i'$
         measured by $L_\infty^\kappa$ norm
         obtained from $|L_i^\kappa|=|\kappa_i-\kappa_i'|$          
         at the interface $\Gamma \lr \alpha \rr$
         and $\Gamma \lr \psi_0 \rr$
         in the test case $T_2$, see   figure \ref{fig24}(top).}
\label{tab9} 
\end{table}
\begin{table}[ht]
\centering 
\begin{tabular}{@{}  l  c c @{}} 
\cline{1-3}
 m.f. & \multicolumn{2}{c}{$\psio$}  \\ [0.5ex] 
\cline{1-3} 
$\Gamma$ &  $\alpha$ & $\psi_0$ \\ [0.5ex] 
\cline{1-3} 
$m_2$    & 9.3597e-1    &  9.3625e-1  \\ 
$m_3$    & 5.9071e-1    &  5.9321e-2  \\
$m_{3.5}$ & 2.0281e-2    &  2.1154e-2  \\ 
$m_4$    & 1.4133e-2    &  1.5113e-2  \\[0.5ex]  
\cline{1-3}  
\end{tabular}
\caption{Convergence 
         of the curvature 
         $\kappa_i'$ obtained
         with the mapping
         function (m.f.) $\psio$ 
         measured by 
         $L_\infty^\kappa$ norm
         obtained from $|L_i^\kappa|=|\kappa_i-\kappa_i'|$          
         at the interface $\Gamma \lr \alpha \rr$
         and $\Gamma \lr \psi_0 \rr$
         in the test case $T_3$, see  figure \ref{fig24}(bottom).}
\label{tab10} 
\end{table}
The shape of the wavy interface 
is preserved on all meshes 
used in the present study
almost independently
from the number of 
re-initialization
steps $N_\tau$.
Although the reduction of the interface
width $\eph$ leads to slower convergence
of the solution to the re-initialization
equation, and slower convergence of its
curvature $\kappa_i'$, this does
not affect
the wavy interface shape.

In summary, 
the results presented 
in tables \ref{tab9}-\ref{tab10}
and in figure \ref{fig24} confirm 
the second-order convergence rate 
of the interface  
curvature $\kappa$
may be achieved within
the second-order accurate
finite volume solver
using the conservative level-set
method and the consistent
re-initialization procedure.

\subsection{Tests with advection}
\label{sec6.2}

The primary aim of 
the studies performed
below is verification
of the new re-initialization
method during
advection of $\alpps$ and $\psio$.
Next,
we solve 
\Eq{eq30a}
with $w_i=u_i$, 
where $u_i$
is the given
divergence free
velocity field.
As we improve the re-initialization method 
first proposed in \cite{olsson05},
similar advection test cases 
are carried out 
for comparison.
In particular, 
we want to investigate
the area (mass) conservation
of the new method
in the case of the variable 
number 
of re-initialization steps 
$N_\tau$. 
Using experience 
gained during the tests
in Section
\ref{sec6.1},
only the mapping function
$\psio$ is used for discretization
of $|\nabla \alpha|$ in \Eq{eq1}
and hence \Eq{eq29} is solved
during the re-initialization step.
The interface width is set to 
$\eph=\sqrt{2}\Delta x/4$
and $\Delta \tau = D/C^2=\eph$.
Present investigations 
are performed in quadratic
domain 
$\Omega = <0,1>\times<0,1>\,m^2$  
on gradually refined grids
$m_i=2^{4+i}\times 2^{4+i}$ 
with the uniform grid nodes 
distribution.

The advection equation 
(\ref{eq30a}) 
is discretized 
in time using
the first-order 
implicit Euler method,
and in space 
using the 
deferred-correction 
method with 
the second-order TVD MUSCL flux 
limiter from \cite{xue98}.
This type 
of spatial 
and temporal
discretization 
is a default 
technique
used in the Fastest solver
for discretization of
advection terms
in all transport 
equations.
As the main goal of
this paper is the
improvement of 
the re-initialization method,
the detailed 
solutions
to numerical issues 
that arise during 
coupling
of equations
(\ref{eq29}) and  
(\ref{eq30a}) 
are left for future 
investigations.

\subsubsection{Rotating circle}
\label{sec6.2.1}
%
In this section, 
a circular interface
$\Gamma$ revolving in the divergence free 
velocity field
$(u_1,u_2)=V_0/L(y-0.5,0.5-x)$
where $V_0=1\,m/s$ and $L=1\,m$
is studied.
Initially at $t=0$,
the center of the circle
with the radius $R=0.15\,m$ 
is located 
in the point $(x_0,y_0)=(0.65\,m,0.5\,m)$.
The time step size during integration
of equation (\ref{eq30a}) is set 
to $\Delta t = 2\pi/N_t$,
where $N_t=360 \cdot i$ 
and $i=1,\ldots,4$ denotes 
the grid number.

\begin{figure}[ht!] \nonumber
 \begin{centering}
 \begin{minipage}{.25\textwidth}
  \centering
  \includegraphics[width=1.\textwidth,height=1.\textwidth,angle=-90]{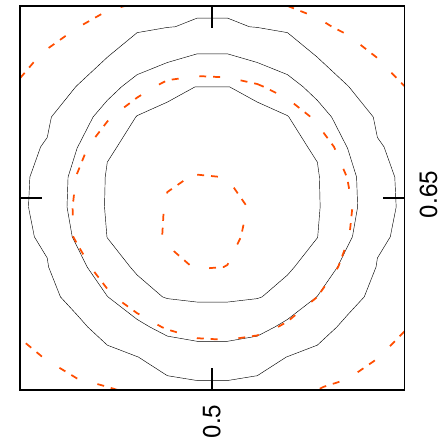}
 \end{minipage}%
 \begin{minipage}{.25\textwidth}
  \centering
  \includegraphics[width=1.\textwidth,height=1.\textwidth,angle=-90]{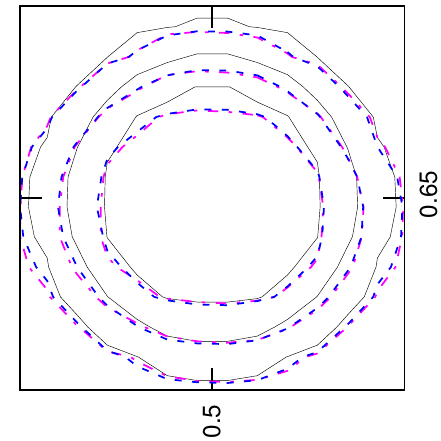}
 \end{minipage}%
 \begin{minipage}{.25\textwidth}
  \centering
  \includegraphics[width=1.\textwidth,height=1.\textwidth,angle=-90]{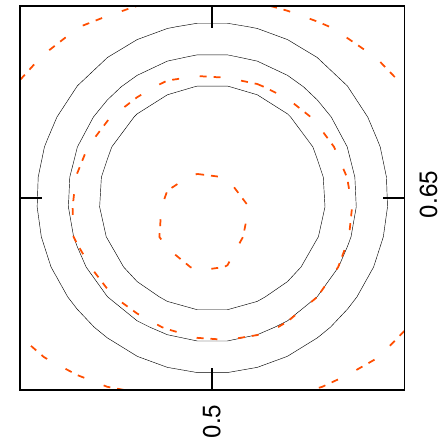}
 \end{minipage}%
 \begin{minipage}{.25\textwidth}
  \centering
  \includegraphics[width=1.\textwidth,height=1.\textwidth,angle=-90]{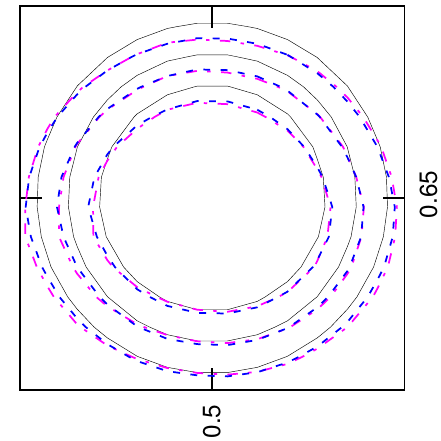}
 \end{minipage}
 \begin{minipage}{.25\textwidth}
  \centering
  \includegraphics[width=1.\textwidth,height=1.\textwidth,angle=-90]{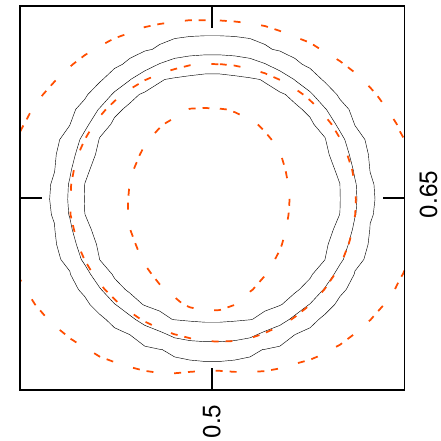}
 \end{minipage}%
 \begin{minipage}{.25\textwidth}
  \centering
  \includegraphics[width=1.\textwidth,height=1.\textwidth,angle=-90]{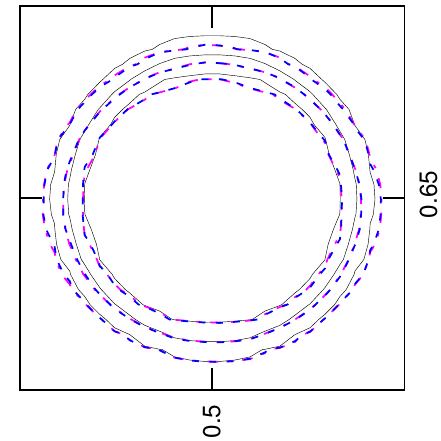}
 \end{minipage}%
 \begin{minipage}{.25\textwidth}
  \centering
  \includegraphics[width=1.\textwidth,height=1.\textwidth,angle=-90]{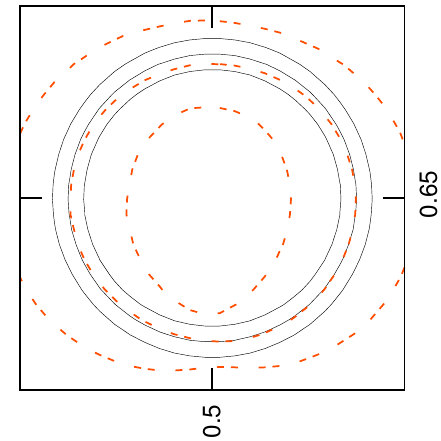}
 \end{minipage}%
 \begin{minipage}{.25\textwidth}
  \centering
  \includegraphics[width=1.\textwidth,height=1.\textwidth,angle=-90]{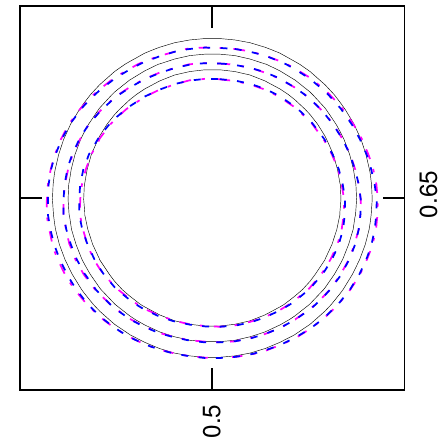}
 \end{minipage}
 \begin{minipage}{.25\textwidth}
  \centering
  \includegraphics[width=1.\textwidth,height=1.\textwidth,angle=-90]{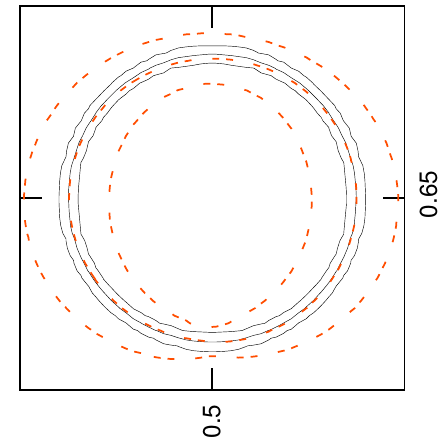}
 \end{minipage}%
 \begin{minipage}{.25\textwidth}
  \centering
  \includegraphics[width=1.\textwidth,height=1.\textwidth,angle=-90]{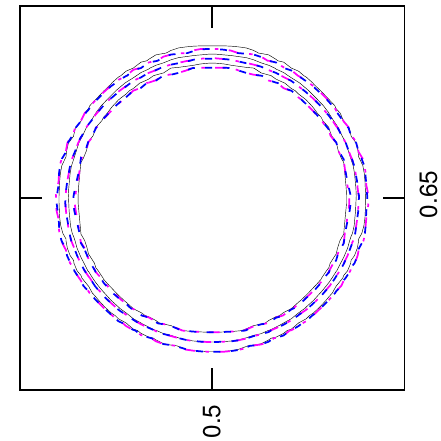}
 \end{minipage}%
 \begin{minipage}{.25\textwidth}
  \centering
  \includegraphics[width=1.\textwidth,height=1.\textwidth,angle=-90]{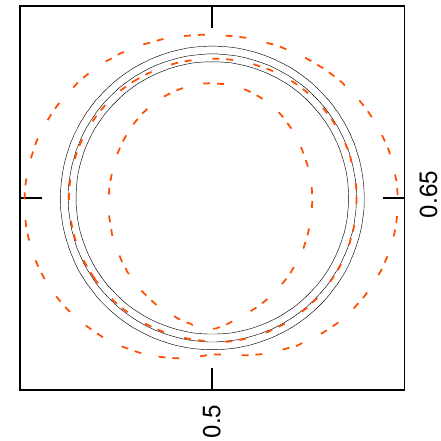}
 \end{minipage}%
 \begin{minipage}{.25\textwidth}
  \centering
  \includegraphics[width=1.\textwidth,height=1.\textwidth,angle=-90]{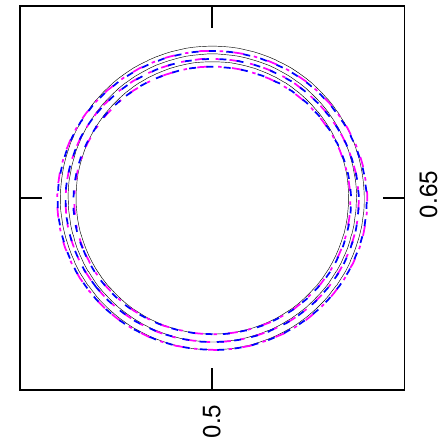}
 \end{minipage}
 \begin{minipage}{.25\textwidth}
  \centering
  \includegraphics[width=1.\textwidth,height=1.\textwidth,angle=-90]{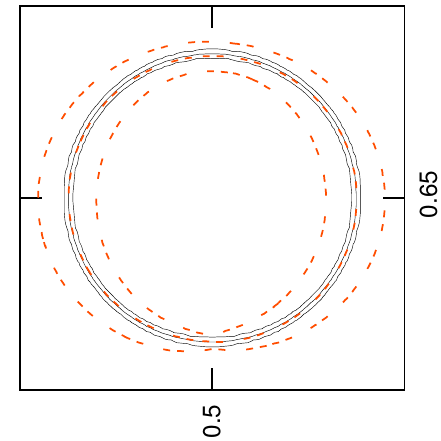}
 \end{minipage}%
 \begin{minipage}{.25\textwidth}
  \centering
  \includegraphics[width=1.\textwidth,height=1.\textwidth,angle=-90]{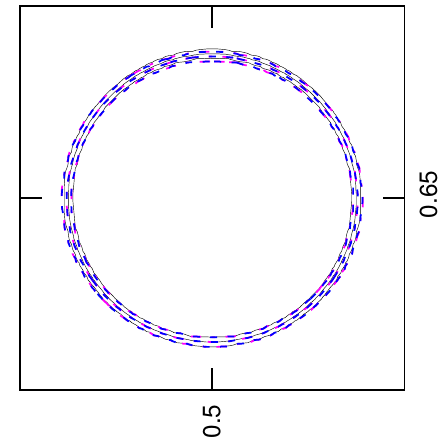}
 \end{minipage}%
 \begin{minipage}{.25\textwidth}
  \centering
  \includegraphics[width=1.\textwidth,height=1.\textwidth,angle=-90]{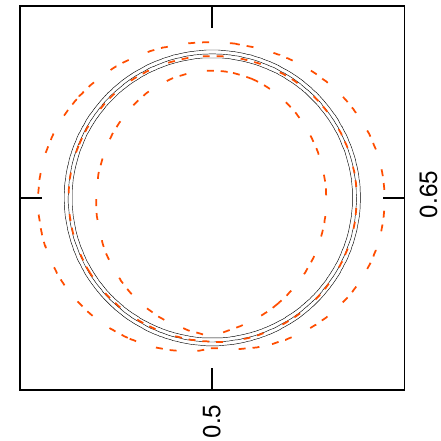}
 \end{minipage}%
 \begin{minipage}{.25\textwidth}
  \centering
  \includegraphics[width=1.\textwidth,height=1.\textwidth,angle=-90]{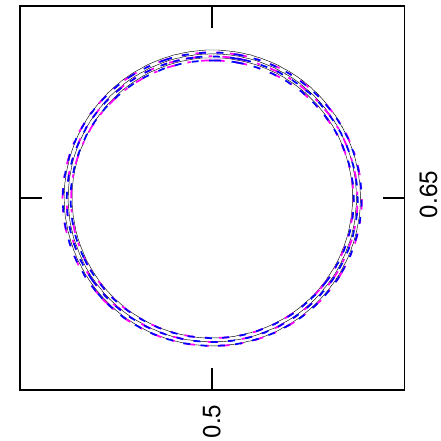}
 \end{minipage}
 \end{centering}
 \caption{\small{The comparison of
                 $\alpps = (0.05,0.5,0.95)$ 
                 iso-lines (two columns left)
                 and $\psio = (r_i^1,0,r_i^2)$  
                 iso-lines (see \Eqs{eq32}{eq33}, two columns right)
                 obtained after one revolution 
                 of the circular interface
                 on grids $m_i$, $i=1,\ldots,4$
                 (from top to bottom) with the
                 initial condition (black solid line).
                 The results were obtained with 
                 $N_\tau=0$ (orange double-dashed line),
                 $N_\tau=1$ (magenta dashed-dotted line) and
                 $N_\tau=32$ (dark blue dashed line)
                 re-initialization steps.}}             
 \label{fig28}
 \end{figure}

In figure \ref{fig28},
convergence of the two interface representations
$\alpps$ and $\psio$ towards the initial condition
is presented.
As previously
observed in \cite{olsson05},
it is clear that an introduction of 
the re-initialization step
after the advection step
improves 
the conservation
of the shape and area of
the advected interface.
In the previous sections
we have shown that 
the new re-initialization
method does not change the
position of the 
stationary interfaces.
For this reason 
small deformations of $\Gamma$
after one revolution may
be attributed to numerical
errors introduced 
by the advection
scheme, see  \Fig{fig28}.
The main difference between
our results and the results presented
in \cite{olsson05} is independence 
of the interface shape 
from the number of re-initialization
steps $N_\tau  \le 32$.

To measure the convergence rate of the
numerical approximation of the interface $\Gamma$
towards analytical solution, the first
order norm is introduced
\blnm
\begin{align}
 \begin{split}
   L_{1,i}^r = \frac{1}{N_i^s} \sum_l^{N_i^s} | \br_{l,i}^{num} -\br_{l,i}^{ext} |,
  \label{eq43b}
 \end{split}
\end{align}
\elnm
where $N_i^s$ is the number
of points in the interface $\Gamma$ 
representation given 
by the iso-lines of 
$\alpps=1/2$ or  $\psio=0$,
$\br_{l,i}^{num} = |\bx_{l,i}^{num}-\bx_0|$,
$\br_{l,i}^{ext} = |\bx_{l,i}^{ext}-\bx_0|$,
where $\bx_{l,i}^{num}$ and 
$\bx_{l,i}^{ext}$ are points 
obtained from the numerical 
and the exact approximation
of the interface $\Gamma$ (given by the initial
condition) on the grid $m_i$.
As in the present numerical scheme,
we use the first-order
discretization method in time.
This convergence rate 
of $\Gamma$ towards 
the initial condition
is also observed 
in \Tab{tab11} for both interface 
representations.
\begin{table}[ht]
\centering 
\begin{tabular}{@{}  l  c c c  c c c @{}} 
\cline{1-7}
$\Gamma$ & \multicolumn{3}{c}{$\alpps$} & \multicolumn{3}{c}{$\psio$} \\ [0.5ex] 
\cline{1-7} 
 $N_\tau$ &   1      & 16          & 32        & 1           & 16         &  32  \\ [0.5ex] 
\cline{1-7} 
$m_1$    & 8.3833e-3 & 8.1417e-3 & 7.9349e-3  &  8.5084e-3 & 8.3689e-3 & 8.1521e-3  \\ 
$m_2$    & 4.4119e-3 & 4.5067e-3 & 4.4682e-3  &  4.3593e-3 & 4.4516e-3 & 4.4166e-3  \\ 
$m_3$    & 2.1589e-3 & 2.2262e-3 & 2.2192e-3  &  2.1621e-3 & 2.2296e-3 & 2.2221e-3  \\
$m_4$    & 1.0916e-3 & 1.1243e-3 & 1.1239e-3 &   1.0923e-3 & 1.1247e-3 & 1.1244e-3  \\ [0.5ex] 
\cline{1-7}  
\end{tabular}
\caption{The first-order convergence of $L_{1,i}^r$ norm,
         computed after one revolution of 
         the circular interface represented by $\Gamma \lr \alpps \rr$
         and $\Gamma \lr \psio \rr$
         with $N_\tau \le 32$ re-initialization steps.}
\label{tab11} 
\end{table}
The convergence rate
of the reconstructed interface $\Gamma$
towards the initial condition
is related to  
details in the coupling
between equations
(\ref{eq29}) and (\ref{eq30a})
and could be  
improved by  introduction
of the higher-order 
discretization schemes.
We note that 
the values of $L_{1,i}^r$ norm
in table \ref{tab11}, remain 
almost the same in spite of
a different number of re-initialization
steps $N_\tau$ used during simulation.

Since the present method  
uses two representations of the interface
$\alpps$ and $\psio$, the area (mass) conservation 
can be estimated in two different ways.
In order to distinguish 
between these two possibilities, 
a cumulative error is defined as
\blnm
\begin{align}
 \begin{split}
 &  E_t = \frac{1}{N_t} \sum_n^{N_t} E_n \quad \textrm{where} \\
 &  E_n = |S_{num}^n-S_{ext}^n| = \int | \alpha_n \lr \psi_0 \rr - \alpha_n^{ext} \lr \psi_0^{ext} \rr | dS,
  \label{eq44}
 \end{split}
\end{align}
\elnm
and $\alpha^{ext}$, $\psi_0^{ext}$
are defined by 
the equations (\ref{eq2}) and (\ref{eq31}),
respectively.
$N_t$ in equation (\ref{eq44}) 
denotes the total number 
of time steps on the grid $m_i$ 
during one revolution of 
the interface $\Gamma$, 
$E_n$ is a difference between 
the numerical approximation of 
the surface $S_{num}^n$
and the exact surface $S_{ext}^n$
which are calculated on each  
time level $n$. 
The instantaneous 
position of the circular
interface center
$(x_c,y_c)$ required 
to compute of $S_{ext}^n$
is obtained from 
\blnm
\be
  x_c  = \frac{\sum_{i,j} \alpha_{i,j} x_{i,j}}{\sum_{i,j} \alpha_{i,j} }, \quad
  y_c  = \frac{\sum_{i,j} \alpha_{i,j} y_{i,j}}{\sum_{i,j} \alpha_{i,j} }.
  \label{eq45}
\ee
\elnm

We found that convergence
rates of the area (mass)
depend on the region of
integration of $\alpps$
in \Eq{eq44}.
 \begin{figure}[ht!] \nonumber
 \includegraphics[angle=-90]{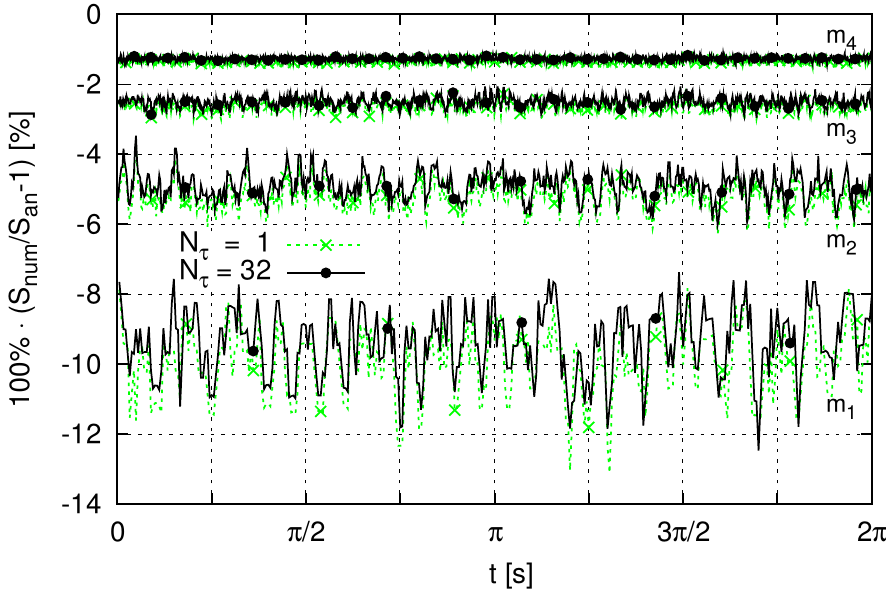}
 \includegraphics[angle=-90]{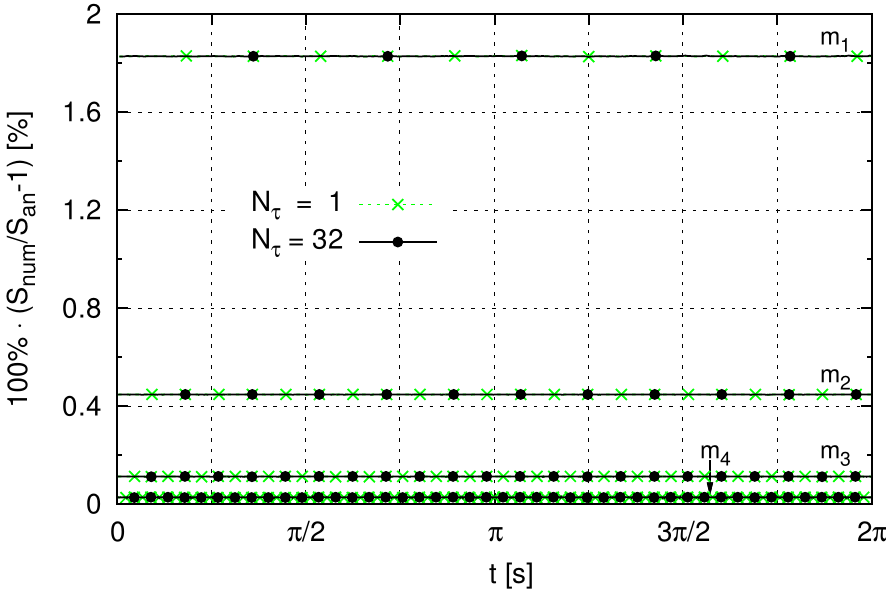}
 \caption{\small{The area conservation  
                 calculated in the regions:
                 (a) $r_1 = \lc x_i \,\big|\, 1-\alpps \ge 0.5  \rc$,
                 (b) $r_2 = \lc x_i \,\big|\, \psio \le 8\eph \rc$ 
                 during one revolution
                 of the circular interface 
                 with $N_\tau = 1$ or $N_\tau = 32$, 
                 the error is 
                 normalized with $S_{an}=\pi R^2$.}}
 \label{fig29}
 \end{figure}
 \begin{figure}[ht!] \nonumber
 \includegraphics[angle=-90]{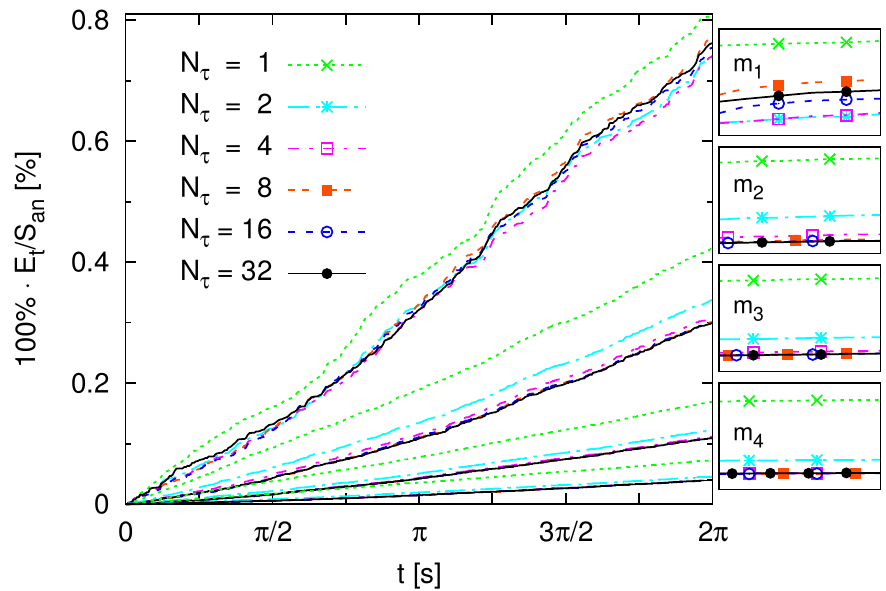}
 \includegraphics[angle=-90]{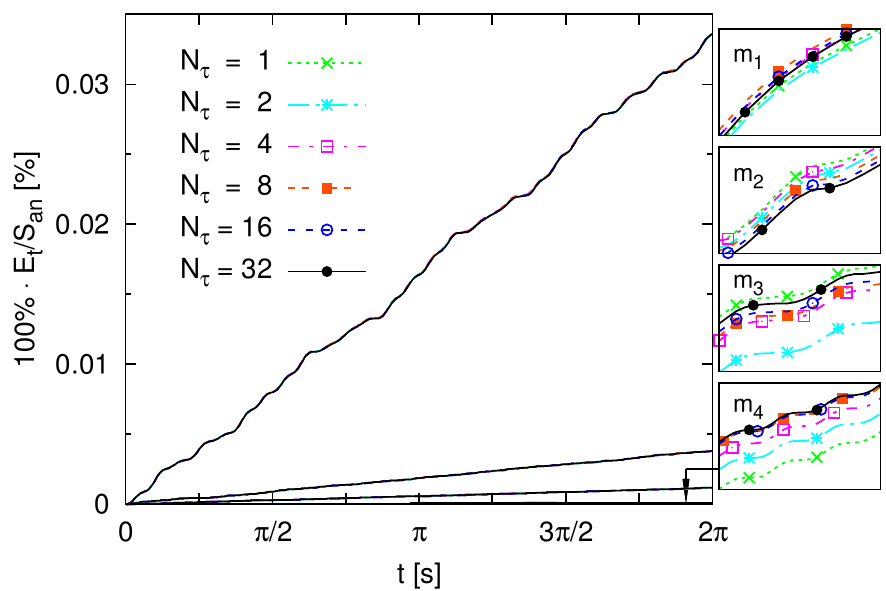}
 \caption{\small{Convergence of the cumulative error $E_t$
                 defined by \Eq{eq44} with the different number of 
                 re-initialization steps $N_\tau \le 32$. 
                 $E_t$ is calculated in the regions: 
                 (a) $r_1 = \lc x_i \,\big|\, 1-\alpps \ge 0.5  \rc$,
                 (b) $r_2 = \lc x_i \,\big|\, \psio \le 8\eph \rc$ 
                 and is normalized with $S_{an}=\pi R^2$.
                 The boxes on the right, 
                 present convergence during last
                 ten time steps $\Delta t$,
                 the number and position of 
                 symbols is arbitrary.}}
 \label{fig30}
 \end{figure}
If the region $r_1=\lc x_i \,\big|\, 1-\alpps \ge 0.5 \rc$ 
or $r_1=\lc x_i \,\big|\, \psio \le 0 \rc$  is chosen, 
then the first-order
convergence rate of area 
is obtained, 
see figures \ref{fig29}(a)-\ref{fig30}(a).
Similar convergence rates of
the Heaviside function were reported
in \cite{olsson05} although therein
the second-order discretization
in time was employed.
When the integration 
is carried out in the region 
$r_2= \lc x_i \,\big|\, \psio \le 8 \eph \rc$
or in the whole computational domain $\Omega$,
then almost the second-order 
area convergence rate
may be deduced 
from the results depicted 
in figure (\ref{fig29})(b). 
In this latter case, convergence
of the interface at each time step
$\Delta t$ with $N_\tau \le 32$ 
(see figure \ref{fig30}(b))
is less evident but
may be also observed 
for example on the grid
$m_4$.
The differences in the magnitude of $E_t$ errors 
visible in figures \ref{fig30}(a)(b) 
may be explained by the more accurate interface
representation with the signed distance function, 
and thus smaller error $E_n$ during 
the area (mass)
computation in 
the region $r_2$, see \Eq{eq44}.

Here and in the following examples,
it becomes clear that $S_{num}$ calculated
in the region  $r_1=\lc x_i \,\big|\, 1-\alpps \ge 0.5 \rc$
indicates how the area (mass) varies during
topological changes of the
interface (stretching, break up, coalescence), 
whereas $S_{num}$ calculated
in the region 
$r_2= \lc x_i \,\big|\, \psio \le 8 \eph \rc$
allows examination of
the total area (mass)
conservation.
Since 
the circular interface is
rotated as a rigid body
in the present case,
the errors integrated in
both regions $r_1$ and $r_2$ 
remain on 
constant levels, 
see \Fig{fig29}.

During rotation of the circular interface
without deformations, 
obtained area convergence rates
are almost independent from
the number of selected re-initialization
steps $N_\tau$  (compare results with
different $N_\tau$ in \Fig{fig29}
and in \Tab{tab11}).
The results presented in figure \ref{fig30}
suggest that in most cases up to four re-initialization
steps should be sufficient  to preserve 
the shape and  area of the interface $\Gamma$.
%

\subsubsection{Vortex test}
\label{sec6.2.2}

To test the new re-initialization
method in a more complex velocity field,
we use 
\blnm
\begin{align}
 \begin{split}
    u_1 =-V_0 sin^2 \lr k x \rr sin \lr 2k y \rr, \quad
    u_2 = V_0 sin^2 \lr k y \rr sin \lr 2k x \rr,
  \label{eq46}
 \end{split}
\end{align}
\elnm
where $k=\pi/L$, $L=1\,m$ and $V_0=1\,m/s$ (similar
to the test case in \cite{olsson05}).
A circle with the radius $R=0.15\,m$ is 
initially located at $(x_0,y_0)=(0.5\,m,0.75\,m)$.
The simulation time is $t=2\,s$, 
as after time $T=1\,s$ the flow
field is reversed so  the
exact solution should be obtained
after the same number of time steps.
In order to satisfy the Courant 
condition $Cr \approx 0.2$, 
the time step size is set to $\Delta t = \Delta x/8$ 
to give the total number of
time steps per revolution 
$N_t=8 T/\Delta x = 8 T N_c$.

\begin{figure}[ht!] \nonumber
 \begin{centering}
 \begin{minipage}{0.25\textwidth}
  \centering
   \includegraphics[width=1.\textwidth,height=1.\textwidth,angle=-90]{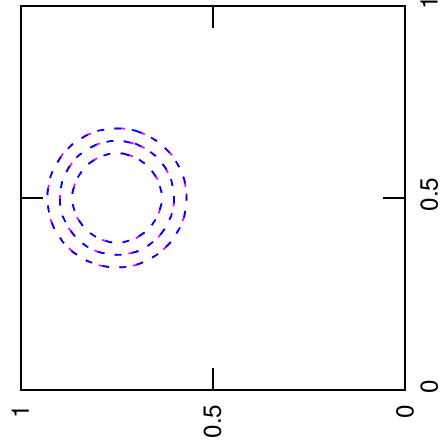}
 \end{minipage}%
 \begin{minipage}{0.25\textwidth}
  \centering
   \includegraphics[width=1.\textwidth,height=1.\textwidth,angle=-90]{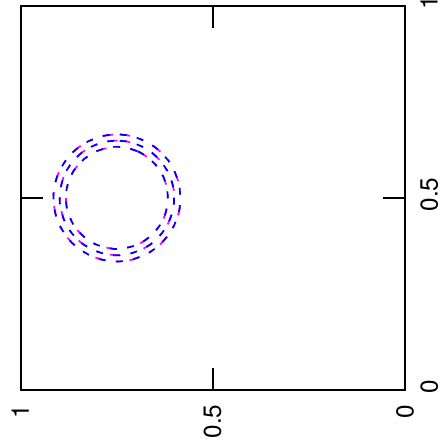}
 \end{minipage}%
 \begin{minipage}{0.25\textwidth}
  \centering
   \includegraphics[width=1.\textwidth,height=1.\textwidth,angle=-90]{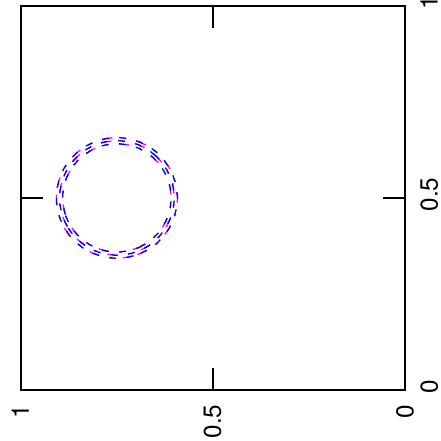}
 \end{minipage}%
 \begin{minipage}{0.25\textwidth}
  \centering
   \includegraphics[width=1.\textwidth,height=1.\textwidth,angle=-90]{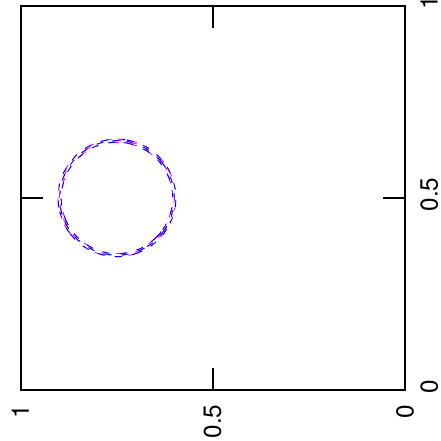}
 \end{minipage}
 \begin{minipage}{0.25\textwidth}
  \centering
   \includegraphics[width=1.\textwidth,height=1.\textwidth,angle=-90]{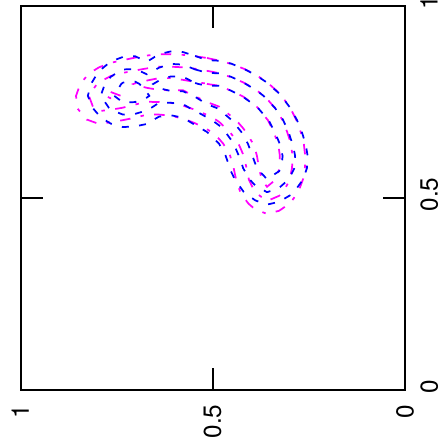}
 \end{minipage}%
 \begin{minipage}{0.25\textwidth}
  \centering
   \includegraphics[width=1.\textwidth,height=1.\textwidth,angle=-90]{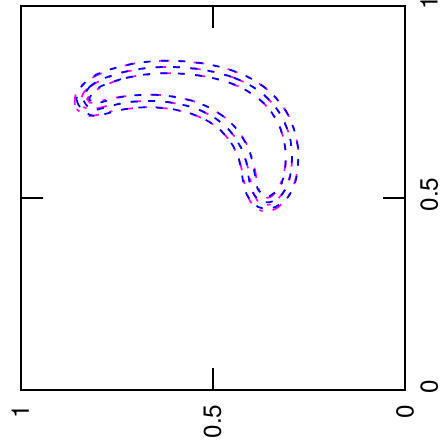}
 \end{minipage}%
 \begin{minipage}{0.25\textwidth}
  \centering
   \includegraphics[width=1.\textwidth,height=1.\textwidth,angle=-90]{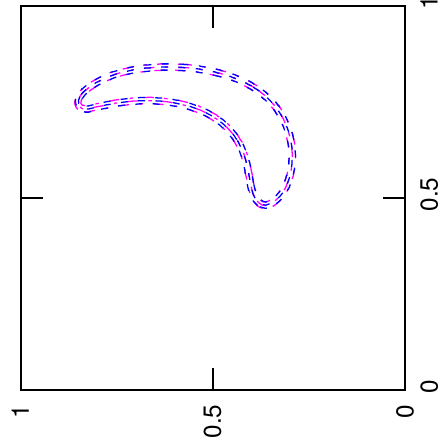}
 \end{minipage}%
 \begin{minipage}{0.25\textwidth}
  \centering
   \includegraphics[width=1.\textwidth,height=1.\textwidth,angle=-90]{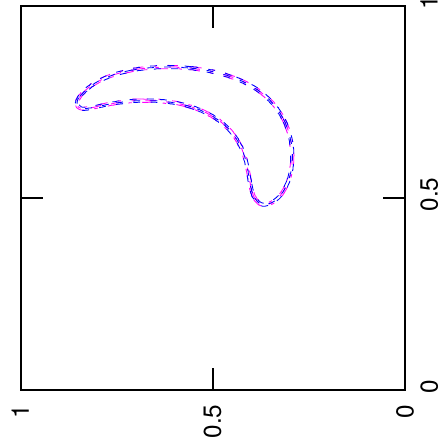}
 \end{minipage}
 \begin{minipage}{0.25\textwidth}
  \centering
   \includegraphics[width=1.\textwidth,height=1.\textwidth,angle=-90]{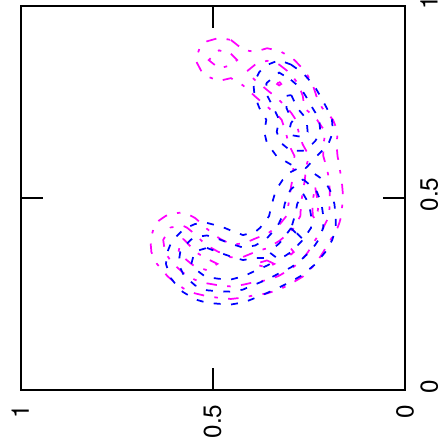}
 \end{minipage}%
 \begin{minipage}{0.25\textwidth}
  \centering
   \includegraphics[width=1.\textwidth,height=1.\textwidth,angle=-90]{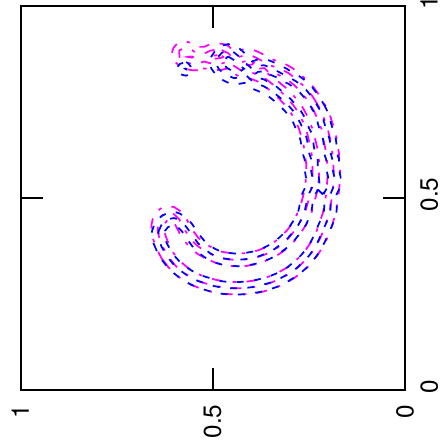}
 \end{minipage}%
 \begin{minipage}{0.25\textwidth}
  \centering
   \includegraphics[width=1.\textwidth,height=1.\textwidth,angle=-90]{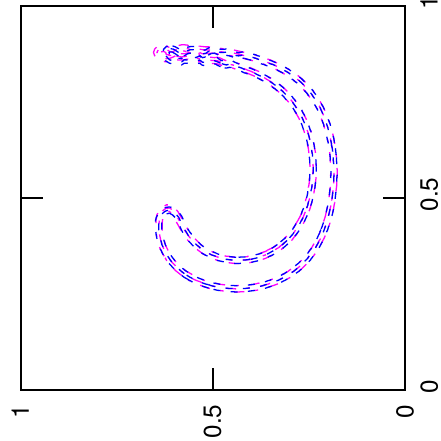}
 \end{minipage}%
 \begin{minipage}{0.25\textwidth}
  \centering
   \includegraphics[width=1.\textwidth,height=1.\textwidth,angle=-90]{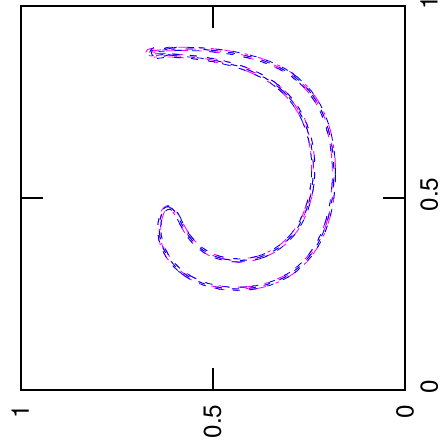}
 \end{minipage}
 \begin{minipage}{0.25\textwidth}
  \centering
   \includegraphics[width=1.\textwidth,height=1.\textwidth,angle=-90]{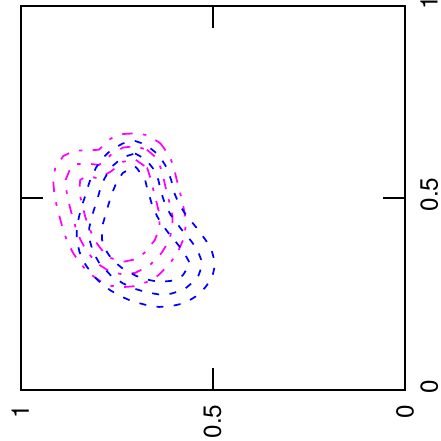}
 \end{minipage}%
 \begin{minipage}{0.25\textwidth}
  \centering
   \includegraphics[width=1.\textwidth,height=1.\textwidth,angle=-90]{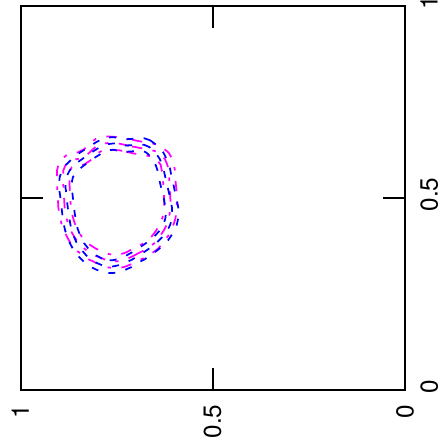}
 \end{minipage}%
 \begin{minipage}{0.25\textwidth}
  \centering
   \includegraphics[width=1.\textwidth,height=1.\textwidth,angle=-90]{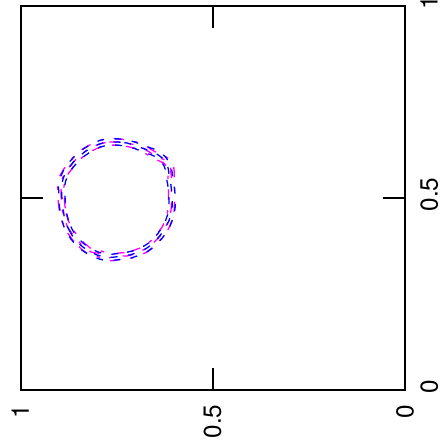}
 \end{minipage}%
 \begin{minipage}{0.25\textwidth}
  \centering
   \includegraphics[width=1.\textwidth,height=1.\textwidth,angle=-90]{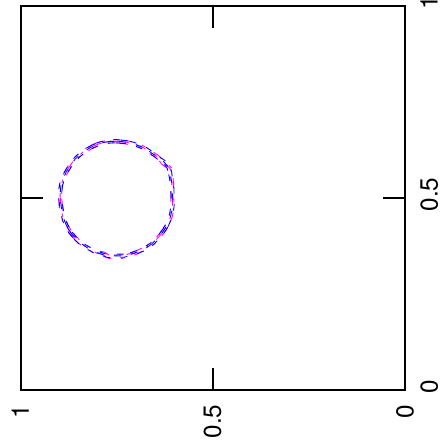}
 \end{minipage}
 \end{centering}
 \caption{\small{Vortex test on four different grids $m_i$, $i=1,\ldots,4$ from 
                 left to right at times $t=0$, $t=0.5$, $t=1.0$, $t=2.0$ from 
                 top to bottom. Figure depicts $\psio = (r_i^1,0,r_i^2)$  
                 iso-lines (see \Eqs{eq32}{eq33}) with $N_\tau=2$ 
                 (magenta, dashed-dotted line)
                 and $N_\tau=16$ (dark blue, dashed line) 
                 re-initialization steps.}}             
 \label{fig31}
 \end{figure}

The iso-lines of $\psio$ 
obtained on grids $m_i$, $i=1,\ldots,4$ with 
$N_\tau=2$ or $N_\tau=16$ are presented 
in figure \ref{fig31}.
For the sake of brevity, 
iso-lines of $\alpps$ are omitted 
as they are almost identical
with iso-lines of $\psio$.
The most important observation in
\Fig{fig31} is that on finer grids, 
the impact of the number 
of re-initialization steps $N_\tau$
on the interface shape 
becomes negligible. 
This result shows convergence
of the new re-initialization method
and is in agreement with 
studies presented in \Fig{fig28}.

When the present results
are compared with 
the results presented
in \cite{olsson05}, 
one notices  due
to the smaller interface width in our
simulations $\eph=\sqrt{2} \Delta x/4 < \Delta x /2$
our results are qualitatively
similar to the results from \cite{olsson05}
on a grid twice as fine.
This illustrates the interface
width $\eph$ is a very important parameter
in the present method; 
$\eph$ decides
not only whether the interface curvature
may be calculated (see Section \ref{sec5.2}) 
but also  governs
the topological changes 
of the interface 
when the grid resolution 
is not sufficient 
to resolve it, see \Fig{fig31}.
 \begin{figure}[ht!] \nonumber
 \begin{centering}
 \begin{minipage}{.33\textwidth}
  \centering
  \includegraphics[width=1.\textwidth,height=1.\textwidth, angle=-90]{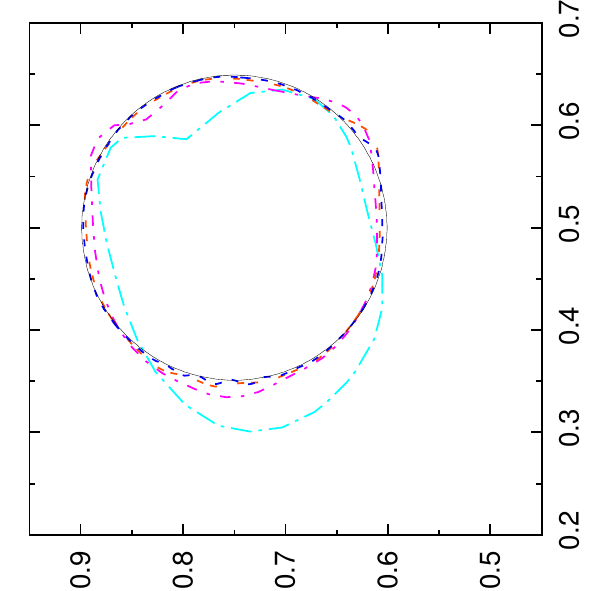}
 \end{minipage}%
 \begin{minipage}{.33\textwidth}
  \centering
  \includegraphics[width=1.\textwidth,height=1.\textwidth, angle=-90]{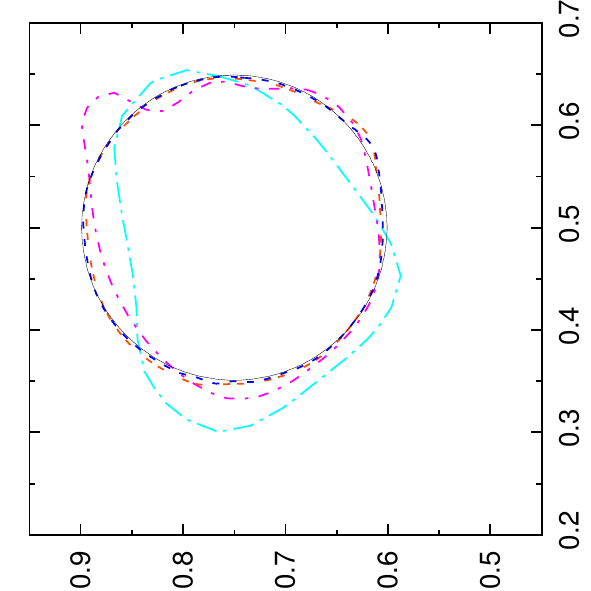}
 \end{minipage}%
 \begin{minipage}{.33\textwidth}
  \centering
  \includegraphics[width=1.\textwidth,height=1.\textwidth, angle=-90]{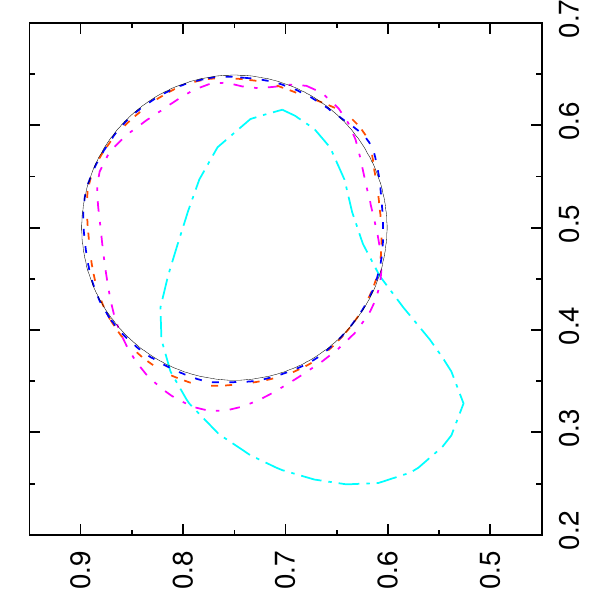}
 \end{minipage}
 \end{centering}
 \caption{\small{Convergence of $\psio=0$ iso-lines on grids $m_i$, $i=1,\ldots,4$
                 towards initial condition (black solid line) in the vortex
                 test case with (a) $N_\tau=2$, (b) $N_\tau=8$, (c) $N_\tau=16$ 
                 re-initialization steps.}}
 \label{fig32}
 \end{figure}
 \begin{figure}[ht!] \nonumber
 \includegraphics[angle=-90]{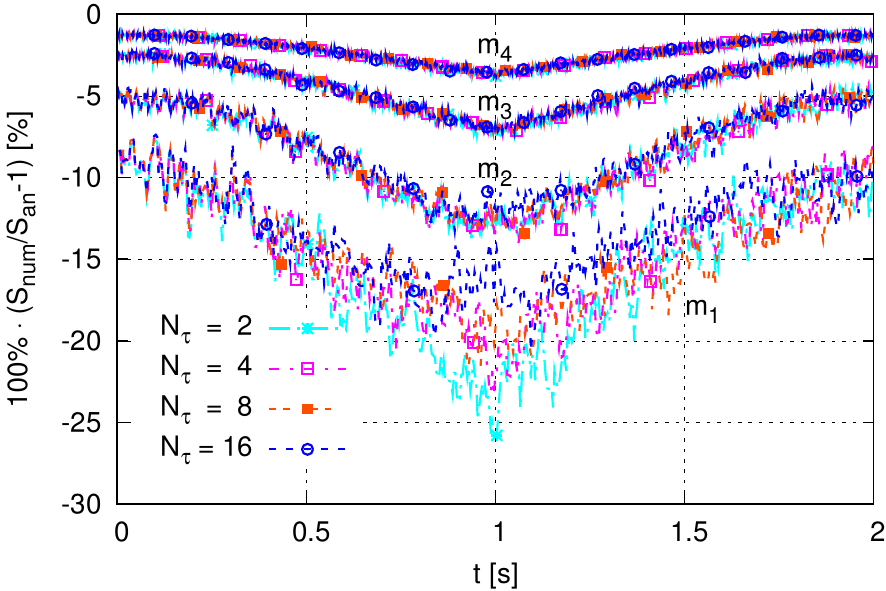} 
 \includegraphics[angle=-90]{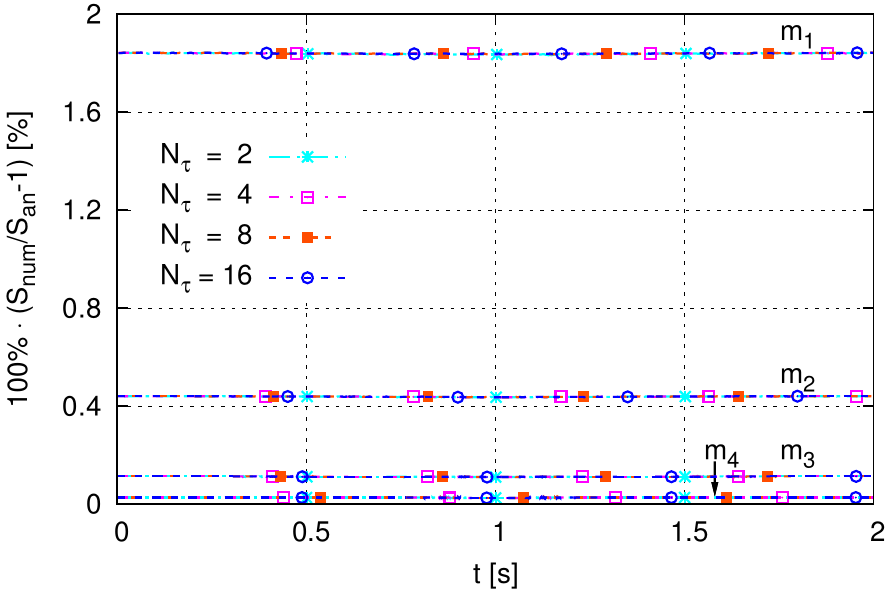} 
 \caption{\small{The area conservation 
                 calculated in the regions: 
                 (a) $r_1 = \lc x_i \,\big|\, 1-\alpps \ge 0.5  \rc$,
                 (b) $r_2 = \lc x_i \,\big|\, \psio \le 8\eph \rc$
                 during vortex test with $N_\tau \le 16$.
                 Error $S_{num}$ defined in \Eq{eq44}
                 is normalized with $S_{an}=\pi R^2$,
                 the number and position of symbols
                 is arbitrary.}}
 \label{fig33}
 \end{figure}
\begin{table}[ht]
\centering 
\begin{tabular}{@{}  l  c c c  c c c @{}} 
\cline{1-7}
 $\Gamma$ & \multicolumn{3}{c}{$\alpps$} & \multicolumn{3}{c}{$\psio$} \\ [0.5ex] 
\cline{1-7} 
 $N_\tau$ &   2      & 8          & 16        & 2           & 8         &  16  \\ [0.5ex] 
\cline{1-7} 
$m_1$    & 2.1424e-1 & 2.3195e-1 & 6.8181e-1  &  2.1434e-1 & 2.3175e-1 &  6.8196e-1  \\ 
$m_2$    & 7.3558e-3 & 1.0785e-2 & 1.1301e-2  &  7.3921e-3 & 1.0773e-2 &  1.1305e-2  \\ 
$m_3$    & 3.2984e-3 & 3.4034e-3 & 3.5041e-3  &  3.3082e-3 & 3.3917e-3 &  3.4896e-3  \\
$m_4$    & 1.6404e-3 & 1.6512e-3 & 1.6444e-3  &   1.6467e-3 & 1.6542e-3 &  1.6433e-3  \\ [0.5ex] 
\cline{1-7}  
\end{tabular}
\caption{Convergence of $L_{1,i}^r$ norm,
         computed at the end of the vortex test
         with the $N_\tau \le 16$ re-initialization steps.
         The circular interface is represented by $\Gamma \lr \alpps \rr$
         and $\Gamma \lr \psio \rr$
         iso-lines, compare with results in figure (\ref{fig32}).}
\label{tab12} 
\end{table}

After rotation of the interface in
the opposite direction  
with the same number of time steps, 
the circular interface $\Gamma$
is reconstructed with 
the first order accuracy,
(see results in figure \ref{fig32} and 
in table \ref{tab12}). 
Similar to the case without interface
deformation, the errors in \Tab{tab12}
do not change much with $N_\tau$.
Convergence towards the initial condition 
may be observed on the gradually refined
grids $m_i$, $i=1,\ldots,4$.

The $L_{1,i}^r$ norm defined in \Eq{eq43b} is not 
the best measure of interface departure 
from its original shape as it does not 
detect oscillations in
the reconstructed interface.
Closer investigations of the results 
obtained on the grid $m_4$ 
(dark blue dashed line in \Fig{fig32})
reveal that the regularity of the 
final interface shape at $t=2\,s$ 
grows with the $N_\tau$ number.
Since re-initialization  improves
$\psio$, the errors of $\bng$
remain smaller when $N_\tau$ is larger
and for this reason oscillations
does not appear in the 
interface shape reconstructed 
at the end of simulations
(compare results on the grid $m_4$ 
in figures \ref{fig32}(a)
and \ref{fig32}(b)(c)).

In figure \ref{fig33}, 
the area conservation
in the vortex test case 
is presented.
Like in rotation of a rigid body, 
we also calculate (integrate) 
errors in two regions 
$r_1$ and $r_2$
in the present case.
Since in this case
the circular interface 
is strongly deformed
through the vortex velocity field,
large variation of the 
error integrated in the region 
$r_1$ is observed, 
see \Fig{fig33}(a).
As the area of 
the interface $\Gamma$
increases when it is deformed, 
the value of the error becomes
smaller and then returns
to its initial level
at the end of the 
simulation.
This result confirms
the present 
re-initialization
method is conservative.
Similar variation may  also be observed
in the case of the error obtained in 
the region $r_2$, 
however the effect 
is very small, and therefore
cannot be observed 
in \Fig{fig33}(b).
The distributions of the errors 
in \Fig{fig33}(a)(b) confirm that at least 
the first-order  convergence rate of
the area (mass) is achieved by 
the present numerical method.
Moreover, we note that
the number of re-initialization 
steps $N_\tau$ performed on each 
time step $\Delta t$ 
has a small impact on 
the calculated area 
when the mesh is sufficiently 
fine
(compare results with 
different $N_\tau$ in figures 
\ref{fig31} and \ref{fig33}).
 \begin{figure}[ht!] \nonumber
 \begin{centering}
 \begin{minipage}{0.25\textwidth}
  \centering
   \includegraphics[width=1.\textwidth,height=1.\textwidth,angle=-90]{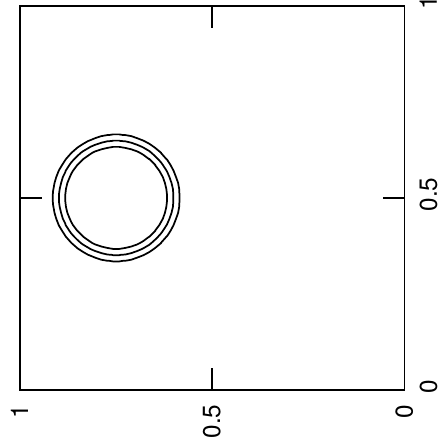}
 \end{minipage}%
 \begin{minipage}{0.25\textwidth}
  \centering
   \includegraphics[width=1.\textwidth,height=1.\textwidth,angle=-90]{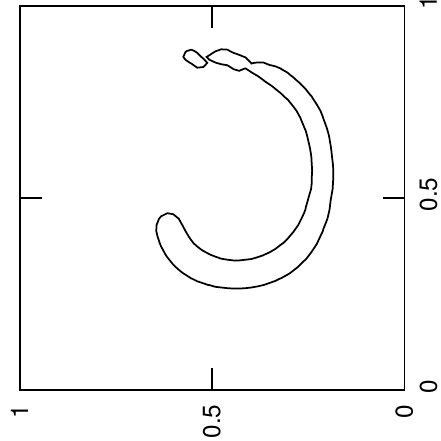}
 \end{minipage}%
 \begin{minipage}{0.25\textwidth}
  \centering
   \includegraphics[width=1.\textwidth,height=1.\textwidth,angle=-90]{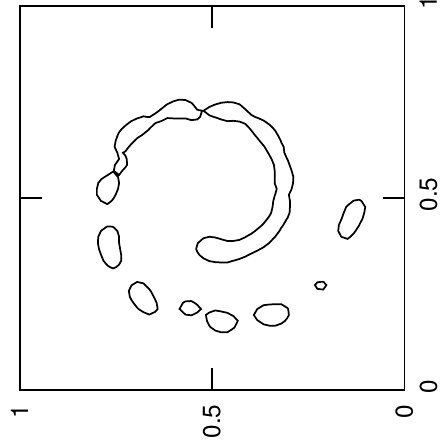}
 \end{minipage}%
 \begin{minipage}{0.25\textwidth}
  \centering
   \includegraphics[width=1.\textwidth,height=1.\textwidth,angle=-90]{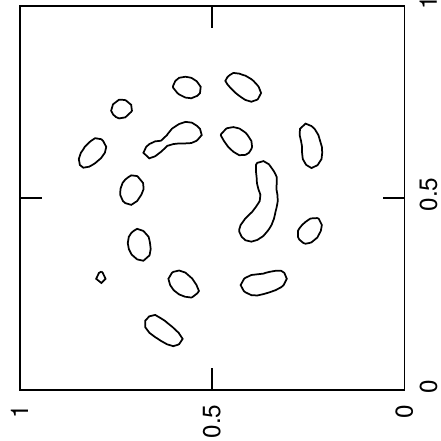}
 \end{minipage}
 \begin{minipage}{0.25\textwidth}
  \centering
   \includegraphics[width=1.\textwidth,height=1.\textwidth,angle=-90]{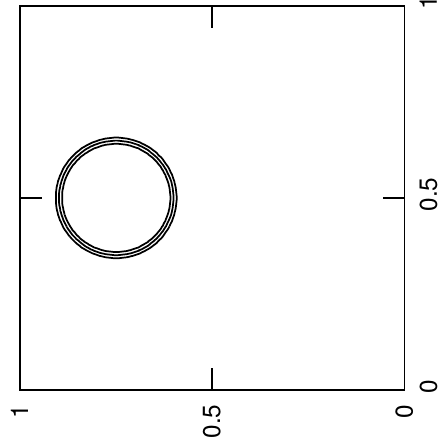}
 \end{minipage}%
 \begin{minipage}{0.25\textwidth}
  \centering
   \includegraphics[width=1.\textwidth,height=1.\textwidth,angle=-90]{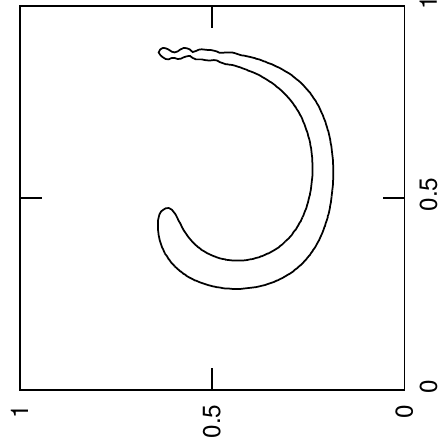}
 \end{minipage}%
 \begin{minipage}{0.25\textwidth}
  \centering
   \includegraphics[width=1.\textwidth,height=1.\textwidth,angle=-90]{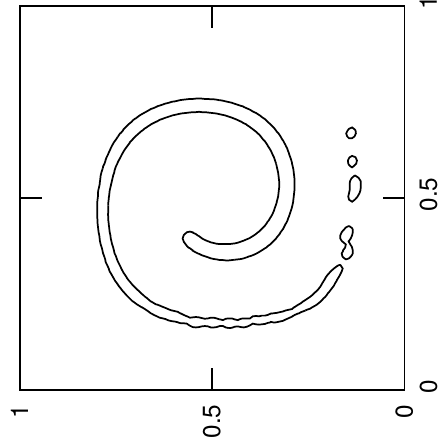}
 \end{minipage}%
 \begin{minipage}{0.25\textwidth}
  \centering
   \includegraphics[width=1.\textwidth,height=1.\textwidth,angle=-90]{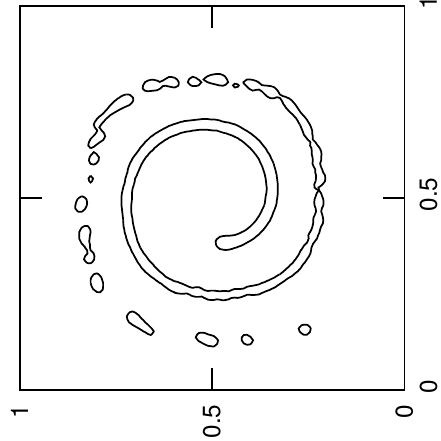}
 \end{minipage}
 \begin{minipage}{0.25\textwidth}
  \centering
   \includegraphics[width=1.\textwidth,height=1.\textwidth,angle=-90]{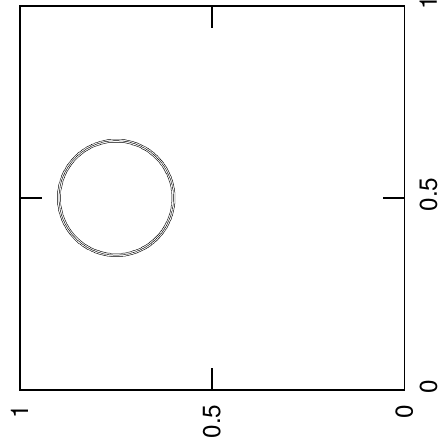}
 \end{minipage}%
 \begin{minipage}{0.25\textwidth}
  \centering
   \includegraphics[width=1.\textwidth,height=1.\textwidth,angle=-90]{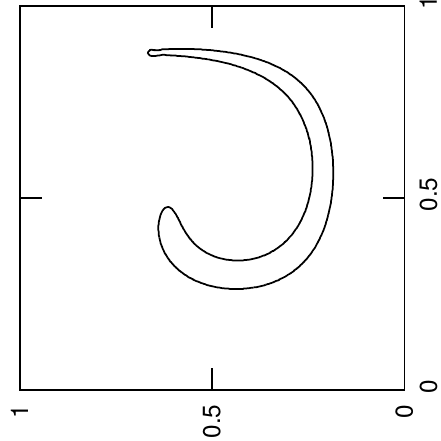}
 \end{minipage}%
 \begin{minipage}{0.25\textwidth}
  \centering
   \includegraphics[width=1.\textwidth,height=1.\textwidth,angle=-90]{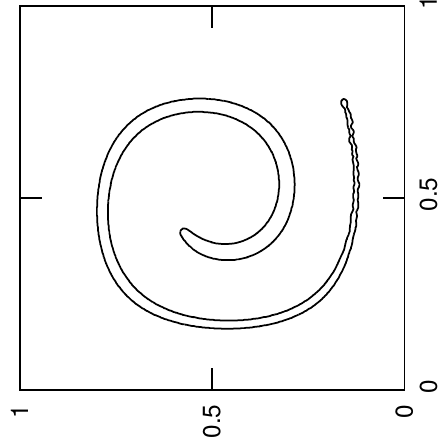}
 \end{minipage}%
 \begin{minipage}{0.25\textwidth}
  \centering
   \includegraphics[width=1.\textwidth,height=1.\textwidth,angle=-90]{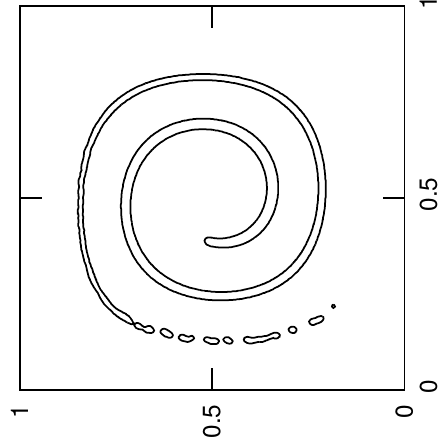}
 \end{minipage}
 \begin{minipage}{0.25\textwidth}
  \centering
   \includegraphics[width=1.\textwidth,height=1.\textwidth,angle=-90]{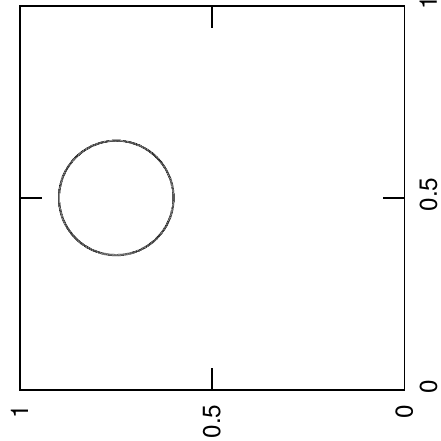}
 \end{minipage}%
 \begin{minipage}{0.25\textwidth}
  \centering
   \includegraphics[width=1.\textwidth,height=1.\textwidth,angle=-90]{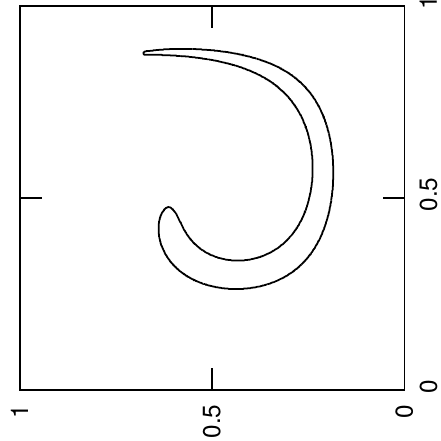}
 \end{minipage}%
 \begin{minipage}{0.25\textwidth}
  \centering
   \includegraphics[width=1.\textwidth,height=1.\textwidth,angle=-90]{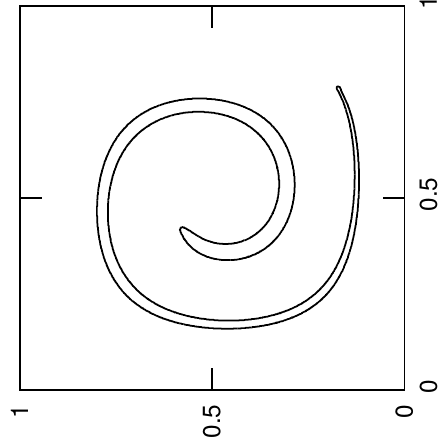}
 \end{minipage}%
 \begin{minipage}{0.25\textwidth}
  \centering
   \includegraphics[width=1.\textwidth,height=1.\textwidth,angle=-90]{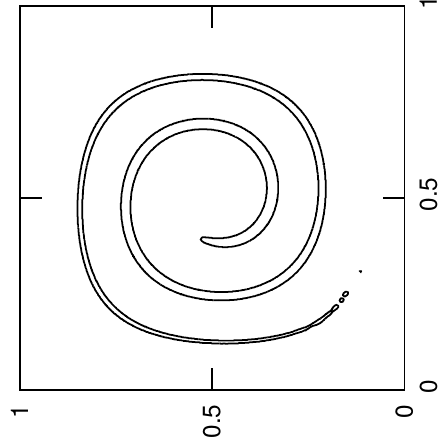}
 \end{minipage}
 \end{centering}
 \caption{\small{$\psio=0$ iso-lines, in the vortex test 
                 on  grids $m_i$, $i=2,\ldots,5$ 
                 (from top to bottom) 
                 after
                 $T=0$, 
                 $T=1$, 
                 $T=2$, 
                 $T=3$ revolutions (from left to right),
                 the width of the interface $\eph=\sqrt{2} \Delta x/4$, 
                 $N_\tau=4$ on each time step $\Delta t$, $\Delta \tau=\eph$.}}             
 \label{fig34}
 \end{figure}
 \begin{figure}[ht!] \nonumber
 \includegraphics[angle=-90]{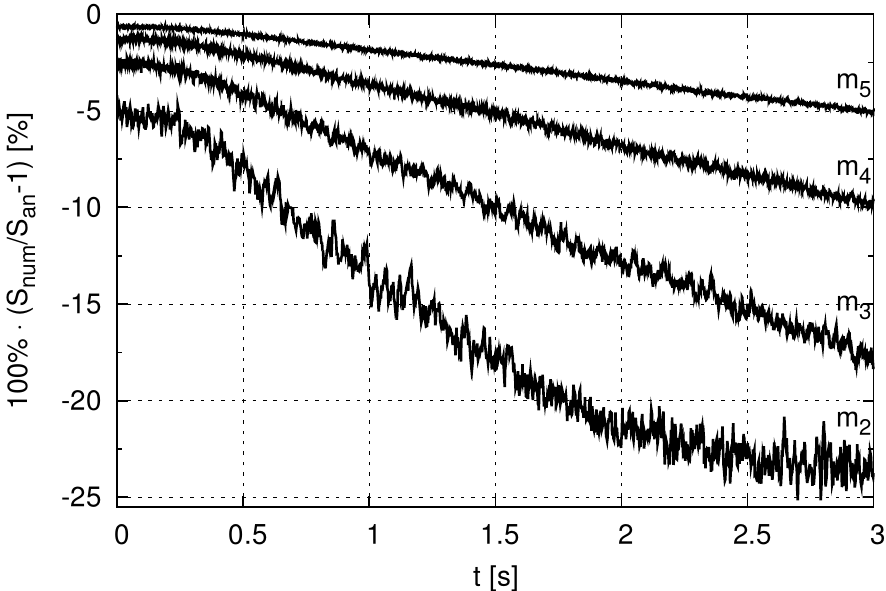}
 \includegraphics[angle=-90]{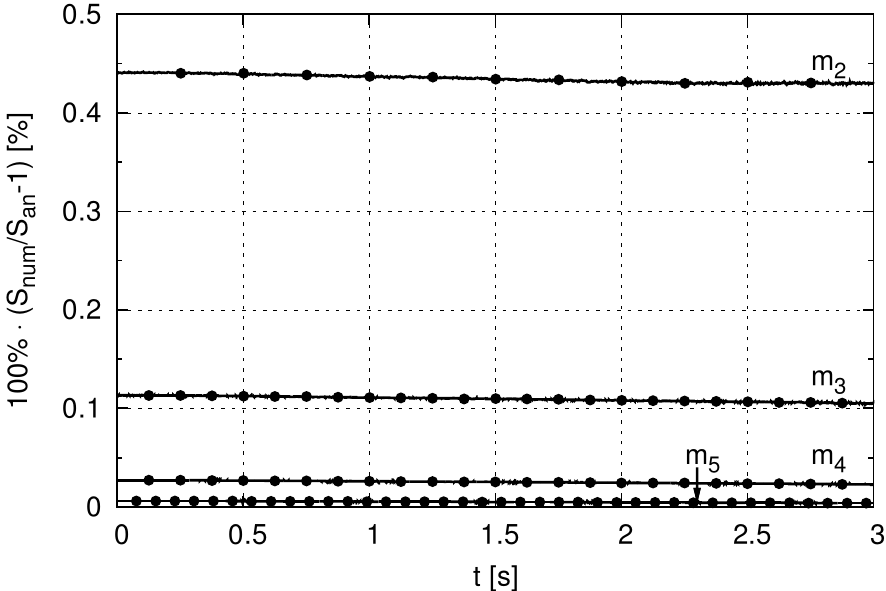}
 \caption{\small{The errors in the area conservation 
                 calculated in the regions 
                 (a) $r_1 = \lc x_i \,\big|\, 1-\alpps \ge 0.5  \rc$,
                 (b) $r_2 = \lc x_i \,\big|\, \psio \le 8\eph \rc$
                 during vortex test with $N_\tau = 4$ on each 
                 time step $\Delta t$.
                 Error $S_{num}$ defined in \Eq{eq44}
                 is normalized with $S_{an}=\pi R^2$.}}
 \label{fig35}
 \end{figure}

To investigate behavior of the new 
numerical method during longer 
integration times  (here $t=3\,s$),
we perform the vortex test without
reversing the velocity field after 
the first revolution.
In this test case, we stretch
the interface until it is broken into
bubbles due to insufficient number 
of grid points required for its
reconstruction.
As our advection and re-initialization 
methods are conservative, 
and the re-initialization method  
keeps the prescribed thickness 
of the interface $\eph$ constant,
this is the only 
possible scenario
which may occur 
when the grid resolution
is insufficient to resolve 
the interface $\Gamma$.
We emphasize that there is no
physical mechanism behind the
abovementioned  break up of the
interface.
This feature of the present
method is direct consequence
of its conservative properties.

Using experience gathered 
in the previous examples, 
we use $N_\tau=4$ re-initialization steps
in each single time step $\Delta t$,
$\eph=\sqrt{2}\Delta x/4$ and $\Delta \tau=\eph$.
The results from this study, obtained
on four grids $m_i$ where $i=2,\ldots,5$,
at four different time moments $T$,
are presented in \Fig{fig34}.
In order to visualize the interface
thickness on each grid $m_i$, at $T=0$ three
iso-lines of  $\psio = (r_i^1,0,r_i^2)$  
(see \Eqs{eq32}{eq33}) are depicted.

We note that the longer 
the integration time is, 
the thiner the interface $\Gamma$
filaments become.
However, 
if the mesh resolution 
is insufficient to resolve
the interface $\Gamma$, 
then the interface break 
up will occur as this is the only
way to preserve the area (mass) 
correctly.

Errors in 
the area conservation 
are depicted 
in figure \ref{fig35}, where, 
as in the previous examples, 
$S_{num}$ is computed
by integration in 
two regions
$r_1$ and $r_2$.
The results in \Fig{fig35} show
that the area is conserved in 
spite of large deformation of
the interface $\Gamma$ on all
grids $m_i$, $i=2,\ldots,5$.
As discussed previously,
the error calculated in 
the region $r_1$
permits tracing  
the impact 
of the interface deformation 
on the area (mass) conservation.
We note that after break up of 
the interface into separate bubbles on 
grid $m_2$ at 
the time about $T=2.5\,s$, 
values of $S_{num}$ integrated 
in region $r_1$
cease to drop 
as the interface 
deformation is stopped
(see \Fig{fig35}(a)).
It is expected 
after large
enough integration 
times $t \gg 3\,s$,
at each grid,
the interface 
will finally 
disintegrate
into droplets 
that can can 
be transported
in the given 
velocity field.

The total error in the
area (mass) conservation, 
obtained by integration of $\alpps$ 
in the region $r_2$,
is almost constant on all meshes used 
in the present study.
The loss of total area (mass) 
is present but negligibly small
as shown in figure \ref{fig35}(b).

\section{Conclusions}
\label{sec7}

In this work,
a new re-initialization
method of the conservative 
level-set function was introduced
and verified.
We have shown 
that the re-initialization
and advection equations 
of the conservative level-set function
$\alpps$ are mathematically equivalent
to the re-initialization and advection equations
of the localized signed distance function $\psio$
(see equations (\ref{eq29}) and (\ref{eq30a})).
It was also proven that 
the RHS of re-initialization equation (\ref{eq29})
is equal to zero when
$|\nabla \psi_0| =1$,
and $\alpps = H \lr \psi_0 \rr$.
These two solutions
to the re-initialization equation (\ref{eq29})
permit computing the interface curvature, 
and controlling its thickness  $\eph < \Delta x/2$ 
which assures the improvement 
of the spatial resolution 
of the new method,
see discussion 
in \Sec{sec5.2}
and the results in \Sec{sec6}.

In the present re-initialization procedure,
the conservative level-set function $\alpps$ 
remains continuous and bounded 
inside the support 
of $\delta \lr \alpha \rr/\eph =\alpha \lr 1-\alpha \rr/\eph$ 
function.
For this reason, 
we do not need to use
the fast marching method,
required in other 
re-initialization techniques
reported in the literature,  
in order to extend
the signed distance function
$\psio$ away from the interface
$\Gamma$.
The continuity
of the solution 
assured
by the new re-initialization method
avoids artificial interface deformations,
which are eventually introduced by 
the flux limiters typically used 
to bound counteracting
compressive and diffusive
fluxes on the RHS of \Eq{eq1}.
In the light of 
the other results presented 
in the literature,
the new re-initialization method 
shows fast convergence 
and in the 1D cases,
reconstructs exact 
behavior of the analytical solution 
to the partial 
differential equation 
(\ref{eq1}), see \Fig{fig4}.

Discretization
of the first and second-order  
spatial derivatives of $\alpps$,
which is consistent 
with the discretization 
of the re-initialization equation
and the theory of distributions,
achieves 
the second-order
convergence rate 
of the interface 
curvature $\kappa$,
see \Eq{eq31k}.
Such level of accuracy 
is obtained in the second-order 
accurate finite volume solver
without  geometrical 
reconstruction of the interface 
or an introduction of the higher-order
essentially non-oscillatory
interpolation schemes in 
the discretization procedure.

The advection tests performed in
Section \ref{sec6.2} show  the
new method conserves the total area 
(mass) with almost
second order accuracy, 
and the shape 
of the advected interface 
is independent from 
the number of the
re-initialization steps
if the spatial 
and temporal resolutions
are sufficient.

The new re-initialization 
method presented herein
may be considered
a backbone of 
the conservative level-set method,
a potentially 
good replacement  
for the compressive 
interface capturing schemes 
commonly used in the 
fast multiphase flow 
solvers.
Its implementation 
is simple and its strengths 
include accurate reconstruction 
of the interface curvature 
and conservation
of the area  (mass) 
without smearing 
the interface.

\section*{Acknowledgments}
The work was funded by the German Research 
Foundation (DFG) in the framework of the project 
"Modeling of turbulence-interface interactions in 
two-fluid systems" WA: 3098/2-1, AOBJ: 595642.


\appendix

\section{Limits of the mapping function}
\label{appA}

In this appendix,
the limits of the
level-set function
$\alpp$ and the mapping
function $\psig$ 
given by \Eq{eq2}
and \Eq{eq8}
are calculated. 
Since in this section we consider analytical 
properties of $\psig$ next
$\epsilon=0$ in equation (\ref{eq8}).
We assume that the number
of grid nodes $N_c \to \infty$ 
and $0 < \gamma \ll \Delta x$
is a small constant. 
The 
analytical solution given 
by \Eq{eq2} may be 
rewritten as 
\blnm
\be
 \alpha =   \frac{1}{2} \ls 1+\tanh{\lr \frac{\psi_0}{2\eph} \rr} \rs
        =    \frac{1}{1+\exp{\lr-2N_c \psi_0  \rr}}
\label{Beq1}
\ee
\elnm
%
%
where $N_c=1/\Delta x$ 
and $\psi_0$ is the signed 
distance function given by \Eq{eq3}.
Additionally, we rearrange 
\Eq{eq13} into 
\blnm
\be
 \psi_1 =  \frac{\exp{\lr N_c \psi_0 \gamma \rr}}{\exp{\lr N_c \psi_0 \gamma \rr}+\exp{\lr-N_c \psi_0 \gamma \rr}} 
        =  \frac{1}{1+\exp{\lr-2N_c \psi_0 \gamma \rr}}.
\label{Beq3}
\ee
\elnm
First, we note that 
if $\gamma \to 0$ in \Eq{Beq3}, $\psi_1 \to 1/2$
as it was observed in \cite{shukla10}.
Let us now consider three 
limits of \Eq{Beq1} and \Eq{Beq3} 
when $N_c \gamma \to \infty$ and
$0 < \gamma \ll \Delta x$:
\begin{enumerate}
 \item{ when $\psi_0<0\quad :\quad \psi_1 \to 0,\, \alpha \to 0,$}
 \item{ when $\psi_0=0\quad :\quad \psi_1=1/2,\, \alpha=1/2,$ } 
 \item{ when $\psi_0>0\quad :\quad \psi_1 \to 1,\, \alpha \to 1 $}
\end{enumerate}
Hence, when the number of grid 
points $N_c  \to \infty$ and $\gamma=const$ 
functions $\psi_1$ and $\alpha$   
become the Heaviside 
function $H \lr \psi_0 \rr$.
The phase indicator
function
$H(\psi_0)$ is
typically discretized in the 
volume of fluid (VOF) methods \cite{trygg11}.

\section{Discretization of the re-initialization equation}
\label{appB}

In this appendix, 
the discretization
of \Eq{eq1} or equivalent
\Eqs{eq29}{eq30} 
in the framework 
of the second order 
accurate finite volume 
method is presented.
Since all terms on the RHS 
of \Eq{eq30} are in 
the form of the mass forces, 
after its integration 
in the control volume $V_P$
one would obtain
\blnm
\begin{align}
  \bs
 \pd{\alpha}{\tau} \Big|_P 
  = \ls \bng \cdot \nabla  \delta \lr \alpha \rr  \lr |\nabla \phi | -  1 \rr  \rs |_P \\
  + \ls \bng \cdot \nabla   \lr |\nabla \phi | -  1 \rr  \delta \lr \alpha \rr  \rs |_P \\
  - \ls \kappa \lr |\nabla \phi | -  1 \rr   \delta \lr \alpha \rr   \rs |_P,
  \es
\label{Aeq3}
\end{align}
\elnm
where $\bng = \nabla \phi/|\nabla \phi|$ 
and $\phi$ is equal to $\alpx$
or is one of the mapping functions
$\psig$, $\psio$.
\Eq{Aeq3} is the 
discretization to
the non-conservative form
of \Eq{eq1},
all terms with sub-script $P$
can be obtained from the values stored in  
the centers of the control volumes. 

In the present work, we use
the conservative discretization
to \Eq{eq1} or equivalent \Eq{eq29}. 
After integration of \Eq{eq29}
in the control volume $V_P$, 
employment
of the Gauss 
theorem and mid-point
rule at the centers
of faces $f$ 
and in the 
centers of
control
volume  
$P$ one obtains
\blnm
\be
  \pd{\alpha}{\tau} \Big|_P = \frac{1}{V_P} 
   \sum_{f=1}^{n_b} \ls \delta \lr \alpha \rr \lr |\nabla \phi | -  1 \rr \bng \cdot \bn \rs_f S_f,
\label{Aeq4}
\ee
\elnm
where $n_b$ 
is  the number 
of neighbors 
of the  control volume $P$,
$\ls \bng\cdot\bn \rs_f$ is
dot product of the normal 
$\bng=\nabla \phi/|\nabla \phi|$ interpolated
to the face $f$ and normal 
$\bn_f$ at the face $f$,
$\deltaa$
is defined by \Eq{eq24};
$|\nabla \phi|$
is computed using
the second-order 
central-difference 
approximation
to the 
$\phi$ 
gradient
components, 
at
face $f=e$ 
this approximation
reads
\blnm
\begin{align}
  \bs
  \pd{\phi}{x_1} \Big|_e &\approx  
   \frac{\lr \phi_{E}-\phi_{P} \rr}{ \Delta x}, \\
  \pd{\phi}{x_2} \Big|_e &\approx
   \frac{\lr \phi_{N}+\phi_{NE}-\phi_{S}-\phi_{SE} \rr}{4 \Delta y}, \\
  \pd{\phi}{x_3} \Big|_e &\approx
  \frac{\lr \phi_{T}+\phi_{TE}-\phi_{B}-\phi_{BE} \rr}{4 \Delta z},
  \es
\label{Aeq5}
\end{align}
\elnm
where subscript $E,N,T,\ldots$ 
represent the centers 
of the neighbor 
control volumes.
The interpolation 
to the faces $f$
of the control volume $P$ is performed
using  second-order accurate
linear interpolation scheme, 
which on  uniform grids 
simplifies to 
$\phi_f=1/2\lr \phi_F+\phi_P\rr$, 
see \cite{peric02,schaefer06}.
Dependent on the test case, 
$\phi$ in \Eq{Aeq5} is calculated
using $\alpx$ or $\psi$ 
obtained from the mapping functions
given by \Eq{eq3} or \Eq{eq8}
with different $\gamma$
values.

\bibliography{mybibfile}

\begin{thebibliography}{10}
\expandafter\ifx\csname url\endcsname\relax
  \def\url#1{\texttt{#1}}\fi
\expandafter\ifx\csname urlprefix\endcsname\relax\def\urlprefix{URL }\fi
\expandafter\ifx\csname href\endcsname\relax
  \def\href#1#2{#2} \def\path#1{#1}\fi

\bibitem{olsson05}
E.~Olsson, G.~Kreiss, A conservative level-set method for two phase flow, J.
  Comput. Phys. 210 (2005) 225--246.

\bibitem{shukla10}
R.~Shukla, C.~Pantano, J.~Freund, An interface capturing method for the
  simulation of mutli-phase compressible flows, J. Compt. Phys. 229 (2010)
  7411--7439.

\bibitem{so11}
K.~So, X.~Hu, N.~Adams, Anti-diffusion method for interface steepening in
  two-phase incompressible flows, J. Compt. Phys. 230 (2011) 5155--5177.

\bibitem{tiwari13}
A.~Tiwari, J.~Freund, C.~Pantano, A diffuse interface model with immiscibility
  preservation, J. Compt. Phys. 252 (2013) 290--309.

\bibitem{balcazar2014}
N.~Balcázar, L.~Jofre, O.~Lehmkuhl, J.~Castro, J.~Rigola, A
  finite-volume/level-set method for simulating two-phase flows on unstructured
  grids, International Journal of Multiphase Flow 64~(0) (2014) 55 -- 72.
\newblock \href
  {http://dx.doi.org/http://dx.doi.org/10.1016/j.ijmultiphaseflow.2014.04.008}
  {\path{doi:http://dx.doi.org/10.1016/j.ijmultiphaseflow.2014.04.008}}.

\bibitem{glasner2001}
K.~Glasner, Nonlinear preconditioning for diffuse interfaces, Journal of
  Computational Physics 174~(2) (2001) 695 -- 711.
\newblock \href {http://dx.doi.org/http://dx.doi.org/10.1006/jcph.2001.6933}
  {\path{doi:http://dx.doi.org/10.1006/jcph.2001.6933}}.

\bibitem{sun2007}
Y.~Sun, C.~Beckermann, Sharp interface tracking using the phase-field equation,
  Journal of Computational Physics 220~(2) (2007) 626 -- 653.
\newblock \href {http://dx.doi.org/http://dx.doi.org/10.1016/j.jcp.2006.05.025}
  {\path{doi:http://dx.doi.org/10.1016/j.jcp.2006.05.025}}.

\bibitem{trygg11}
G.~Tryggvason, R.~Scardovelli, S.~Zaleski, Direct Numerical Simulations of
  Gas-Liquid Multiphase Flows, Cambridge University Press, 2011.

\bibitem{osher1988}
S.~Osher, J.~A. Sethian, Fronts propagating with curvature-dependent speed:
  {A}lgorithms based on {H}amilton-{J}acobi formulations, Journal of
  Computational Physics 79~(1) (1988) 12 -- 49.
\newblock \href
  {http://dx.doi.org/http://dx.doi.org/10.1016/0021-9991(88)90002-2}
  {\path{doi:http://dx.doi.org/10.1016/0021-9991(88)90002-2}}.

\bibitem{sussman94}
M.~Sussman, P.~Smereka, S.~J. Osher, A level set approach for computing
  solutions to incompressible two-phase flows, J. Comput. Phys 114 (1994)
  146--159.

\bibitem{sussman98}
M.~Sussman, E.~Fatemi, P.~Smereka, S.~Osher, An improved level set method for
  incompressible two-phase flows, Computers and Fluids 27~(5–6) (1998) 663 --
  680.
\newblock \href
  {http://dx.doi.org/http://dx.doi.org/10.1016/S0045-7930(97)00053-4}
  {\path{doi:http://dx.doi.org/10.1016/S0045-7930(97)00053-4}}.

\bibitem{osher03}
S.~Osher, R.~Fedkiw, {L}evel {S}et {M}ethods and {D}ynamic {I}mplicit
  {S}urfaces, Springer Verlag, INC. New-York, 2003.

\bibitem{hartmann2010}
D.~Hartmann, M.~Meinke, W.~Schr{\"o}der, The constrained reinitialization
  equation for level set methods, Journal of Computational Physics 229~(5)
  (2010) 1514 -- 1535.
\newblock \href {http://dx.doi.org/http://dx.doi.org/10.1016/j.jcp.2009.10.042}
  {\path{doi:http://dx.doi.org/10.1016/j.jcp.2009.10.042}}.

\bibitem{rocca14}
G.~D. Rocca, G.~Blanquart, Level set reinitialization at a contact line,
  Journal of Computational Physics 265~(0) (2014) 34 -- 49.
\newblock \href {http://dx.doi.org/http://dx.doi.org/10.1016/j.jcp.2014.01.040}
  {\path{doi:http://dx.doi.org/10.1016/j.jcp.2014.01.040}}.

\bibitem{olsson07}
E.~Olsson, G.~Kreiss, S.~Zahedi, A conservative level set method for two phase
  flow {II}, Journal of Computational Physics 225~(1) (2007) 785 -- 807.
\newblock \href {http://dx.doi.org/http://dx.doi.org/10.1016/j.jcp.2006.12.027}
  {\path{doi:http://dx.doi.org/10.1016/j.jcp.2006.12.027}}.

\bibitem{mccaslin2014}
J.~O. McCaslin, O.~Desjardins, A localized re-initialization equation for the
  conservative level set method, Journal of Computational Physics 262~(0)
  (2014) 408 -- 426.
\newblock \href {http://dx.doi.org/http://dx.doi.org/10.1016/j.jcp.2014.01.017}
  {\path{doi:http://dx.doi.org/10.1016/j.jcp.2014.01.017}}.

\bibitem{gottlieb98}
S.~Gottlieb, S.~C.-W., Total variation diminishing {R}unge-{K}utta schemes,
  Mathematics of Computation 67 (1998) 73--85.

\bibitem{bronstein2012}
I.~Bronstein, K.~Semendajew, G.~Musiol, H.~M{\"u}lig, Tachenbuch der
  Mathematik, Harri Deutsch, 2012.

\bibitem{ubbink99}
O.~Ubbink, R.~Issa, Method for capturing sharp fluid interfaces on arbitrary
  meshes, J. Comput. Phys. 153 (1999) 26--50.

\bibitem{waclawczyk06}
T.~Wac{\l}awczyk, T.~Koronowicz, Modelling of the free surface flows with
  high-resolution schemes, Chemical and Process Engineering 27 (2006) 783--802.

\bibitem{waclawczyk08_2}
T.~Wac{\l}awczyk, T.~Koronowicz, Comparison of {CICSAM} and {HRIC} high
  resolution schemes for interface capturing, J. Theoretical and Applied
  Mechanics 46 (2008) 325--345.

\bibitem{waclawczyk08_3}
T.~Wac{\l}awczyk, T.~Koronowicz, Remarks on prediction of wave drag using {VOF}
  method with interface capturing approach, Archives of Civil and Mechanical
  Engineering 8~(1) (2008) 5 -- 14.
\newblock \href
  {http://dx.doi.org/http://dx.doi.org/10.1016/S1644-9665(12)60262-3}
  {\path{doi:http://dx.doi.org/10.1016/S1644-9665(12)60262-3}}.

\bibitem{mwaclawczyk11}
M.~Wac{\l}awczyk, M.~Oberlack, Closure proposals for the tracking of
  turbulence-agitated gas-liquid interfaces in stratified flows, Int. J.
  Multiphase Flow 37 (2011) 967--976.

\bibitem{twaclawczyketal14}
T.~Wac{\l}awczyk, M.~Wac{\l}awczyk, S.~V. Kraheberger, Modelling of
  turbulence-interface interactions in stratified two-phase flows, Journal of
  Physics: Conference Series 530.

\bibitem{waclawczyk2015}
M.~Wac{\l}awczyk, T.~Wac{\l}awczyk, A priori study for the modelling of
  velocity-interface correlations in the stratified air–water flows,
  International Journal of Heat and Fluid Flow 52~(0) (2015) 40 -- 49.
\newblock \href
  {http://dx.doi.org/http://dx.doi.org/10.1016/j.ijheatfluidflow.2014.11.004}
  {\path{doi:http://dx.doi.org/10.1016/j.ijheatfluidflow.2014.11.004}}.

\bibitem{peric02}
J.~H. Ferziger, M.~Peri{\'c}, Computational {M}ethods for {F}luid {D}ynamics,
  Springer Verlag, Berlin Heidelberg New York, 2002.

\bibitem{schaefer06}
M.~Sch{\"a}fer, {C}omputational {E}ngineering, {I}ntroduction to {N}umerical
  {M}ethods, Springer-Verlag Berlin Heidelberg New York, 2006.

\bibitem{xue98}
L.~Xue, Entwicklung eines effizenten parallelen l{\"o}sungsalgorithmus zur
  dreidimenisonalen simulation komplexer turbulenter str{\"o}mungen, Ph.D.
  thesis, Doktorarbeit, Technische Universit{\"a}t Berlin (1998).

\end{thebibliography}

\end{document}